# GROUPS AS GRAPHS

**W. B. Vasantha Kandasamy**
**Florentin Smarandache**

**2009**

# GROUPS AS GRAPHS

**W. B. Vasantha Kandasamy**
e-mail: **vasanthakandasamy@gmail.com**
web: **http://mat.iitm.ac.in/~wbv**
**www.vasantha.in**

**Florentin Smarandache**
e-mail: **smarand@unm.edu**

**2009**

# CONTENTS

# PREFACE

Through this book, for the first time we represent every finite group in the form of a graph. The authors choose to call these graphs as identity graph, since the main role in obtaining the graph is played by the identity element of the group.

This study is innovative because through this description one can immediately look at the graph and say the number of elements in the group G which are self-inversed. Also study of different properties like the subgroups of a group, normal subgroups of a group, p-sylow subgroups of a group and conjugate elements of a group are carried out using the identity graph of the group in this book. Merely for the sake of completeness we have defined similar type of graphs for algebraic structures like commutative semigroups, loops, commutative groupoids and commutative rings.

This book has four chapters. Chapter one is introductory in nature. The reader is expected to have a good background of algebra and graph theory in order to derive maximum understanding of this research.

The second chapter represents groups as graphs. The main feature of this chapter is that it contains 93 examples with diagrams and 18 theorems. In chapter three we describe commutative semigroups, loops, commutative groupoids and commutative rings as special graphs. The final chapter contains 52 problems.

Finally it is an immense pleasure to thank Dr. K. Kandasamy for proof-reading and Kama and Meena without whose help the book would have been impossibility.

W.B.VASANTHA KANDASAMY
FLORENTIN SMARANDACHE

**Chapter One**

# INTRODUCTION TO SOME BASIC CONCEPTS

This chapter has two sections. In section one; we introduce some basic and essential properties about rooted trees. In section two we just recall the definitions of some basic algebraic structures for which we find special identity graphs.

1.1 Properties of Rooted Trees

In this section we give the notion of basic properties of rooted tree.

DEFINITION 1.1.1: *A tree in which one vertex (called) the root is distinguished from all the others is called a rooted tree.*

*Example 1.1.1:*

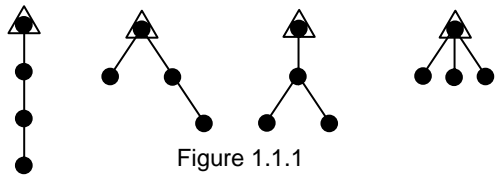

Figure 1.1.1

Figure 1.1.1 gives rooted trees with four vertices.

We would be working with rooted trees of the type.

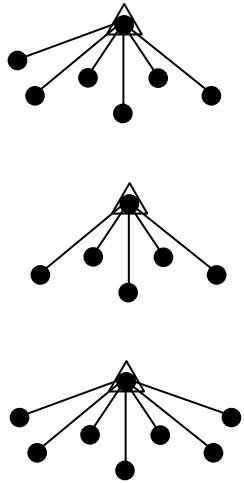

Figure 1.1.2

We will also call a vertex to be the center of the graph if every vertex of the graph has an edge with that vertex; we may have more than one center for a graph.

In case of a complete graph $K_n$ we have n centers.

We have four centers for $K_4$.

$K_3$ has 3 centers

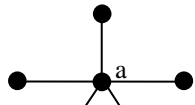

Figure 1.1.3

a is a center of the graph.

For rooted trees the special vertex viz. the root is the center.

Cayley showed that every group of order n can be represented by a strongly connected digraph of n vertices.

However we introduce a special identity graph of a group in the next chapter. As identity plays a unique role in the graph of group we choose to call the graph related with the group as the identity graph of the group G.

For more about Cayley graph and graphs in general refer any standard book on graph theory.

1.2 Basic Concepts

In this section we just recall some basic notions about some algebraic structures to make this book a self contained one.

**DEFINITION 1.2.1:** *A non empty set S on which is defined an associative binary operation * is called a semigroup; if for all a, b ∈ S, a * b ∈ S.*

*Example 1.2.1:* $Z^+$ = {1, 2, …} is a semigroup under multiplication.

*Example 1.2.2:* Let $Z_n$ = {0, 1, …, n – 1} is a semigroup under multiplication modulo n. n ∈ $Z^+$.

*Example 1.2.3:* S(2) = {set of all mappings of (1, 2) to itself is a semigroup under composition of mappings}. The number of elements in S(2) is $2^2 = 4$.

*Example 1.2.4:* S(n) = {set of all mappings of (1, 2, 3, …, n) to itself is a semigroup under composition of mappings}, called the symmetric semigroup. The number of elements in S(n) is $n^n$.

Now we proceed on to recall the definition of a group.

**DEFINITION 1.2.2:** *A non empty set G is said to form a group if on G is defined an associative binary operation * such that*

1. $a, b \in G$ then $a * b \in G$
2. There exists an element $e \in G$ such that $a * e = e * a = a$ for all $a \in G$.
3. For every $a \in G$ there is an element $a^{-1}$ in $G$ such that $a * a^{-1} = a^{-1} * a = e$ (existence of inverse in $G$).

*A group $G$ is called an abelian or commutative if $a * b = b * a$ for all $a, b, \in G$.*

***Example 1.2.5:*** Let $G = \{1, -1\}$, $G$ is a group under multiplication.

***Example 1.2.6:*** Let $G = Z$ be the set of positive and negative integers. $G$ is a abelian group under addition.

***Example 1.2.7:*** Let $Z_n = \{0, 1, 2, …, n-1\}$; $Z_n$ is an abelian group under addition modulo n. $n \in N$.

***Example 1.2.8:*** Let $G = Z_p \setminus \{0\} = \{1, 2, …, p-1\}$, p a prime number $G = Z_p \setminus \{0\}$ is a group under multiplication of even order ($p \neq 2$)

***Example 1.2.9:*** Let $S_n = \{$group of all one to one mappings of $(1, 2, …, n)$ to itself$\}$; $S_n$ is a group under the composition maps. $o(S_n) = \underline{n}$. $S_n$ is called the permutation group or symmetric group of degree n.

***Example 1.2.10:*** Let $A_n$ be the set of all even permutations. $A_n$ is a subgroup of $S_n$ called the alternating subgroup of $S_n$, $o(An) = \underline{n}/2$.

***Example 1.2.11:*** Let $D_{2n} = \{a, b \mid a^2 = b^n = 1; bab = a\}$; $D_{2n}$ is the dihedral group of order 2n. $D_{2n}$ is not abelian ($n \in N$).

***Example 1.2.12:*** Let $G = \langle g \mid g^n = 1 \rangle$, $G$ is the cyclic group of order n, $G = \{1, g, g^2, …, g^{n-1}\}$.

***Example 1.2.13:*** Let $G = G_1 \times G_2 \times G_3 = \{(g_1, g_2, g_3) \mid g_i \in G_i; 1 \leq i \leq 3\}$ where $G_i$'s groups $1 \leq i \leq 3$. $G$ is a group.

**DEFINITION 1.2.3:** *Let (G, \*) be a group. H a proper subset of G. If (H, \*) is a group then we call H to be subgroup of G.*

***Example 1.2.14:*** Let $Z_{10} = \{0, 1, 2, \ldots, 9\}$ be the group under addition modulo 10. $H = \{0, 2, 4, 6, 8\}$ is a proper subset of $Z_{10}$ and H is a subgroup of G under addition modulo 10.

***Example 1.2.15:*** Let $D_{26} = \{a, b / a^2 = b^6 = 1, bab = a\}$ be the dihedral group of order 12. $H = \{1, b, b^2, b^3, b^4, b^5\}$ is a subgroup of $D_{2.6}$. Also $H_1 = \{1, ab\}$ is a subgroup of $D_{26}$.

For more about properties of groups please refer Hall Marshall (1961). Now we proceed on to recall the definition of the notion of Smarandache semigroups (S-semigroups).

**DEFINITION 1.2.4:** *Let $(S_i, o)$ be a semigroup. Let H be a proper subset of S. If (H, o) is a group, then we call (S, o) to be a Smarandache semigroup (S-semigroup).*

We illustrate this situation by some examples.

***Example 1.2.16:*** Let S(7) = {The mappings of the set (1, 2, 3, …, 7) to itself, under the composition of mappings} be the semigroup. $S_7$ the set of all one to one maps of (1, 2, 3, …, 7) to itself is a group under composition of mappings.
Clearly $S_7$ is a subset of S(7). Thus S(7) is a S-semigroup.

***Example 1.2.17:*** Let $Z_{12} = \{0, 1, 2, \ldots, 11\}$ be the semigroup under multiplication modulo 12. Take $H = \{1, 11\} \subseteq Z_{12}$, H is a group under multiplication modulo 12. Thus $Z_{12}$ is a S-semigroup.

***Example 1.2.18:*** Let $Z_{15} = \{0, 1, 2, \ldots, 14\}$ be the semigroup under multiplication modulo 15. Take $H = \{1, 14\} \subseteq Z_{15}$, H is a group. Thus $Z_{15}$ is a S-semigroup. $P = \{5, 10\} \subseteq Z_{15}$ is group of $Z_{15}$. So $Z_{15}$ is a S-semigroup.

**DEFINITION 1.2.5:** *Let G be a non commutative group. For h, g $\in$ G there exist x $\in$ G such that g = x h $x^{-1}$, then we say g and h are conjugate with each other.*

For more about this concept please refer I.N.Herstein (1975).

Now we proceed onto define groupoids.

**DEFINITION 1.2.6:** *Let G be a non empty set. If * be a binary operation of G such that for all a, b $\in$ G, a * b $\in$ G and if in general a * (b*c) $\neq$ (a * b) * c for a, b, c $\in$ G. Then we call (G, *) to be a groupoid. We say (G, *) is commutative if a*b = b * a for all a, b $\in$ G.*

*Example 1.2.19:* Let G be a groupoid given by the following table.

| * | 0 | 1 | 2 | 3 | 4 |
|---|---|---|---|---|---|
| 0 | 0 | 2 | 4 | 1 | 3 |
| 1 | 2 | 4 | 1 | 3 | 0 |
| 2 | 4 | 1 | 3 | 0 | 2 |
| 3 | 1 | 3 | 0 | 2 | 4 |
| 4 | 3 | 0 | 2 | 4 | 0 |

G is a commutative groupoid.

*Note:* We say a groupoid G has zero divisors if a * b = 0 for a, b $\in$ G \ {0} where o $\in$ G.

We say if e $\in$ G such that a * e = e * a = a for all a $\in$ G then we call G to be a monoid or a semigroup with identity. If for a $\in$ G there exists b $\in$ G such that a * b = b * a = e then we say a is a unit in G.

*Example 1.2.20:* Let G = {0, 1, 2, …, 9} define '*' on G by a * b = 8a + 4b (mod 10), a, b $\in$ G \ {0}. (G, *) is a groupoid.

We can have classes of groupoids built using $Z_n$.

Next we proceed onto define loops.

**DEFINITION 1.2.7:** *A non empty set L is said to form a loop if on L is defined a binary non associative operation called the product denoted by * such that*
1. *For all a, b ∈ L, a * b ∈ L.*
2. *There exists an element e ∈ L such that a * e = e * a = a for all a ∈ L. e is called the identity element of L.*
3. *For every ordered pair (a, b) ∈ L × L there exists a unique pair (x, y) ∈ L such that ax = b and ya = b.*

We give some examples.

*Example 1.2.21:* Let L = {e, 1, 2, 3, 4, 5}. The loop using L is given by the following table

| * | e | 1 | 2 | 3 | 4 | 5 |
|---|---|---|---|---|---|---|
| e | e | 1 | 2 | 3 | 4 | 5 |
| 1 | 1 | e | 5 | 4 | 3 | 2 |
| 2 | 2 | 3 | e | 1 | 5 | 4 |
| 3 | 3 | 5 | 4 | e | 2 | 1 |
| 4 | 4 | 2 | 1 | 5 | e | 3 |
| 5 | 5 | 4 | 3 | 2 | 1 | e |

*Example 1.2.22:* Let L = {e, 1, 2, 3, …, 7}. The loop is given by the following table.

| * | e | 1 | 2 | 3 | 4 | 5 | 6 | 7 |
|---|---|---|---|---|---|---|---|---|
| e | e | 1 | 2 | 3 | 4 | 5 | 6 | 7 |
| 1 | 1 | e | 4 | 7 | 3 | 6 | 2 | 5 |
| 2 | 2 | 6 | e | 5 | 1 | 4 | 7 | 3 |
| 3 | 3 | 4 | 7 | e | 6 | 2 | 5 | 1 |
| 4 | 4 | 2 | 5 | 1 | e | 7 | 3 | 6 |
| 5 | 5 | 7 | 3 | 6 | 2 | e | 1 | 4 |
| 6 | 6 | 5 | 1 | 4 | 7 | 3 | e | 2 |
| 7 | 7 | 3 | 6 | 2 | 5 | 1 | 4 | e |

We now proceed onto define a special class of loops.

**DEFINITION 1.2.8:** *Let $L_n(m) = \{e, 1, 2, ..., n\}$ be a set where $n > 3$, $n$ is odd and $m$ is a positive integer such that $(m, n) = 1$ and $(m -1, n) = 1$ with $m < n$. Define on $L_n(m)$ a binary operation 'o' as follows.*

1. *$e \circ i = i \circ e = i$ for all $i \in L_n(m)$*
2. *$i^2 = i \circ i = e$ for $i \in L_n(m)$*
3. *$i \circ j = t$ where $t = (mj - (m - 1)i) \pmod{n}$ for all $i, j \in L_n(m)$; $i \neq j$; $i \neq e$ and $j \neq e$, then $L_n(m)$ is a loop under the operation o.*

We illustrate this by some example.

*Example 1.2.23:* Let $L_7(3) = \{e, 1, 2, …, 7\}$ be a loop given by the following table.

| o | e | 1 | 2 | 3 | 4 | 5 | 6 | 7 |
|---|---|---|---|---|---|---|---|---|
| e | e | 1 | 2 | 3 | 4 | 5 | 6 | 7 |
| 1 | 1 | e | 4 | 7 | 3 | 6 | 2 | 5 |
| 2 | 2 | 6 | e | 5 | 1 | 4 | 7 | 3 |
| 3 | 3 | 4 | 7 | e | 6 | 2 | 5 | 1 |
| 4 | 4 | 2 | 5 | 1 | e | 7 | 3 | 6 |
| 5 | 5 | 7 | 3 | 6 | 2 | e | 1 | 4 |
| 6 | 6 | 5 | 1 | 4 | 7 | 3 | e | 2 |
| 7 | 7 | 3 | 6 | 2 | 5 | 1 | 4 | e |

*Example 1.2.24:* Let $L_5(3)$ be the loop given by the following table.

| o | e | 1 | 2 | 3 | 4 | 5 |
|---|---|---|---|---|---|---|
| e | e | 1 | 2 | 3 | 4 | 5 |
| 1 | 1 | e | 4 | 2 | 5 | 3 |
| 2 | 2 | 4 | e | 5 | 3 | 1 |
| 3 | 3 | 2 | 5 | e | 1 | 4 |
| 4 | 4 | 5 | 3 | 1 | e | 2 |
| 5 | 5 | 3 | 1 | 4 | 2 | e |

Now we proceed on to recall the definition of rings and S-rings

**DEFINITION 1.2.9:** *Let $(R, +, o)$ be a nonempty set R with two closed binary operations + and o defined on it.*

1. *$(R, +)$ is an abelian group.*
2. *$(R, o)$ is a semigroup*
3. *$a \; o \; (b + c) = a \; o \; b + b \; o \; c$ for all $a, b, c \in R$. We call R a ring. If $(R, o)$ is a semigroup with identity (i.e., a monoid) then we say $(R, +, o)$ is a ring with unit.*

*If $a \; o \; b = b \; o \; a$ for all $a, b \in R$ then we say $(R, +, o)$ is a commutative ring.*

*Example 1.2.25:* Let $(Z, +, \times)$ is a ring, Z the set of integers.

*Example 1.2.26:* $(Q, +, \times)$ is a ring, Q the set of rationals.

*Example 1.2.27:* $Z_{30} = \{0, 1, 2, …, 29\}$ is the ring of integers modulo 30.

We recall the definition of a field.

**DEFINITION 1.2.10:** *Let $(F, +, o)$ be such that F is a nonempty set with 0 and 1. F is a field if the following conditions hold good.*

1. *$(F, +)$ is an abelian group.*
2. *$(F \setminus \{0\}, o)$ is an abelian group*
3. *$a \; o \; (b + c) = a \; o \; b + a \; o \; c$ and*
   *$(a + b) \; o \; c = a \; o \; c + b \; o \; c$*

*for all $a, b, c \in F$.*

*Example 1.2.28:* $(Q, +, \times)$ is a field, known as the field of rationals.

*Example 1.2.29:* $Z_5 = \{0, 1, 2, 3, 4\}$ is a field, prime finite field of characteristic 5.

***Example 1.2.30:*** $F = \dfrac{Z_2[x]}{\langle x^2 + x + 1 \rangle = I}$ is a quotient ring which is a finite field of characteristic two.

Now we proceed onto recall the definition of a Smarandache ring.

**DEFINITION 1.2.11:** *Let $(R, +, o)$ be a ring we say $R$ is a Smarandache ring (S-ring) if $R$ contains a proper subset $P$ such that $(P, +, o)$ is a field.*

We illustrate this situation by some simple examples.

***Example 1.2.31:*** Let $Q[x]$ be a polynomials ring. $Q[x]$ is a S-ring for $Q \subseteq Q[x]$ is a field. So $Q[x]$ is a S-ring.

***Example 1.2.32:*** Let

$$M_{2\times 2} = \left\{ \begin{pmatrix} a & b \\ c & d \end{pmatrix} \middle| a, b, c, d \in Q \right\},$$

$M_{2\times 2}$ is a ring with respect to matrix addition and multiplication. But

$$P = \left\{ \begin{pmatrix} a & 0 \\ 0 & 0 \end{pmatrix} \middle| a \in Q \right\}$$

is a proper subset of $M_{2\times 2}$ which is a field. Thus $M_{2\times 2}$ is a S-ring.

***Example 1.2.33:*** Let $R = Z_{11} \times Z_{11} \times Z_{11}$ be the ring with component wise addition and multiplication modulo 11. $P = Z_{11} \times \{0\} \times \{0\}$ is a field contained in R. Hence R is a S-ring.

It is important to note that in general all rings are not S-rings.

***Example 1.2.34:*** Z be the ring of integers. Z is not a S-ring.

Chapter Two

# GROUPS AS GRAPHS

Here we venture to express groups as graphs. From the structure of the graphs try to study the properties of groups. To describe the group in terms of a graph we exploit the notion of identity in group so we call the graph associated with the group as identity graph. We say two elements x, y in the group are adjacent or can be joined by an edge if x.y = e (e, identity element of G). Since we have in group x.y = y.x = e we need not use the property of commutatively. It is by convention every element is adjoined with the identity of the group G. If G = {g, 1 | $g^2$ = 1} then we represent this by a line as $g^2$ = 1. This is the convention we use when trying to represent a group by a graph. The vertices corresponds to the elements of the group, hence the order of the group G corresponds to the number of vertices in the identity graph.

***Example 2.1:*** Let $Z_2$ = {0, 1} be the group under addition modulo 2. The identity graph of $Z_2$ is as 1+1 = 0(mod 2), 0 is the identity of $Z_2$.

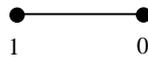

Figure 2.1

***Example 2.2:*** Let $Z_3 = \{0, 1, 2\}$ be the group under addition modulo three. The identity graph of $Z_3$ is

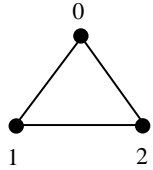

Figure 2.2

***Example 2.3:*** Let $Z_4 = \{0, 1, 2, 3\}$ be the group under addition modulo four. The identity graph of the group $Z_4$ is

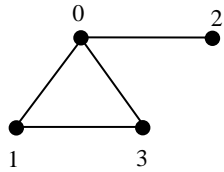

Figure 2.3

***Example 2.4:*** Let $G = \langle g \mid g^6 = 1 \rangle$ the cyclic group of order 6 under multiplication.

The identity graph of G is

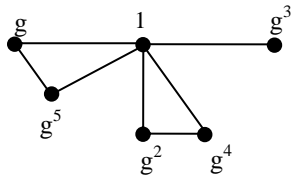

Figure 2.4

***Example 2.5:*** Let

$$S_3 = \left\{ e = \begin{pmatrix} 1 & 2 & 3 \\ 1 & 2 & 3 \end{pmatrix}, \ p_1 = \begin{pmatrix} 1 & 2 & 3 \\ 1 & 3 & 2 \end{pmatrix}, \right.$$

$$p_2 = \begin{pmatrix} 1 & 2 & 3 \\ 3 & 2 & 1 \end{pmatrix}, \; p_3 = \begin{pmatrix} 1 & 2 & 3 \\ 2 & 1 & 3 \end{pmatrix},$$

$$p_4 = \begin{pmatrix} 1 & 2 & 3 \\ 2 & 3 & 1 \end{pmatrix} \text{ and } p_4 = \begin{pmatrix} 1 & 2 & 3 \\ 3 & 1 & 2 \end{pmatrix} \Big\}$$

be the symmetric group of degree three. The identity graph associated with $S_3$ is

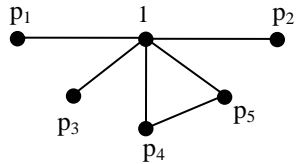

Figure 2.5

We see the groups $S_3$ and G are groups of order six but the identity graphs of $S_3$ and G are not identical.

***Example 2.6:*** Let $D_{2.3} = \{a, b \mid a^2 = b^3 = 1; b \, a \, b = a\}$ be the dihedral group. $o(D_{2.3}) = 6$.

The identity graph associated with $D_{2.3}$ is given below.

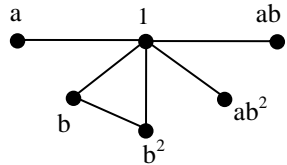

Figure 2.6

We see the identity graph of $D_{2.3}$ and $S_3$ are identical i.e., one and the same.

***Example 2.7:*** Let $G = \langle g \mid g^8 = 1 \rangle$. The identity graph of the cyclic group G is given by

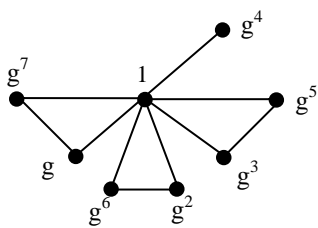

Figure 2.7

***Example 2.8:*** Let $G = \langle h \mid h^7 = 1 \rangle$ be the cyclic group of order 7. The identity graph of G is

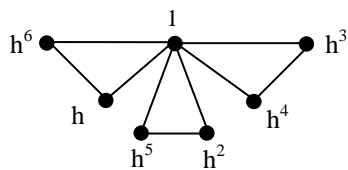

Figure 2.8

***Example 2.9:*** Let $Z_7 = \{0, 1, 2, \ldots, 6\}$, the group of integers modulo 7 under addition. The identity graph of $Z_7$ is

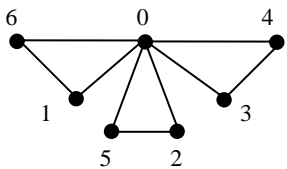

Figure 2.1.9

It is interesting to observe that the identity graph of $Z_7$ and G in example 2.8 are identical.

***Example 2.10:*** Let $G = \langle g \mid g^{12} = 1 \rangle$ be the cyclic group of order 12.
The identity graph of G is

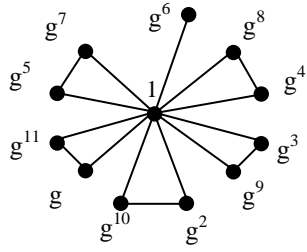

Figure 2.10

*Example 2.11:* Let

$$A_4 = \left\{ e = \begin{pmatrix} 1 & 2 & 3 & 4 \\ 1 & 2 & 3 & 4 \end{pmatrix}, h_1 = \begin{pmatrix} 1 & 2 & 3 & 4 \\ 2 & 1 & 4 & 3 \end{pmatrix}, \right.$$

$$h_2 = \begin{pmatrix} 1 & 2 & 3 & 4 \\ 4 & 3 & 2 & 1 \end{pmatrix}, h_3 = \begin{pmatrix} 1 & 2 & 3 & 4 \\ 3 & 4 & 1 & 2 \end{pmatrix},$$

$$h_4 = \begin{pmatrix} 1 & 2 & 3 & 4 \\ 1 & 3 & 4 & 2 \end{pmatrix}, h_5 = \begin{pmatrix} 1 & 2 & 3 & 4 \\ 1 & 4 & 2 & 3 \end{pmatrix},$$

$$h_6 = \begin{pmatrix} 1 & 2 & 3 & 4 \\ 3 & 2 & 4 & 1 \end{pmatrix}, h_7 = \begin{pmatrix} 1 & 2 & 3 & 4 \\ 4 & 2 & 1 & 3 \end{pmatrix},$$

$$h_8 = \begin{pmatrix} 1 & 2 & 3 & 4 \\ 2 & 4 & 3 & 1 \end{pmatrix}, h_9 = \begin{pmatrix} 1 & 2 & 3 & 4 \\ 4 & 1 & 3 & 2 \end{pmatrix},$$

$$h_{10} = \begin{pmatrix} 1 & 2 & 3 & 4 \\ 2 & 3 & 1 & 4 \end{pmatrix}, h_{11} = \left. \begin{pmatrix} 1 & 2 & 3 & 4 \\ 3 & 1 & 2 & 4 \end{pmatrix} \right\}$$

be the alternating group of order 12.

The identity graph associated with $A_4$ is

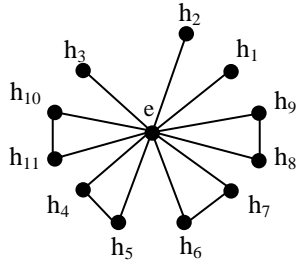

Figure 2.11

It is interesting to observe both $A_4$ and G are of order 12 same order but the identity graph of both $A_4$ and G are not identical.

We now find the identity graph of $Z_{12}$, the set of integers modulo 12.

***Example 2.12:*** Let $Z_{12} = \{0, 1, 2, \ldots, 11\}$ be the group under addition modulo 12.

The identity graph of $Z_{12}$ is

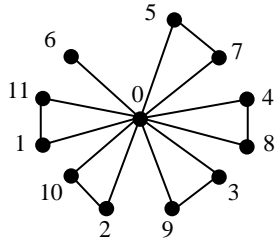

Figure 2.12

We see the identity graph of $Z_{12}$ and G given in example 2.10 are identical.

Now we see the identity graph of $D_{2.6}$.

***Example 2.13:*** The dihedral group $D_{2.6} = \{a, b \mid a^2 = b^6 = 1, b\,a\,b = a\} = \{1, a, b, ab, ab^2, ab^3, ab^4, ab^5, b^2, b^3, b^4, b^5\}$.

The identity graph of $D_{2.6}$ is

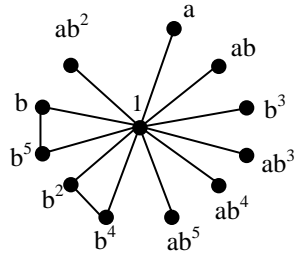

Figure 2.13

We see the identity graph of $D_{26}$ is different from that of $A_4$, $Z_{12}$ and G, though $o(D_{26}) = 12$.

*Example 2.14:* The identity graph of the cyclic group $G = \langle g \mid g^{14} = 1 \rangle$ is as follows

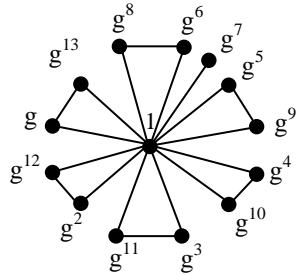

Figure 2.14

*Example 2.15:* The identity graph of the dihedral group $D_{2.7} = \{a, b \mid a^2 = b^7 = 1, bab = a\}$ is as follows:

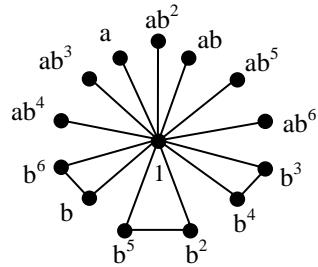

Figure 2.15

We see $o(D_{2.7}) = o(G) = 14$, but the identity graph of $D_{2.7}$ and G are not identical.

***Example 2.16:*** Let $Z_{17} = \{0, 1, 2, \ldots, 16\}$ be the group under addition modulo 17.

The identity graph of $Z_{17}$ is

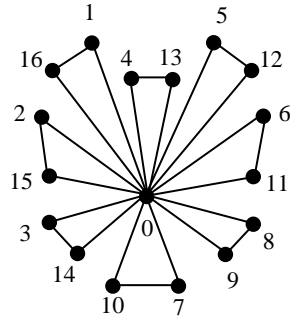

Figure 2.16

***Example 2.17:*** Let $G = Z_{17} \setminus \{0\} = \{1, 2, \ldots, 16\}$ be the group under multiplication modulo 17.

The identity graph associated with G is

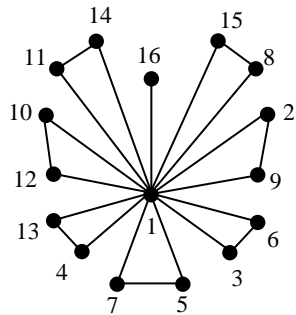

Figure 2.17

***Example 2.18:*** Let $G' = \langle g \mid g^{16} = 1 \rangle$ be the cyclic group of order 16.

The identity graph of G′.

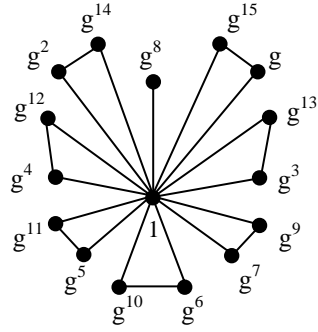

Figure 2.18

We see the identity graphs of G and G′ are identical.

***Example 2.19:*** Let $\overline{G} = H \times K = \{1, g \mid g^2 = 1\} \times \langle h \mid h^8 = 1\rangle = \{(1, 1) (1, h), (1, h^2), (1, h^3), (1, h^4), (1, h^5), (1, h^6), (1, h^7), (g, h), (g, h^2), (g, h^3), (g, h^4), (g, h^5), (g, h^6), (g, 1), (g, h^7)\}$.

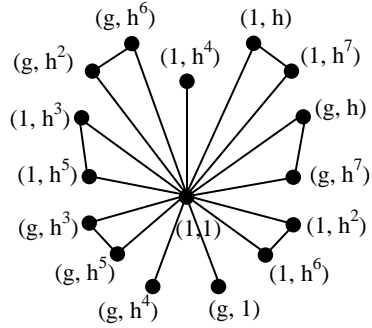

Figure 2.19

We see $\overline{G} = H \times K$ is of order 16 but the identity graph of $\overline{G}$ is different from that of G and G′ given in examples

***Example 2.20:*** Let $K = P \times Q = \{1, g, g^2, g^3 \mid g^4 = 1\} \times \{1, h, h^2, h^3 \mid h^4 = 1\}$ be the group of order 16. $K = \{(1, 1), (1, h), (1, h^2), (1, h^3), (g, h) (g, 1) (g, h^2), (g, h^3), (g^2, 1) (g^2, h^2), (g^2, h) (g^2, h^3) (g^3, 1) (g^3, h) (g^3, h^2) (g^3, h^3)\}$.

The identity graph of K is as follows.

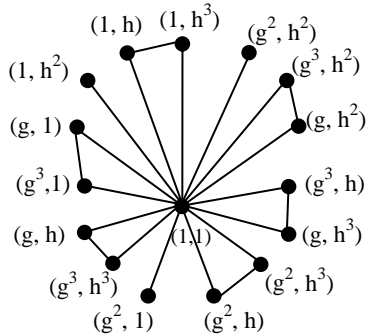

Figure 2.20

$o(K) = o(\overline{G}) = 16$.

We see the identity graphs of K and $\overline{G}$ are identical. $\overline{G}$ given in example 2.19.

***Example 2.21:*** Let $D_{2.8} = \{a, b \mid a^2 = b^8 = 1, bab = a\} = \{1, a, b, b^2, b^3, b^4, b^5, b^6, b^7, ab, ab^2, ab^3, ab^4, ab^5, ab^6, ab^7\}$ be the dihedral group of order 16. The identity graph of $D_{2.8}$ is

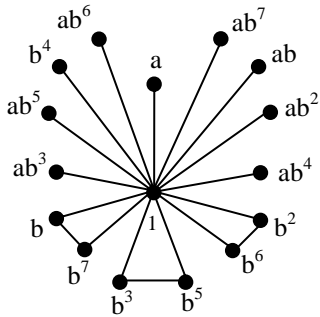

Figure 2.21

We see $o(D_{28}) = o(\overline{G}) = o(K) = o(G) = 16$.

But the identity graph of the group $D_{2.8}$ is distinctly different from that of the groups $\overline{G}$ and K given in examples.

***Example 2.22:*** Let $G = \{g \mid g^2 = 1\} \times \{h \mid h^3 = 1\} \times \{p \mid p^2 = 1\}$
$= \{(1\ 1\ 1), (g\ 1\ 1), (1\ h\ 1), (1\ h^2\ 1), (g\ h\ 1), (g\ h^2\ 1), (1\ 1\ p), (g\ h\ p), (g\ 1\ p), (g\ h^2\ p), (1\ h\ p), (1\ h^2\ p)\}$ be a group of order 12. The identity graph of G is

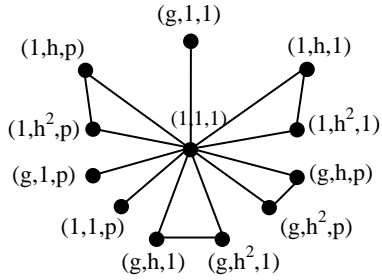

Figure 2.22

The identity graph of G is identical with that of $A_4$.

***Example 2.23:*** Let $G = \{g \mid g^2 = 1\} \times \{h \mid h^6 = 1\} = \{(1\ 1), (g\ 1), (1\ h), (1\ h^2), (1\ h^3), (1\ h^4), (1\ h^5), (g\ h), (g\ h^2), (g\ h^3), (g\ h^4), (g\ h^5)\}$ be a group of order 12. The identity graph of H is as follows.

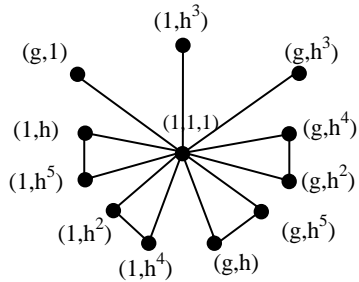

Figure 2.23

This is identical with the identity graph of the group $A_4$.

***Example 2.24:*** Let $G = Z_{30} = \{0, 1, 2, \ldots, 29\}$ be the group under addition modulo 30.

The identity graph associated with G is as follows.

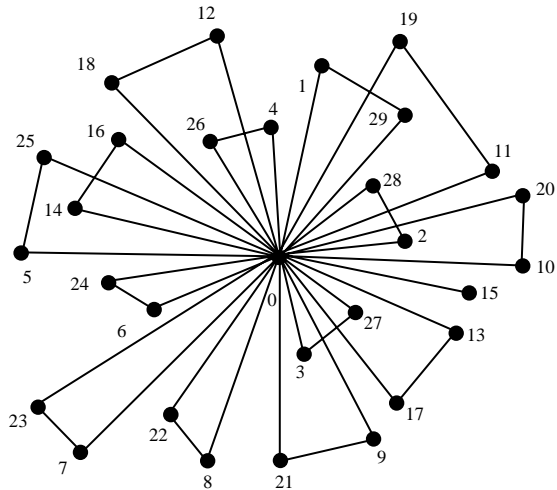

Figure 2.24

Now we can express every group G as the identity graph consisting of lines and triangles emerging from the identity element of G. Lines give the number of self inversed elements in the group. Triangles represent elements that are not self inversed.

Now we proceed on to describe the subgroup of a group by a identity graph.

**DEFINITION 2.1:** *Let G be a group. H a subgroup of G then the identity graph drawn for the subgroup H is known as the identity special subgraph of G (special identity subgraph of G).*

*Example 2.25:* Let $G = \{g \mid g^8 = 1\}$ be a group of order 8. The subgroups of G are $H_1 = \{1, g^4\}$, $H_2 = \{1, g^2, g^4, g^6\}$, $H_3 = \{1\}$ and $H_4 = G$.

The identity graphs of $H_1$, $H_2$, $H_3$ and $H_4 = G$ is as follows: The identity special subgraph of $H_3$ is just a point,

1

which is known as the trivial identity graph.

The identity graph of $H_1$ is

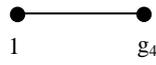

Figure 2.25.1

i.e., a line graph as $g^4$ is a self inversed element of G.

The identity graph of $H_2$ is

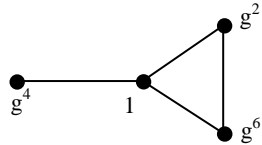

Figure 2.25.2

The identity graph of $H_4 = G$ is as follows.

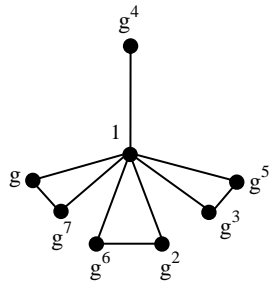

Figure 2.25.3

We see clearly the identity graphs of $H_1$, $H_2$ and $H_3$ are also identity subgraphs of G.

***Example 2.26:*** Let $D_{2.7} = \{a, b \mid a^2 = b^7 = 1, bab = a\}$ be the dihedral group of order 14.
   The identity graph of $D_{2.7}$ is as follows.

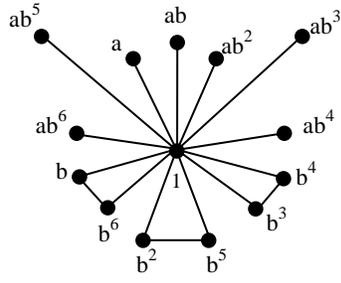

Figure 2.26

The subgroups of $D_{27}$ are $H_o = \{1, a\}$ whose identity graph is

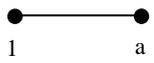

$H_i = \{1, ab^i\}$ is a subgroup of $D_{27}$ with its identity graph as

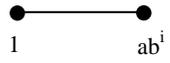

such subgroups for i = 1, 2, 3, 4, 5 and 6 are given by

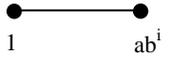

$H_7 = \{1, b, b^2, b^3, b^4, b^5, b^6\}$ is a subgroup of $D_{26}$. The identity graph associated with $H_7$ is

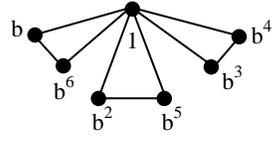

    We see the identity graphs of the subgroups are special identity subgraphs of the identity graph of $D_{27}$.
    However it is interesting to note that all subgraphs of an identity graph need not correspond with a subgroup. We have for every subgroup H of a group G a special identity subgraph of the identity graph, however the converse is not true.

**THEOREM 2.1:** *Let G be a group. $G_i$ denote the identity graph related to G. Every subgroup of G has an identity graph which is a special identity subgraph of $G_i$ and every identity subgraph of $G_i$ need not in general be associated with a subgroup of G.*

*Proof:* Given G is a group. $G_i$ the related identity graph of G. Suppose H is a subgroup of G then since H itself is a group and H a subset of G the identity graph associated with H will be a special identity subgraph of $G_i$.

Conversely if $H_i$ is a identity subgraph of $G_i$ then we may not in general have a subgroup associated with it.

Let $G = \{g \mid g^{11} = 1\}$ be the cyclic group of order 11. The identity graph associated with G be $G_i$ which is as follows:

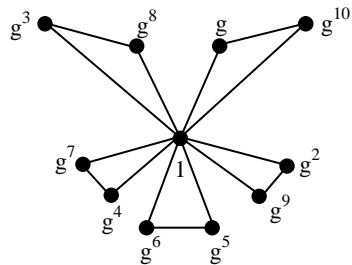

Figure 2.27

This has

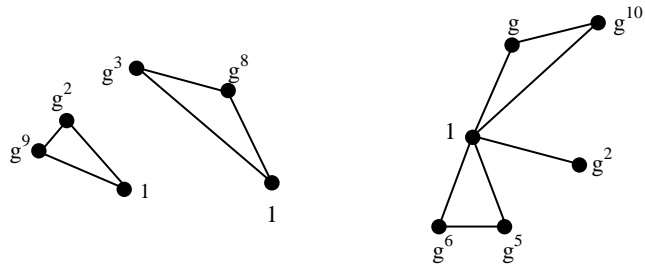

as some identity subgraphs of $G_i$. Clearly no subgroups can be associated with them as G has no proper subgroups, as $o(G) = 11$, a prime. Hence the theorem.

**Example 2.27:** Let $A_4 = \{1, h_1, h_2, h_3, h_4, h_5, h_6, h_7, h_8, h_9, h_{10}, h_{11}\}$ be the alternating subgroup of $S_4$. The identity graph $G_i$ associated with $A_4$ is as follows:

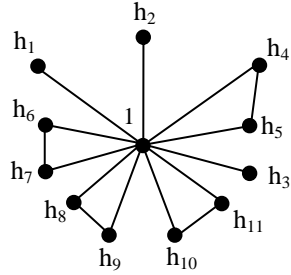

Figure 2.28

The subgroups of $A_4$ are $\{1\} = P_1$, $P_2 = A_4$, $P_3 = \{1, h_2\}$, $P_4 = \{1, h_1\}$ $P_5 = \{1, h_3\}$, $P_6 = \{1, h_1, h_2, h_3\}$, $P_7 = \{1, h_4, h_5\}$, $P_8 = \{1, h_6, h_7\}$, $P_9 = \{1, h_8, h_9\}$ and $P_{10} = \{1, h_{10}, h_{11}\}$.

The special identity subgraph $H_1$ of $P_1$ is

1

The special identity subgraph of $P_2$ is given in figure 2.28. The special identity subgraph of the subgroup $P_3$ is

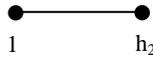

1         $h_2$

The special identity subgraph of the subgroup $P_4$ is

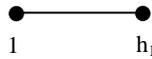

1         $h_1$

The special identity subgraph of the subgroup $P_5$ is

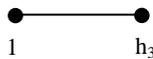

1         $h_3$

The special identity subgraph of $P_6$ is

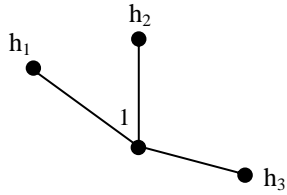

The special identity subgraph of $P_7$, the subgroup of $A_4$

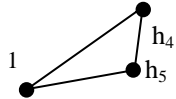

The special identity subgraph of the subgroup $P_8$ is

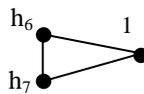

The special identity subgraph of the subgroup $P_9$ is

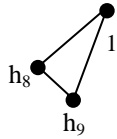

The special identity subgraph of the subgroup of $P_{10}$ is

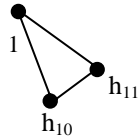

However we see some of the subgraphs of the identity graph of $G_i$ are

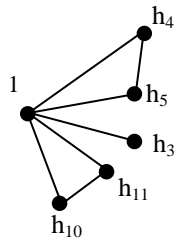

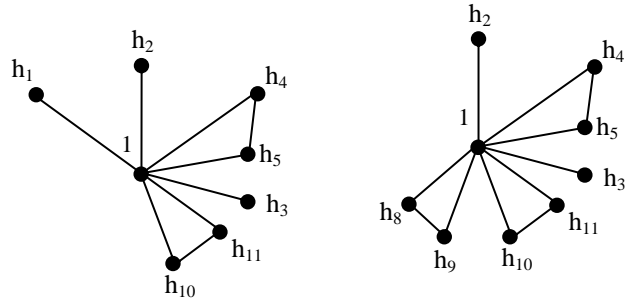

These have no subgroups of $A_4$ being associated with it.

*Example 2.28:* Let $G = \langle g \mid g^5 = 1 \rangle$ be the cyclic group of order 5. The identity graph $G_i$ associated with G is

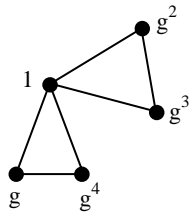

Figure 2.29

Clearly this has no identity special subgraph.

*Example 2.29:* Let $G = \langle g \mid g^{19} = 1 \rangle$ be the cyclic group of order 19. The identity graph associated with G is

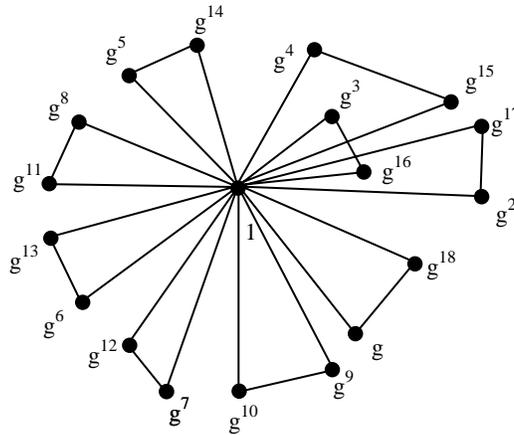

Figure 2.30

This too has no identity special subgraphs.

***Example 2.30:*** Let $G = \langle g \mid g^{15} = 1 \rangle$ be the cyclic group of order 15. The identity graph of G, is given in the following.

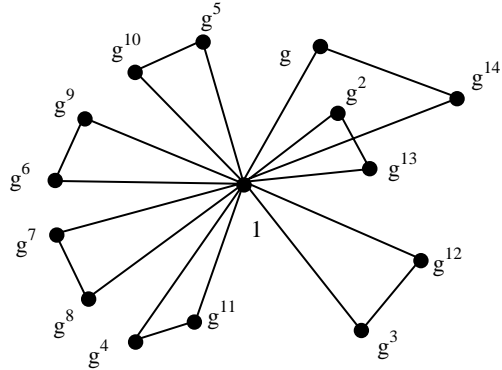

Figure 2.31

The subgroups of G are $H_1 = \{1\}$, $H_2 = G$, $H_3 = \{1, g^3, g^6, g^9, g^{12}\}$, $H_4 = \{1, g^5, g^{10}\}$. The associated special identity subgraphs of $G_i$ are

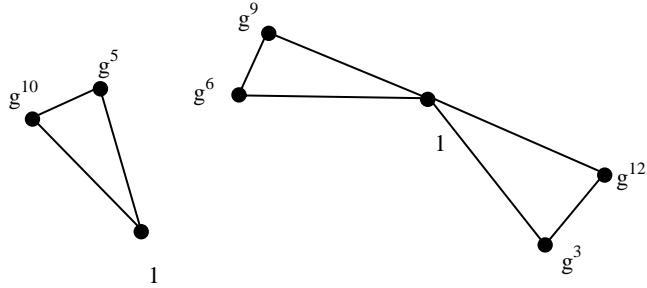

**THEOREM 2.2:** *If $G = \langle g \mid g^p = 1 \rangle$ be a cyclic group of order p, p a prime. Then the identity graph formed by G has only triangles infact $(p – 1) / 2$ triangles.*

*Proof:* Given $G = \langle g \mid g^p = 1 \rangle$ is a cyclic group of order p, p a prime. G has no proper subgroups. So no element in G is a self inversed element i.e., for no $g_i$ in G is such that $(g_i)^2 = 1$. For by Cauchy theorem G cannot have elements of order two. Thus for every $g^i$ in G there exists a unique $g^j$ in G such that $g^i g^j = 1$ i.e., for every $g^i$ the $g^j$ is such that $j = (p - i)$ so from this the elements $1, g^i, g^{p-i}$ form a triangle. Hence the identity graph will not have any line any graphs. Thus a typical identity graph of these G will be of the following form.

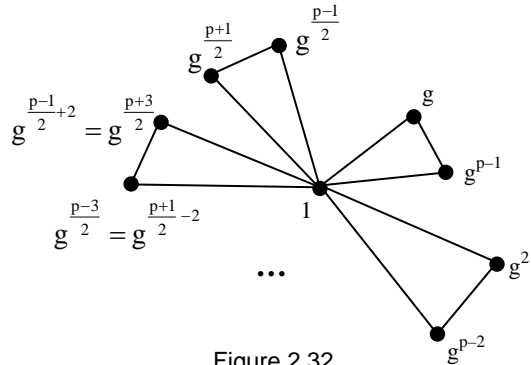

Figure 2.32

Hence the theorem.

**COROLLARY 2.1:** *If G is a cyclic group of odd order then also G has the identity graph $G_i$ which is formed only by triangles with no lines.*

*Proof:* Let $G = \{g \mid g^n = 1\}$, n is a odd number. Then the identity graph associated with G is as follows.

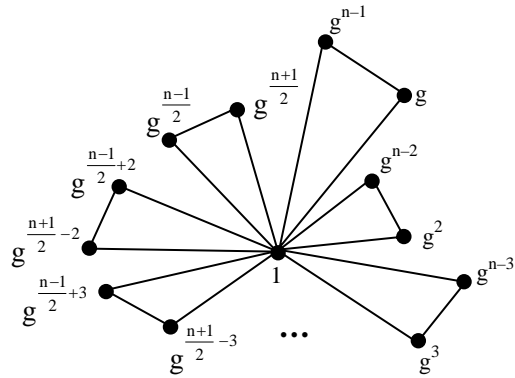

Figure 2.33

We illustrate this by examples and show $G_i$ mentioned in the corollary has special identity subgraphs.

***Example 2.31:*** Let $G = \langle g \mid g^{13} = 1 \rangle$ be the cyclic group of order 13. The identity graph $G_i$ associated with G is as follows:

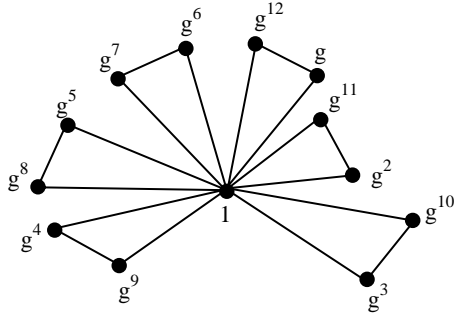

Figure 2.34

This has no special identity subgraphs.

***Example 2.32:*** Let $G = \langle g \mid g^{16} = 1 \rangle$ be the cyclic group of order 16. The identity graph of G is $G_i$ given by

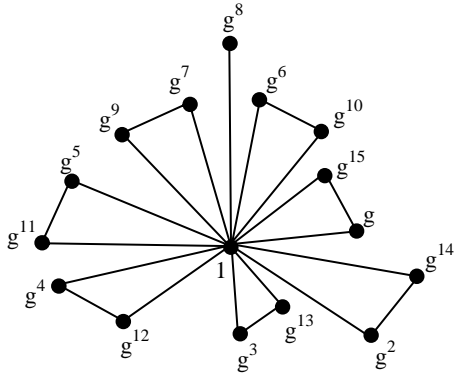

Figure 2.35

The subgroups of G are $H_1 = \{g^2, g^4, g^6, g^8, g^{10}, g^{12}, g^{14}, 1\}$, $H_2 = \{g^4, g^8, g^{12}, 1\}$ and $H_3 = \{1, g^8\}$.

The special identity subgraph of $H_1$ is

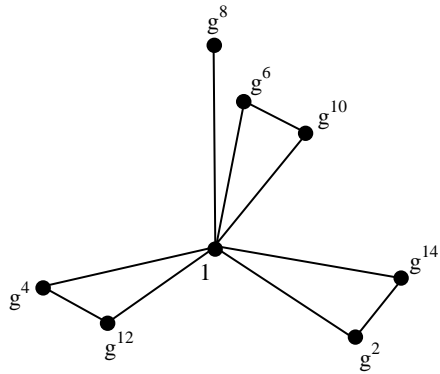

The special identity subgraph of $H_2$ is

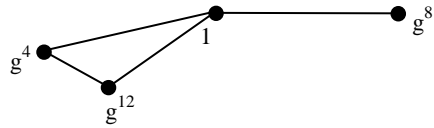

The special identity subgraph of $H_3$ is

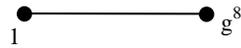

***Example 2.33:*** Let $G = \langle g \mid g^9 = 1 \rangle$ be the cyclic group of order 9. The identity graph of G is

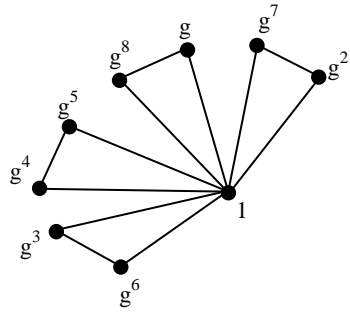

Figure 2.36

The subgroup of G is H = {1, $g^3$, $g^6$}. The special identity subgraph of H is

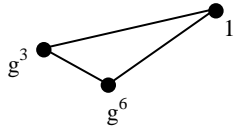

**THEOREM 2.3:** *If $G = \langle g \mid g^n = 1 \rangle$ be a cyclic group of order n, n an odd number then the identity graph $G_i$ of G is formed with $(n - 1)/2$ triangles.*

*Proof*: Clear from figure 2.32 and theorem 2.2.

**THEOREM 2.4:** *Let $G = \langle g \mid g^m = 1 \rangle$ be a cyclic group of order m; m an even number. Then the identity graph $G_i$ has $(m - 2)/2$ triangles and a line.*

*Proof:* Given $G = \langle g \mid g^m = 1 \rangle$ is a cyclic group of order m where m is even. The identity graph $G_i$ associated with G is given below:

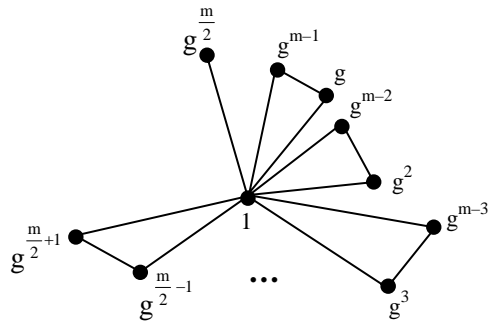

Figure 2.37

If it easily verified that exactly $(m - 1)/2$ triangles are present and only a line connecting 1 to $g^{m/2}$ for $g^{m/2}$ is a self inversed element of G.

*Example 2.34:* Let $G = \langle g \mid g^6 = 1 \rangle$ be the cyclic group of order 6. The identity graph associated with G is given by

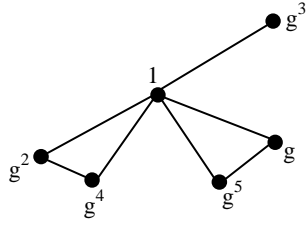

Figure 2.38

***Example 2.35:*** Let $G = \langle g \mid g^4 = 1 \rangle$ be the cyclic group of order 4. The identity graph associated with G is

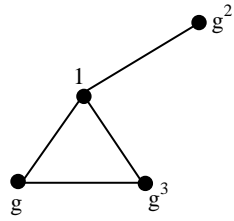

Figure 2.39

***Example 2.36:*** Let $G = \langle g \mid g^{10} = 1 \rangle$ the cyclic group of order 10. The identity graph associated with G is

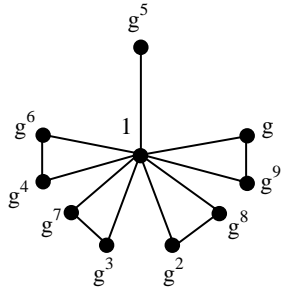

Figure 2.40

***Example 2.37:*** Let $G = \{g \mathbin{/} g^4 = 1\} \times \{g_1 \mathbin{/} g_1^6 = 1\} = \{(1, 1),$ $(g, 1), (g, g_1), (g, g_1^2), (g, g_1^3), (g, g_1^4), (g, g_1^5), (g^2, g_1), (g^2, 1),$ $(g^2, g_1^2), (g^2, g_1^3), (g^2, g_1^4), (g^2, g_1^5), (g^3, 1), (g^3, g_1), (g^3, g_1^2), (g^3,$ $g_1^3), (g^3, g_1^4), (g^3, g_1^5), (1, g_1), (1, g_1^2), (1, g_1^3), (1, g_1^4), (1, g_1^5)\}$ be the group of order 24. The identity graph of G is as follows.

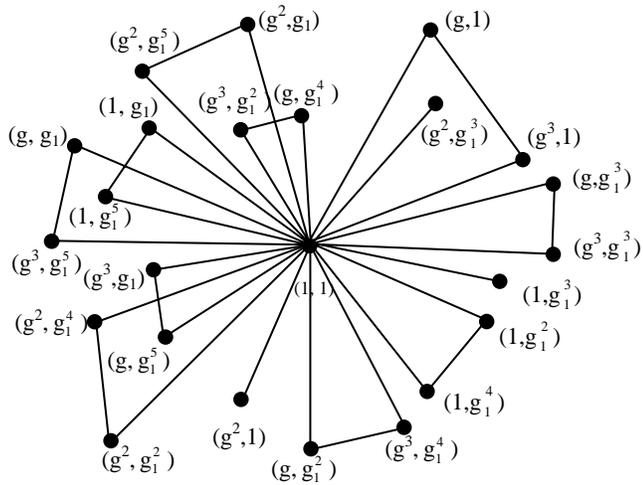

Figure 2.41

This group is not cyclic has twenty four elements. This has three lines and the rest are triangles.

***Example 2.38:*** Let $D_{2.10} = \{a, b \,/\, a^2 = b^{10} = 1, bab = a\} = \{1, a, b, b^2, b^3, b^4, b^5, b^6, b^7, b^8, b^9, ab, ab^2, ab^3, ab^4, ab^5, ab^6, ab^7, ab^8, ab^9\}$ be the dihedral group. The identity graph of $D_{2,10}$ is as follows.

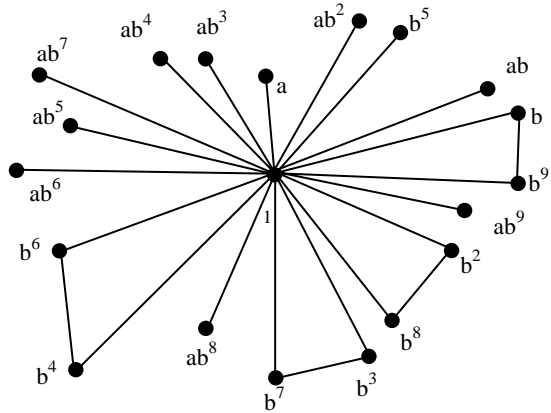

Figure 2.42

We can now bring an analogue of the converse of the Lagrange's theorem.

**THEOREM 2.5:** *Let G be a finite group. $G_i$ be the identity graph related with G. Just as for every divisor d of the order of G we do not have subgroups of order d we can say corresponding to every identity subgraph $H_i$ of the identity graph $G_i$ we may not in general have a subgroup of G associated with $H_i$.*

*Proof:* This can be proved only by an example. Let $G = D_{27} = \{a, b / a^2 = b^7 = 1, bab = a\} = \{1, a, b, b^2, b^3, b^4, b^5, b^6, ab, ab^2, ab^3, ab^4, ab^5, ab^6\}$ be the dihedral group of order 14. Let $G_i$ be the identity graph associated with $D_{27}$.

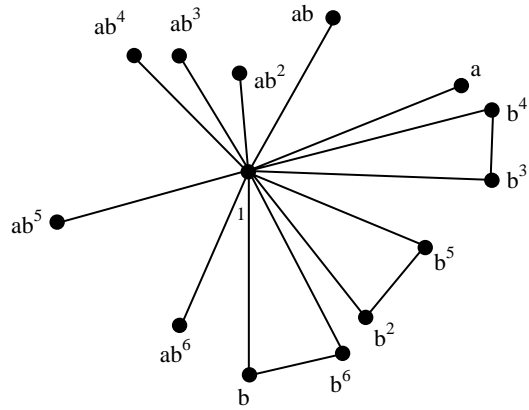

Figure 2.43

Let $H_i$ be the subgraph

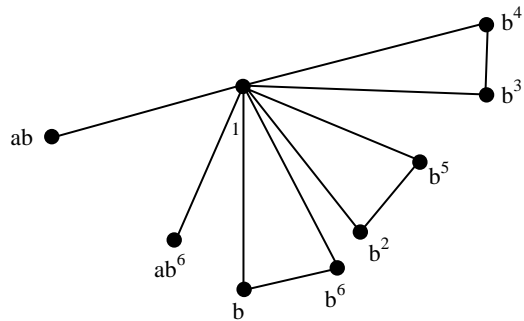

Clearly $H_i$ is the identity subgraph of $G_i$ but the number of vertices in $H_i$ is 9, i.e., the order of the subgroup associated with $H_i$ if it exists is of 9 which is an impossibility as $(14, 9) = 1$ i.e., $9 \nmid 14$.

Thus we have for a given identity graph $G_i$ of a finite group G to each of the subgraphs $H_i$ of $G_i$ we need not in general have a subgroup H of G associated with it. Hence the theorem.

**COROLLARY 2.2:** *Take $G = \langle g \mid g^p = 1 \rangle$, be a cyclic group of order p, p a prime, $G_i$ the identity graph with p vertices associated with it. G has several identity subgraphs but no subgroup associated with it.*

**COROLLARY 2.3:** *Let G be a finite group of order n. $G_i$ the identity graph of G with n vertices. Let H be a subgroup of G of order m (m < n). Suppose $H_i$ is the identity subgraph associated with H having m vertices. Every subgraph of $G_i$ with m vertices need not in general be the subgraph associated with the subgroup H.*

*Proof:* The proof is only by an example.

Let $D_{2.11} = \{a, b \mid a^2 = b^{11} = 1, bab = a\} = \{1, a, b, b^2, b^3, b^4, b^5, b^6, b^7, b^8, b^9, b^{10}, ab, ab^2, ab^3, ab^4, ab^5, ab^6, ab^7, ab^8, ab^9, ab^{10}\}$ be the dihedral group of order 22. Let $G_i$ be the identity graph associated with $D_{2.11}$ which has 22 vertices.

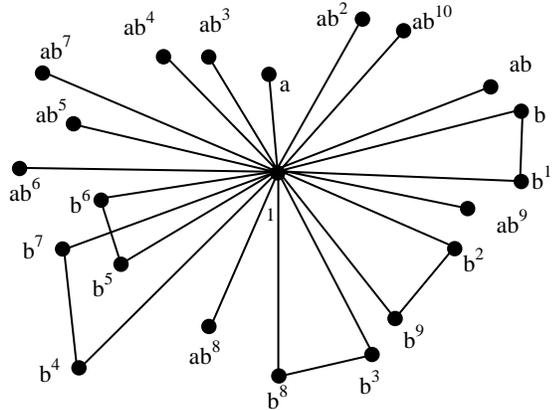

Figure 2.44

Let H = {1, b, $b^2$, …, $b^9$, $b^{10}$} be the subgroup of $D_{2.11}$ of order 11. The identity graph associated with H be $H_i$, $H_i$ is a subgraph of $G_i$ with 11 vertices.

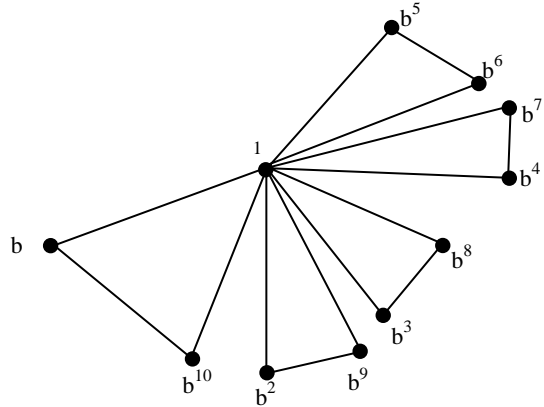

Take $P_i$ a subgraph of $G_i$ with 11 vertices viz.

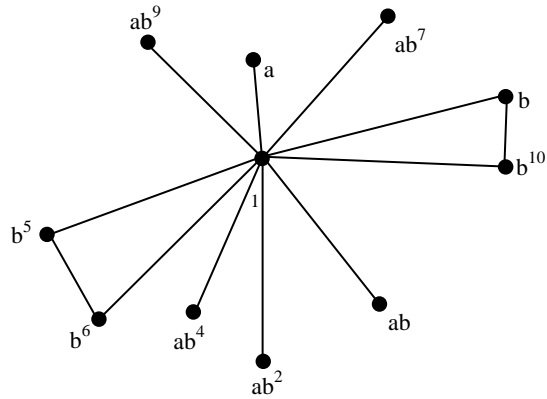

$P_i$ is a subgraph of $G_i$ but $P_i$ has no subgroup of $D_{2.11}$ associated with it.

Clearly $D_{2.11}$ has only one subgroup H of order 11 and only the graph $H_i$ is a special identity subgraph with 11 vertices and $G_i$ has no special identity subgraph with 11 vertices. However it can have identity subgraphs with 11 vertices.

In view of the above corollary we can have the following nice characterization of the dihedral groups and its subgroup of order 2p where p is a prime.

**THEOREM 2.6:** *Let $D_{2p} = \{1, a, b \mid a^2 = b^p = 1; bab = a\} = \{1, a, b, b^2, ..., b^{p-1}, ab, ab^2, ..., ab^{p-1}\}$ be the dihedral group of order 2p, p a prime. Let $G_i$ be the identity graph of $D_{2p}$ with 2p vertices. $H_i$ be a special identity subgraph of $G_i$ with p-vertices. Then $H_i$ is formed by $(p - 1)/2$ triangles which will have p vertices. This $H_i$ is unique special identity subgraph of $G_i$.*

*Proof:* Given $D_{2p}$ is a dihedral group of order 2p, p a prime so order of $D_{2p}$ is 2p and $D_{2p}$ has one and only one subgroup of order p.

Let $G_i$ be the identity graph associated with $D_{2p}$. It has 2p vertices and it is of the following form.

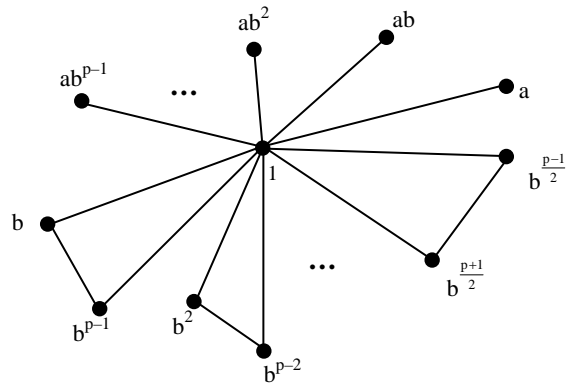

Figure 2.45

Thus the graph $G_i$ has p lines and p/2 triangles which comprises of (p + 1) vertices for the p lines and p vertices for the triangles the central vertex 1 is counted twice so the total vertices = vertices contributed by the p lines + vertices contributed by the p/2 triangles – (one vertex), which is counted twice.

This add ups to $p + 1 + p - 1 = 2p$. Thus the special identity subgraph of $G_i$ is one formed by the p/2 triangles given by which has p vertices given by the subgroup $H = \{1, b, b^2, \ldots, b^{p-1}\}$. Hence the claim. This subgroup H is unique so also is the special identity subgraph $H_i$ of $G_i$.

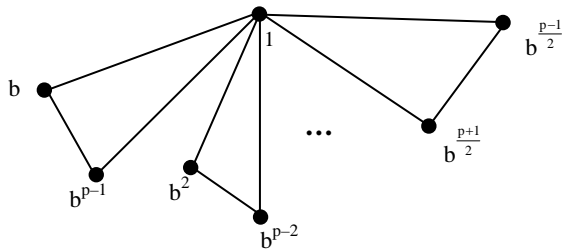

*Example 2.39:* Let $G = Z_4 \times Z_3 = \{(0, 0), (0, 1), (0, 2), (1, 0), (1, 1), (1, 2), (2, 0), (2, 1), (2, 2), (3, 0), (3, 1), (3, 2)\}$ be the direct product group under component wise addition of two additive groups $Z_4$ and $Z_3$.

The identity subgraph of G is

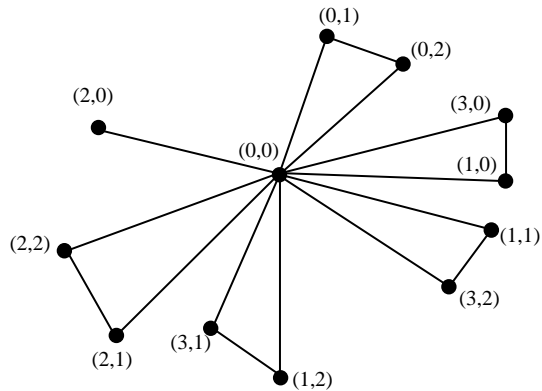

Figure 2.46

The subgroup $H = \{(0, 0), (0, 1), (2, 0), (2, 1), (2, 2), (0, 2)\}$ of G is of order 6.

The special identity subgraph associated with H is given by

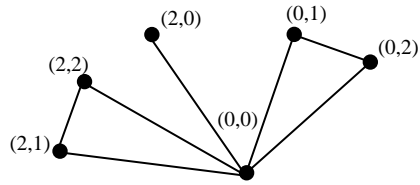

However the following

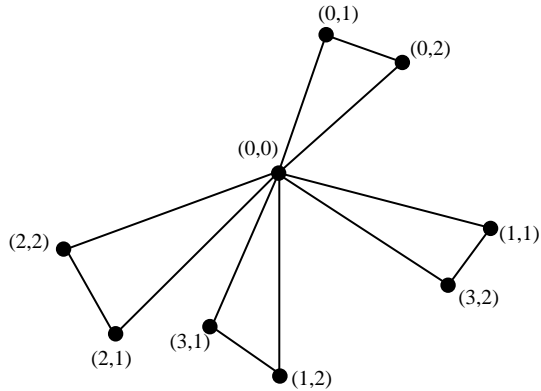

subgraph does not yield a subgroup in G hence it cannot be a special identity subgraph of G.

**DEFINITION 2.2:** *Let G be a group, $G_i$ the identity graph associated with G. H be a normal subgroup of G; then $H_i$ the special identity subgraph of H is defined to be the special identity normal subgroup of $G_i$.*

*If G has no normal subgroups then we define the identity graph $G_i$ to be a identity simple graph. Thus if G is a simple group the identity graph associated with G is a identity simple graph.*

*Example 2.40:* Let $G = \langle g \mid g^{13} = 1 \rangle$ be the cyclic group of order 13. Clearly G is a simple group.

The identity graph associated with G be $G_i$ which is as follows:

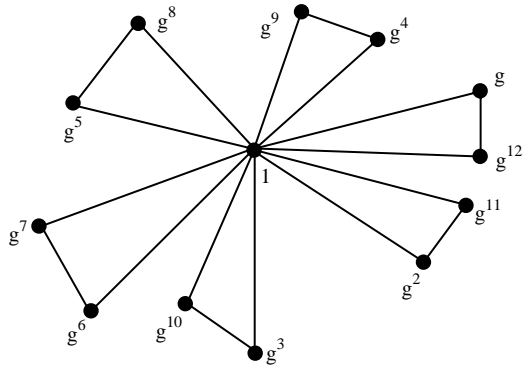
Figure 2.47

Clearly $G_i$ is a identity simple graph.

**Example 2.41:** Let $G = \langle g \mid g^{15} = 1 \rangle$ be a group of order 15. The normal subgroups of G are $H = \{1, g^3, g^6, g^9, g^{12}\}$ and $K = \{1, g^5, g^{10}\}$. The identity graph of G be $G_i$ which is as follows:

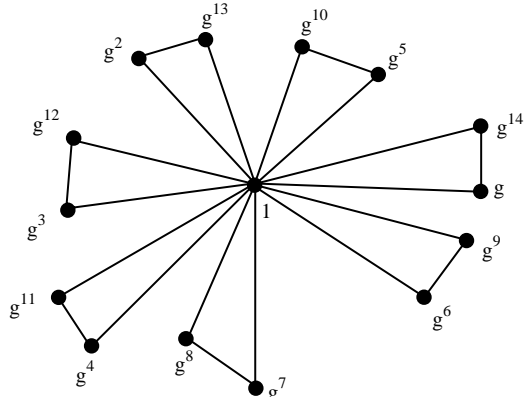
Figure 2.48

The special normal identity subgraph $H_i$ of $G_i$ is given in figure.

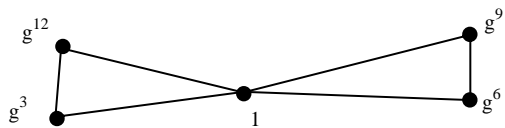

The special normal identity subgraph $K_i$ of $G_i$ is

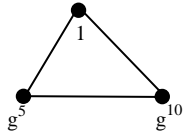

We now discuss of the vertex or (and) edge colouring of the identity graph related with a group in the following way. First when the group G is finite which we are interested mostly; we say subgroup colouring of the vertex (or edges) if we colour the verices (or edges) of those subgroups $H_i$ of G such that $H_i \not\subseteq H_j$ or $H_j \not\subseteq H_i$. We take only subgroups of G and not the subgroups of subgroups of G.

We now define the special clique of a group G.

**DEFINITION 2.3:** *Let G be a group $S = \{H_1, ..., H_n \mid H_i$'s are subgroups of G such that $H_i \not\subseteq H_j$ or $H_j \not\subseteq H_i$; $H_i \cap H_j = \{e\}$ with $G = \cup H_i$ for $i, j \in \{1, 2, ..., n\}$; i.e., the subgroup $H_i$ is not a subgroup of any of the subgroups $H_j$ for $j=1, 2, ..., j$; $i \neq j$ true for every $i$, $i=1, 2, ..., n\}$.*

*We call S a clique of the group G, if G contains a clique with n element and every clique of G has at most n elements. If G has a clique with n elements then we say clique $G = n$ if $n = \infty$ then we say clique $G = \infty$. We assume that all the vertices of each subgroup $H_i$ is given the same colour. The identity element which all subgroups have in common can be given any one of the colours assumed by the subgroups $H_i$.*

*The map $C: S \to T$ such that $C(H_i) \neq C(H_j)$ when ever the subgroups $H_i$ and $H_j$ are adjacent and the set T is the set of available colours.*

*All that interests us about T is its size; typically we seek for the smallest integer k such that S has a k-colouring, a vertex colouring $C: S \to \{1, 2, ..., k\}$. This k is defined to be the special identity chromatic number of the group G and is denoted by $\chi(S)$. The identity graph $G_i$ with $\chi(S) = k$ is called k-chromatic; if $\chi(S) \leq k$ and the group G is k-colourable.*

We illustrate this by a few examples.

***Example 2.42:*** Let

$$S_3 = \left\{ \begin{pmatrix} 1 & 2 & 3 \\ 1 & 2 & 3 \end{pmatrix} = e, \begin{pmatrix} 1 & 2 & 3 \\ 1 & 3 & 2 \end{pmatrix} = p_1, \right.$$

$$\begin{pmatrix} 1 & 2 & 3 \\ 3 & 2 & 1 \end{pmatrix} = p_2, \begin{pmatrix} 1 & 2 & 3 \\ 2 & 1 & 3 \end{pmatrix} = p_3,$$

$$\left. \begin{pmatrix} 1 & 2 & 3 \\ 2 & 3 & 1 \end{pmatrix} = p_4 \text{ and } \begin{pmatrix} 1 & 2 & 3 \\ 3 & 1 & 2 \end{pmatrix} = p_5 \right\}$$

be the symmetric group of $S_3$. $S = \{H_1 = \{e, p_1\}, H_2 = \{e, p_2\}, H_3 = \{e, p_3\}, H_4 = \{p_4, p_5, e\}$.

The identity graph $G_i$ of $S_3$ is

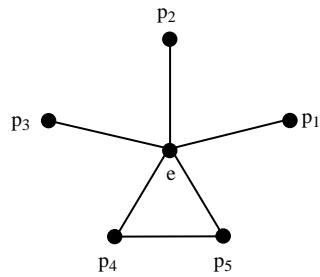

Figure 2.49

Here two colours are sufficient for $p_3$ and $p_1$ can be given one colour and $p_2$, $p_4$ and $p_5$ another colour so $k = 2$, for the identity graph given in figure 2.49 for the group $S_3$.

Thus $\chi(S) = 2$.

***Example 2.43:*** Let $G = \{1, 2, 3, 4, 5, 6\} = Z_7 \setminus \{0\}$ be the group under multiplication modulo 7. The subgroups of G are $H_1 = \{1, 6\}, \{1, 2, 4\} = H_2$. We see $G \neq H_1 \cup H_2$. Thus we see the group elements 3 and 5 are left out.

However $G_i$ the identity graph associated with G is

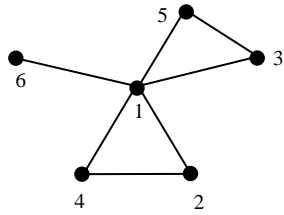

We from our definition we cannot talk of S hence of χ (S).

**Example 2.44:** Let $G = Z_{12}$, the group under addition modulo 12. 0 is the identity element of $Z_{12}$.

The identity graph $G_i$ associated with G is given by the following figure.

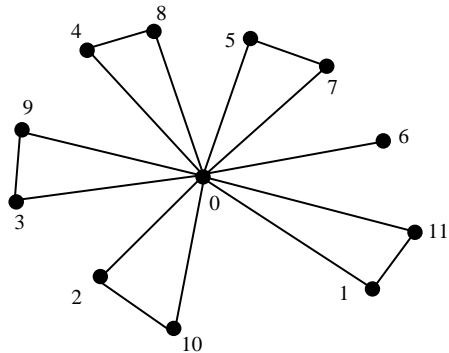

Figure 2.50

The subgroups of $G = Z_{12}$ are $H_1 = \{0, 6\}$, $H_2 = \{2, 4, 6, 8, 10, 0\}$ and $H_3 = \{0, 3, 6, 9\}$.

We see the elements 1 or 5 or 7 or 11 cannot be in any one of the subgroups. So S for $G = Z_{12}$ cannot be formed as union of subgroups.

**Example 2.45:** Let $G = \langle g \mid g^{18} = 1 \rangle$ be the group of order 18, The identity graph $G_i$ of G is given below.

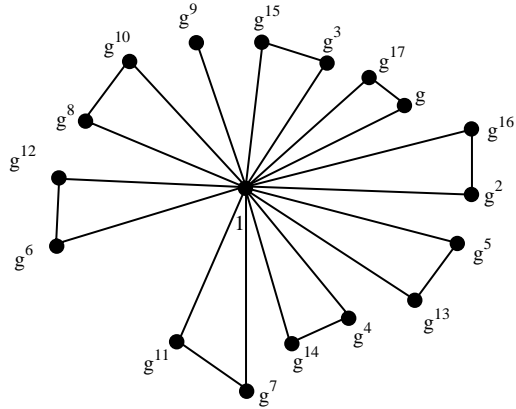

Figure 2.51

The subgroups of G are $H_1 = \{1, g^9\}$, $H_2 = \{1, g^2, g^4, g^6, g^8, g^{10}, g^{12}, g^{14}, g^{16}\}$, $H_3 = \{1, g^3, g^6, g^9, g^{12}, g^{15}\}$, $H_4 = \{1, g^6, g^{12}\}$. Clearly $\{g, g^5, g^7, g^{11}, g^{13}, g^{17}\}$ do not form any part of any of the subgroups. Thus for this G also we cannot find any S. So the question of $\chi(S)$ is impossible.

Now we proceed on to find more examples.

***Example 2.46:*** Let $D_{2.8} = \{a, b \,/\, a^2 = b^8 = 1; bab = a\}$ be the dihedral group of order 16. Thus $D_{28} = \{1, a, b, b^2, b^3, b^4, b^5, b^6, b^7, ab, ab^2, ab^3, ab^4, ab^5, ab^6, ab^7\}$ and its identity graph is as follows.

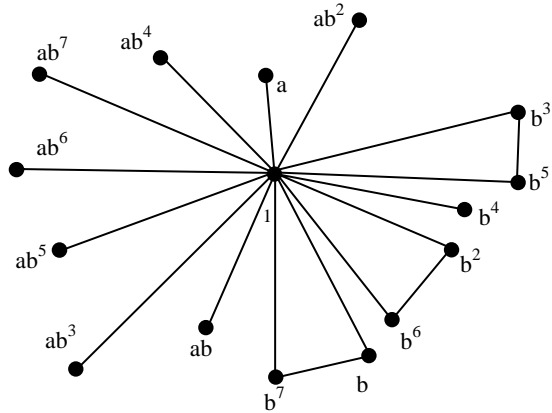

Figure 2.52

The subgroups of $D_{2.6}$ which can contribute to S are as follows.
$$H_1 = \{1, b, b^2, b^3, b^4, b^5, b^6, b^7\}$$
$$H_2 = \{1, a\}, H_3 = \{1, ab\},$$
$$H_4 = \{1, ab^2\}, H_5 = \{1, ab^3\},$$
$$H_6 = \{1, ab^4\}, H_7 = \{1, ab^5\},$$
$$H_8 = \{1, ab^6\} \text{ and } H_9 = \{1, ab^7\}.$$

Clearly $S = \{H_1, H_2, H_3, H_4, H_5, H_6, H_7, H_8, H_9\}$ is such that $D_{2.8} = \bigcup_{i=1}^{9} H_i$. Now $\chi(S) = 3$.

*Example 2.47:* Let

$$A_4 = \left\{ e = \begin{pmatrix} 1 & 2 & 3 & 4 \\ 1 & 2 & 3 & 4 \end{pmatrix}, p_1 = \begin{pmatrix} 1 & 2 & 3 & 4 \\ 2 & 1 & 4 & 3 \end{pmatrix}, \right.$$

$$p_2 = \begin{pmatrix} 1 & 2 & 3 & 4 \\ 3 & 4 & 1 & 2 \end{pmatrix}, p_3 = \begin{pmatrix} 1 & 2 & 3 & 4 \\ 4 & 3 & 2 & 1 \end{pmatrix},$$

$$p_4 = \begin{pmatrix} 1 & 2 & 3 & 4 \\ 1 & 3 & 4 & 2 \end{pmatrix}, p_5 = \begin{pmatrix} 1 & 2 & 3 & 4 \\ 1 & 4 & 2 & 3 \end{pmatrix},$$

$$p_6 = \begin{pmatrix} 1 & 2 & 3 & 4 \\ 3 & 2 & 4 & 1 \end{pmatrix}, p_7 = \begin{pmatrix} 1 & 2 & 3 & 4 \\ 4 & 2 & 1 & 3 \end{pmatrix},$$

$$p_8 = \begin{pmatrix} 1 & 2 & 3 & 4 \\ 2 & 4 & 3 & 1 \end{pmatrix}, p_9 = \begin{pmatrix} 1 & 2 & 3 & 4 \\ 4 & 1 & 3 & 2 \end{pmatrix},$$

$$\left. p_{10} = \begin{pmatrix} 1 & 2 & 3 & 4 \\ 2 & 3 & 1 & 4 \end{pmatrix}, p_{11} = \begin{pmatrix} 1 & 2 & 3 & 4 \\ 3 & 1 & 2 & 4 \end{pmatrix} \right\}$$

be the alternating group of $S_4$. Let $G_i$ be the identity graph associated with $A_4$ which is given below.

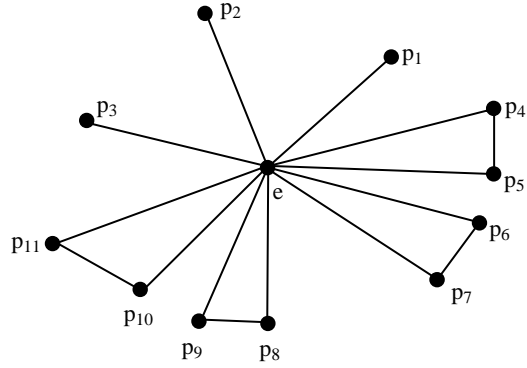

Figure 2.53

$S = \{H_1 = \{e, p_1, p_2, p_3\}, H_2 = \{p_4, p_5, e\}, H_3 = \{e, p_6, p_7\}, H_4 = \{e, p_8, p_9\}, H_5 = \{p_{10}, p_{11}, e\}\}. = \{P_1 = \{e, p_1\}, P_2 = \{e, p_2\}, P_3 = \{e, p_3\}, P_4 = \{e, p_5, p_4\}, P_5 = \{e, p_6, p_7\}, P_6 = \{e, p_8, p_9\}, P_7 = \{e, p_{10}, p_{11}\}\}$.

We see we have 3 colouring for these subgroups.

$$G = \bigcup_{i=1}^{5} H_i = \bigcup_{i=1}^{7} P_i .$$

Now we have seen from the examples some groups have a nice representation i.e., when S is well defined in case of some groups S does not exist. This leads us to define a new notion called graphically good groups and graphically bad groups.

**DEFINITION 2.4:** *Let G be a group of finite or infinite order. S = {$H_1, \ldots, H_n$; subgroups of G such that $H_i \not\subseteq H_j$, $H_j \not\subseteq H_i$ if $1 \le j$, $j \le n$ and $G = \bigcup_{i=1}^{n} H_i$ } i.e., the clique of the group exists and G is colourable then we say G is a graphically good group. If G has no clique then we call G a graphically bad group and the identity subgraph of $G_i$ for its subgroups is called the special bad identity subgraph $G_i$ of G.*

We illustrate by a few examples these situations.

***Example 2.48:*** Let $D_{2.6} = \{a, b \,/\, a^2 = b^6 = 1, bab = a\} = \{1, a, b, b^2, \ldots, b^5, ab, ab^2, \ldots, ab^5\}$ be the dihedral group of order 12. $S = \{H_1 = \{1, a\}, H_2 = \{1, ab\}, H_3 = \{1, ab^2\}, H_4 = \{1, ab^3\}, H_5 =$ $\{1, ab^4\}$, $H_6 = \{1, ab^5\}$, $H_7 = \{1, b, b^2, \ldots, b^5\}\}$ where $G = \bigcup_{i=1}^{7} H_i$  $H_i \not\subset H_j$; $1 \le i, j \le 7\}$. The special identity graph $G_i$ of $D_{2.6}$ is

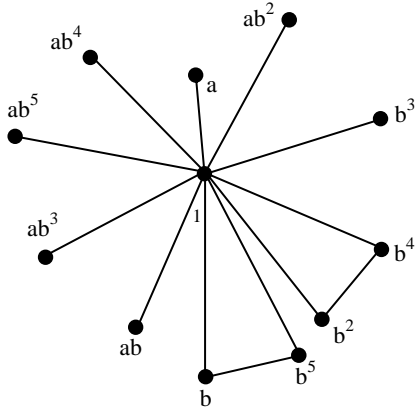

Figure 2.54

Then $\chi(S) = 3$. By putting $b$, $b^2$, $b^3$, $b^4$, $b^5$ one colour say red, $ab^2$ and $ab$–blue colour, $ab^3$ and $ab^4$–red, $ab^3$–blue and $a$–white we can use a minimum of three colours to colour the group $D_{2.6}$. Thus $D_{2.6}$ is graphically a good graph.

***Example 2.49:*** Let $G = \{g \mid g^{12} = 1\}$ be the cyclic group of order 12. The subgroups of $G$ are $H_1 = \{1, g^2, g^4, g^6, g^8 \text{ and } g^{10}\}$, $H_2 = \{1, g^3, g^6, g^9\}$. The elements of $G$ which is not found in the subgroups $H_1$ and $H_2$ are $\{g, g^5, g^7, g^{11}\}$.

The special identity graph of $G$ is given in figure 2.55.

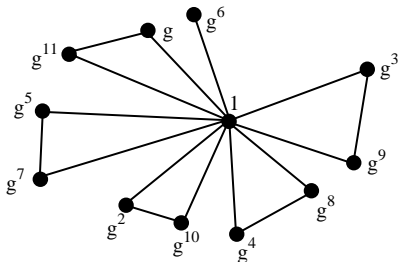

Figure 2.55

The special bad subgroup of G associated with the subgroups $H_1$ and $H_2$ of G is

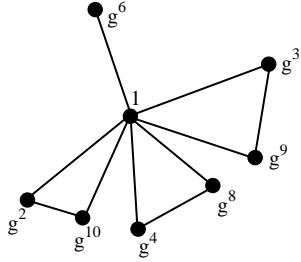

$H_1 = \{1, g^6, g^3, g^9\}$ and $H_2 = \{1, g^2, g^4, g^6, g^8, g^{10}\}$ where $H_1 \cap H_2 = \{1, g^6\}$.

We cannot give any colouring less than two. Now one interesting problem arises $g^6$ should get which colour if we colour the two subgroups of G. Such situations must be addressed to.

This situation we cannot colour so we say the identity graphs are impossible to be coloured so we call the special graphs related with these groups to be graphically not colourable groups.

**Example 2.50:** Let $G = \langle g \mid g^6 = 1 \rangle = \{1, g, g^2, g^3, g^4, g^5\}$ be the cyclic group of order 6.

The special identity graph associated with G is $G_i$ which is given below.

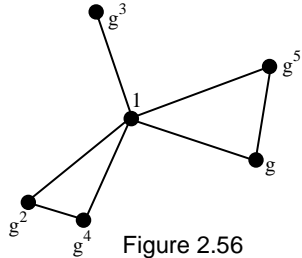

Figure 2.56

The subgroups of G are $H_1 = \{1, g^3\}$ and $H_2 = \{1, g^2, g^4\}$. Clearly $G \neq H_1 \cup H_2$; The bad graph associated with G is given by the following figure.

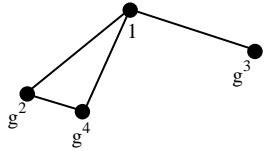

This group needs atleast two colours and the group is graphically bad group or has a bad graph for its subgroup as S does not exist.

***Example 2.51:*** Let $G = \langle g \mid g^8 = 1 \rangle = \{1, g, g^2, g^3, g^4, g^5, g^6, g^7\}$ be the cyclic group of order 8. The special identity graph of G denoted by $G_i$ is as follows.

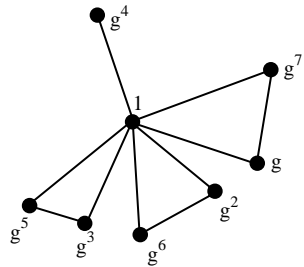

Figure 2.57

The subgroups of G are $H_1 = \{1, g^4\}$, $H_2 = \{1, g^2, g^4, g^6\}$. Thus S does not exist as S has only one subgroup, viz., $H_2$.

So the question of finding a minimal colour does not exist as it has only one subgroup.

Now this is yet a special and interesting study.

***Example 2.52:*** Let $G = \{g \mid g^9 = 1\} = \{1, g, g^2, g^3, g^4, g^5, g^6, g^7, g^8\}$ be the cyclic group of order 9.

The special identity graph of G is given by $G_i$ which is as follows:

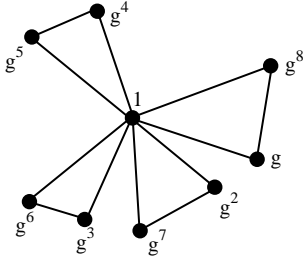

Figure 2.58

Now the subgroup of G is $H_1 = \{1, g^3, g^6\}$. So G cannot have S associated with it and hence the question of colouring G does not arise.

***Example 2.53:*** Let $G = \langle g \mid g^{16} = 1 \rangle = \{1, g, g^2, g^3, \ldots, g^{15}\}$ be the cyclic group of order 16. The identity special graph $G_i$ associated with G is given below.

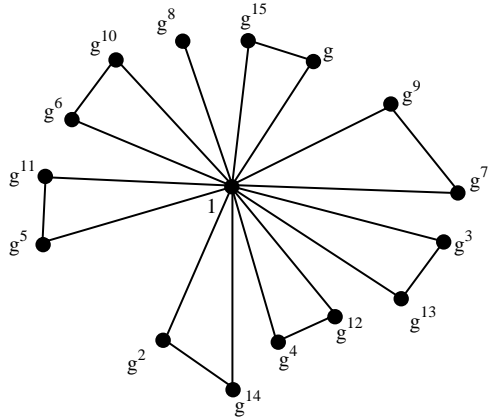

Figure 2.59

The subgroups of G are $H_1 = \{1, g^8\}$, $H_2 = \{1, g^4, g^8, g^{12}\}$ and $H_3 = \{1, g^2, g^4, g^6, g^8, g^{10}, g^{12}, g^{14}\}$. Clearly $H_1 \subseteq H_2 \subseteq H_3$. So $H_3$ is the only subgroup of G and we cannot find S.

Thus colouring of G does not arise. In view of this we define single special identity subgraph of colourable subgroup of a group.

**DEFINITION 2.5:** *Let G be group of finite order. If G has only one subgroup $H_i$ such that $H_i \not\subseteq H_j$ for only j, j ≠ i, i.e., G has one and only one maximal subgroup (i.e., we say a subgroup H of G is maximal if $H \subseteq K \subseteq G$. K any other subgroup containing H then either K = H or K = G). We see then in such case we have to give only one colour to the special identity subgraph of subgroup.*

*We call this situation as a single colourable special identity subgraph, hence a single colourable bad group. Clearly these groups are badly colourable groups. But every bad colourable special identity subgraph in general is not a single special colourable special identity subgraph.*

***Example 2.54:*** Let $G = \langle g \mid g^{25} = 1 \rangle = \{1, g, g^2, \ldots, g^{24}\}$ be the cyclic group of order 25. The subgroups of G is $H = \{g^5, g^{10}, g^{15}, g^{20}, 1\}$; i.e., G has one and only one subgroup. The special identity graph $G_i$ of G is as follows:

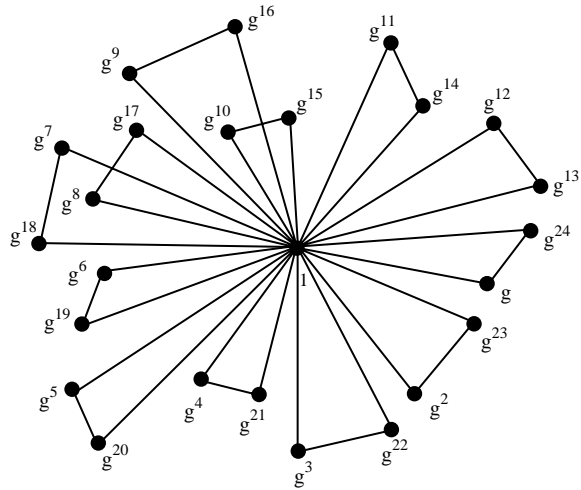

Figure 2.60

The special identity subgraph associated with the subgroup H of G is i.e., single (one) colourable special identity subgraph i.e., G is a single colourable bad group.

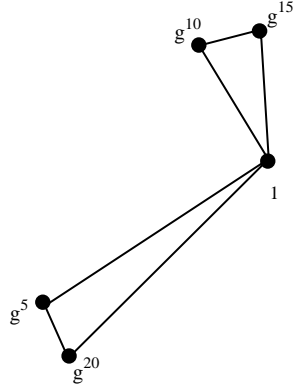

Now we give a theorem which guarantees the existence of single colourable bad groups.

**THEOREM 2.7:** *Let* $G = \langle g \mid g^{p^2} = 1 \rangle$ *where p is a prime, be a cyclic group of order $p^2$. G is a single colourable bad group.*

*Proof:* Let $G = \langle g \mid g^{p^2} = 1 \rangle$ = {1, g, $g^2$, ... $g^{p^2-1}$} be a cyclic group of order $p^2$, p a prime. The only subgroup of G is H = {$g^p$, $g^{2p}$, $g^{3p}$, ..., $g^{(p-1)p}$, 1}. Clearly order of H is p. Further G has no other subgroups. The special graph of G is formed by $(p^2 - 1)/2$ triangles centered around 1. The special identity subgroup of H is formed by $(p - 1) / 2$ triangles centered around 1.

Thus G is a one colourable bad group.

*Example 2.55:* Let $G = \langle g \mid g^{32} = 1 \rangle$ be a cyclic group of order 32. The largest subgroup of G is given by H = {1, $g^2$, $g^4$, $g^6$, $g^8$, ..., $g^{30}$} which is of order 16 and has no other subgroup K such that K $\not\subset$ H.

Thus G is a one colourable graphically bad group.

**THEOREM 2.8:** *Let $G = \langle g \mid g^{p^n} = 1 \rangle$ where p is a prime $n \geq 2$ be a cyclic group of order $p^n$. G is a one colourable graph bad group.*

*Proof:* Given G is a cyclic group of prime power order. To show G has only one subgraph H such that there is no other subgroup K in G with $K \subseteq H$ or $H \subseteq K$. i.e., G has one and only one maximal subgroup.

Now the maximal group $H = \{1, g^p, g^{2p}, \ldots, g^{(n-1)p}\}$ which is of order $p^{n-1}$.

Thus G is only one colourable as g has no S associated with it we see it is a bad graph group. Now as G has only one subgroup G is a uniquely colourable graph bad group. Thus we have a class of groups which are single colourable graph bad groups.

*Example 2.56:* Let $G = \{g \mid g^{10} = 1\}$ be a cyclic group of order 10. The subgroups of G are $H_1 = \{g^2, g^4, g^6, g^8, 1\}$ and $H_2 = \{1, g^5\}$. Clearly S does not exist as $G \neq H_1 \cup H_2$ so G cannot be a good graph group.

The special identity graph associated with G is as follows

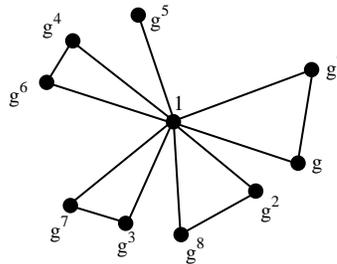

Figure 2.61

The subgraph of the subgroups of G is given by

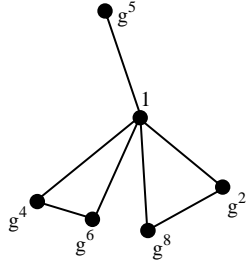

Clearly $H_i$ is minimum two colourable. For $g^5$ vertex is given one colour and the vertices $g^2$, $g^4$, $g^6$ and $g^8$ are given another colour. However $G_i \neq H_i$ as the vertices $\{g, g^9, g^3, g^7\}$ cannot be associated with a group.

***Example 2.57:*** Let $G = \{g \mid g^{30} = 1\}$ be a cyclic group of order 30. The special identity graph $G_i$ associated with G is as follows.

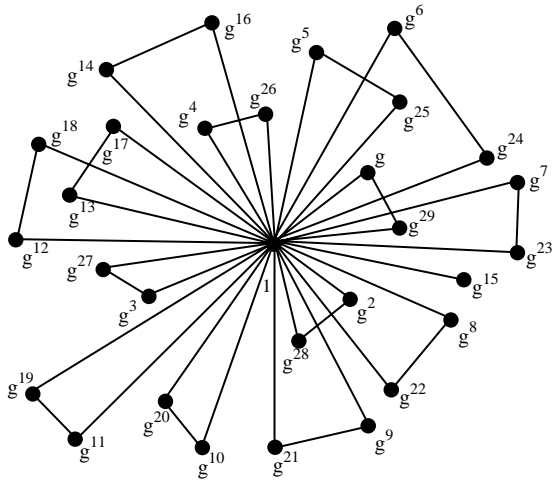

Figure 2.62

The subgroups of G are $H_1 = \{1, g^{15}\}$, $H_2 = \{g^{10}, g^{20}, 1\}$, $H_3 = \{1, g^5, g^{10}, g^{15}, g^{20}, g^{25}\}$, $H_4 = \{1, g^6, g^{12}, g^{18}, g^{24}\}$, $H_5 = \{1, g^2, g^4, g^6, g^8, \ldots, g^{26}, g^{28}\}$ and $H_6 = \{1, g^3, g^6, g^9, g^{12}, g^{15}, g^{18}, g^{21}, g^{24}, g^{27}\}$.

Clearly $G \neq \cup H_i$.

The maximal subgroups are $H_3$, $H_5$ and $H_6$. We cannot give one colour to $g^6$, $g^{12}$, $g^{10}$, $g^{18}$, $g^{20}$, $g^{24}$.

Hence it is impossible to colour the subgroups properly. In view of this we are not in a position to say whether there exists three colourable bad groups.

*Example 2.58:* Let $G = \{g \mid g^{18} = 1\}$ be a cyclic group of order 18. The subgroups of G are $H_1 = \{1, g^9\}$, $H_2 = \{g^3, g^6, g^9, g^{12}, g^{15}, 1\}$, $H_3 = \{g^2, g^4, g^6, g^8, g^{10}, g^{12}, g^{14}, g^{16}, 1\}$ and $H_4 = \{1, g^6, g^{12}\}$. We have two subgroups $H_2$ and $H_3$ such that $H_2 \not\subseteq H_3$ or $H_3 \not\subseteq H_2$. However $H_2 \cap H_3 = \{1, g^6, g^{12}\}$. So how to colour these subgroups even as bad groups. The special identity graph $G_i$ associated with G is as follows:

Figure 2.63

Since we have common elements between $H_2$ and $H_3$. We cannot assign any colour to $g^6$ and $g^{12}$. Hence this is not a colourable bad group.

Now we proceed onto show that we have a class of two colourable bad groups.

**THEOREM 2.9:** *Let $G = \langle g \mid g^n = 1 \rangle$ where $n = pq$ with p and q two distinct primes be a cyclic group of order n. Then G is a two colourable bad group.*

*Proof:* Let $G = \{1, g, g^2, ..., g^{n-1}\}$ be the cyclic group of degree n. The maximal of G are $H_1 = \{1, g^p, g^{2p}, ..., g^{(p-1)q}\}$ and $H_2 = \{1, g^q, g^{2q}, ..., g^{(q-1)p}\}$ be subgroups of G.

Clearly both $H_1$ and $H_2$ are maximal subgroups of G. However $G \neq H_1 \cup H_2$ and $H_1 \cap H_2 = \{e\}$. Thus G is a two colourable bad group.

*Example 2.59:* Let $G = \{g \mid g^{35} = 1\}$ be a cyclic group of order 35. The two maximal subgroups of G are $H_1 = \{1, g^5, g^{10}, g^{15}, g^{20}, g^{25}, g^{30}\}$ and $H_2 = \{1, g^7, g^{14}, g^{21}, g^{28}\}$. $H_1 \cap H_2 = \{1\}$ and $G \neq H_1 \cup H_2$. These can be coloured with 2 colours.

The special identity subgraph associated with the subgroups is given by 5 triangles with centre one which is as follows:

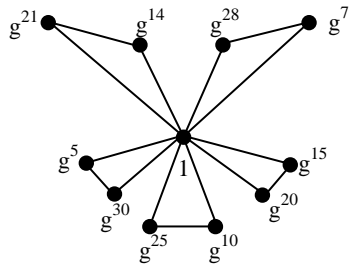

Figure 2.64

The vertices $g^5$, $g^{30}$, $g^{25}$, $g^{10}$, $g^{20}$ and $g^{15}$ are given one colour and $g^{21}$, $g^{14}$, $g^{27}$ and $g^7$ are given another colour. Thus G is two colourable special graph bad group.

*Example 2.60:* Let $G = \langle g \mid g^{22} = 1 \rangle$ be the cyclic group of order 22. The two maximal subgroups of G are $H_1 = \{1, g^{11}\}$ and $H_2 = \{g^2, g^4, g^6, g^8, g^{10}, g^{12}, g^{14}, g^{16}, g^{18}, g^{20}, 1\}$

Thus G is a two colourable special identity subgraph bad group.

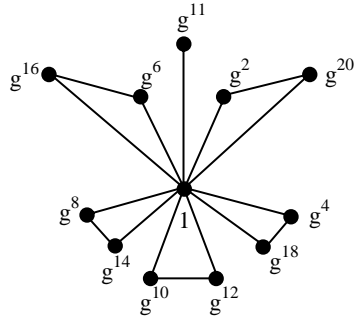

Figure 2.65

The vertex $\{g^{11}\}$ is given one colour and the vertices $\{g^2, g^4, g^6, g^8, g^{10}, g^{12}, g^{14}, g^{16}, g^{18}$ and $g^{20}\}$ are given another colour.

Now we proceed onto study the colouring problem of normal subgroups of a group G.

**DEFINITION 2.6:** *Let G be a group. Suppose $N = \{N_1, \ldots, N_n\}$ are normal subgroups of G such that $G = \bigcup_i N_i$, $N_i \cap N_j = \{e\}$ then we call N the normal clique of G and has atmost n elements, then clique $G = n$. If the sizes of the clique are not bounded, then define normal clique $G = \infty$.*

*The chromatic number of the special identity graph $G_i$ of G denoted by $\chi(G)$, is the minimum k for which $N_1, \ldots, N_n$ accepts k colours where one colour is given to the vertices of a subgroup $N_i$ of G for $i = 1, 2, \ldots, n$.*

*We call such groups are k – colourable normal good groups.*

However the authors find it as an open problems to find a groups G which can be written as a union of normal subgroups $H_i$ such $G = \bigcup_i H_i$ with $H_i \cap H_j = \{1\}$, 1 the identity element of G.

We give illustration before we define more notions in this direction.

***Example 2.61:*** Let $G = \{g \mid g^6 = 1\} = \{1, g, g^2, g^3, g^4, g^5\}$ be a cyclic group of order six. The normal subgroups of G are $H_1 =$ $\{1, g^3\}$ and $H_2 = \{1, g^2, g^4\}$. Clearly the elements of $H_1$ i.e., $g^3$ can be given one colour and the vertices of $H_2$ viz. $\{1, g^2, g^4\}$ can be given another colour thus two colours are sufficient to colour the normal subgroups of G. However the elements $\{g, g^5\}$ do not form any part of the normal subgroups so cannot be given any colour.

***Example 2.62:*** Let $G = \{g \mid g^{12} = 1\} = \{1, g, g^2, g^3, g^4, g^5, g^6, g^7, g^8, g^9, g^{10}, g^{11}\}$ be a group of order 12. The normal subgroups of G are $H_1 = \{1, g^2, \ldots, g^{10}\}$ and $H_2 = \{1, g^3, g^6, g^9\}$ other normal subgroups are $H_3 = \{1, g^6\}$, $H_4 = \{1, g^4, g^8\}$. We see $H_3$ and $H_4$ are contained in $H_1$ and $H_3 \subseteq H_2$ and further $H_1 \cap H_2 = \{1, g^6\}$, so the normal subgroups cannot be coloured as $g^6$ cannot simultaneously get two colours.

In view of all these examples we propose the following definition.

**DEFINITION 2.7:** *Let G be a group, if G has no normal subgroups $N_i$ such that $G = \bigcup_{i=1}^{n} N_i$ and $N_i \cap N_j = \{1\}$, if $i \neq j$ then we say G is a k-colourable normal bad group $0 \leq k \leq n$. If $k = 1$ then we say G is a one colourable normal bad group. If $k = 2$ then we say G is a 2-colourable normal bad group. If $k = t$ then we say G is a t-colourable normal bad group. If $k = 0$ we say G is a 0-colourable normal bad group.*

We illustrate these situations by explicit examples.

***Example 2.63:*** $G = \{g \mid g^p = 1\}$, p a prime be a cyclic group of order p. We see G is simple i.e., G has no normal subgroups, so G is 0-colourable normal bad group.

***Example 2.64:*** Let $A_3$ be the alternating subgroup of $S_3$. $A_3$ is also 0-colourable normal bad group.

***Example 2.65:*** Let $A_n$, $n \geq 5$ be the alternating subgroup of $S_n$. $A_n$ is also 0-colourable normal bad group.

**THEOREM 2.10:** The class of 0-colourable normal bad groups is non empty.

*Proof:* Take all alternating subgroups $A_n$ of $S_n$, $n \neq 4$, $n \geq 3$; and as $A_n$ is simple so $A_n$'s are 0-colourable normal bad groups. Similarly all groups G of order p, p a prime is a 0-colourable normal bad group as G has no proper subgroups. Hence the theorem.

*Example 2.66:* Let $S_5$ be the symmetric group of order $\underline{|5}$. $S_5$ has only one proper normal subgroup viz. $A_5$ so $S_5$ is a 1-colourable normal bad group.

*Example 2.67:* Let $G = \{g \mid g^9 = 1\}$ be a cyclic group of order 9. $G = \{1, g, g^2, g^3, g^4, g^5, g^6, g^7, g^8\}$ and $H = \{1, g^3, g^6\}$ is the only normal subgroup of G. So G is one colourable normal bad group.

*Example 2.68:* Let $G = \{g \,/\, g^{32} = 1\} = \{1, g, g^2, \ldots, g^{31}\}$ be the cyclic group of order 32. The only normal subgroup of G is $H = \{1, g^2, g^4, \ldots, g^{30}\}$. Clearly every other normal subgroup K of G is properly contained in H. Thus G is one colourable normal bad group.

Inview of these examples we have the following theorem which gurantees the existence of a class of one colourable normal bad groups.

**THEOREM 2.11:** *Let* $G = \left\langle g \mid g^{p^n} = 1 \right\rangle$ *where p is any prime and* $n \geq 2$, *G is a one colourable normal bad group.*

*Proof:* To show G is a one colourable normal bad group, it is enough if we show G has only one normal subgroup H and all other subgroups are contained in H. Take $H = \{1, g^p, g^{p^2}, \ldots, g^{p^{n-1}}\} = \langle g^p \rangle$; This is the largest normal subgroup of G. Clearly any other subgroup of G is contained in H. Thus G has only one normal subgroup such that all other subgroups of G are contained in it. Hence G is a one colourable normal bad group, as all the vertices of the group H are given one colour i.e., all the elements of H are given only one colour.

Thus G is a one colourable normal bad group.

**THEOREM 2.12:** *Let $S_n$ be the symmetric group of degree n. $S_n$ is a one colourable normal bad group.*

*Proof:* Let $S_n$ be the symmetric group of degree n. $A_n$ be the alternating normal subgroup of $S_n$. We know $S_n$ has only one normal subgroup viz. $A_n$. So $S_n$ is a one colourable normal bad group.

Thus the class of one colourable bad groups is non empty.

Consider the following example.

***Example 2.69:*** $A_4$ be the alternating subgroup of $S_4$. $A_4$ is one colourable normal bad group.

***Example 2.70:*** Let $G = \{g \mid g^{26} = 1\}$ be a cyclic group of order 26. G has only two normal subgroups $H_1 = \{1, g^{13}\}$ and $H_2 = \{1, g^2, g^4, ..., g^{24}\}$.

Infact G has no other subgroups. But $G \neq H_1 \cup H_2$ so G is not a k-colourable normal good group, but only a 2-colourable normal bad group.

This theorem shows the existence of 2-colourable normal bad groups.

**THEOREM 2.13:** *Let $G = \{g \mid g^{pq} = 1\}$ where p and q are primes $p \neq q$, be a cyclic group of order pq. G is a 2-colourable normal bad group.*

*Proof:* Given $G = \{1, g, g^2, ..., g^{pq-1}\}$ be the cyclic group of order pq where p and q are primes; $p \neq q$. The two normal subgroups of G are $H_1 = \{1, g^p, p^{2p}, ..., g^{p(q-1)}\}$ and $H_2 = \{1, g^q, p^{2q}, ..., g^{q(p-1)}\}$. Thus the group G is a 2-colourable normal bad group.

***Example 2.71:*** Let $G = Z_n = \{0, 1, 2, \ldots, pqr - 1 = n - 1\}$, (where p, q and r are distinct primes) be the group of order n. Take $n = 30$, $Z_{30} = \{0, 1, 2, \ldots, 29\}$. The normal subgroups of $Z_{30}$ are $H_1 = \{0, 2, 4, 6, 8, 10, 12, 14, 16, 18, 20, 22, 24, 26, 28\}$, $H_2 = \{0, 3, 6, 9, 12, 15, 18, 21, 24, 27\}$ and $H_3 = \{0, 5, 10, 15, 20, 25\}$. Thus $Z_{30}$ has only 3 normal subgroups. G is not a three colourable normal bad group.

***Example 2.72:*** $Z_{210} = \{0, 1, 2, \ldots, 209\}$ is an additive group order 210. $H_1 = \langle 2 \rangle$, $H_2 = \langle 3 \rangle$, $H_3 = \langle 5 \rangle$ and $H_4 = \langle 7 \rangle$ are the four normal subgroups of $Z_{210}$. $Z_{210}$ is a not a four colourable normal bad group.

***Example 2.73:*** Let $G = Z_3 \times Z_2 \times Z_5 = \{ (a, b, c) / a \in Z_3, b \in Z_2$ and $c \in Z_5\}$ a group got as a the external direct product of the groups under modulo addition. G is a group of order 30. G has three normal subgroups. $H_1 = \{(000), (100) (200)\}$, $H_2 = \{(000), (010)\}$ and $H_3 = \{(000), (001), (002), (003), (004)\}$.

These are the normal subgroups of G. However $G \neq \bigcup_{i=1}^{3} H_i$ $H_i \cap H_j = (000)$ if $i \neq j$.

The special identity subgraph of G is as follows.

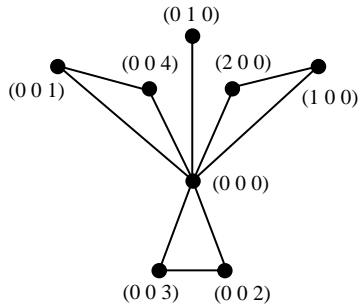

Figure 2.66

Clearly minimum 3 colours are needed to colour this figure. Thus one colour is given to the vertex (010) another colour to the vertices {(001), (002), (003) and (004)}. Yet another colour different from other two are given to the vertices (100) and (200). (000) can be given any one of the three colours.

Thus G is a 3-colourable normal bad groups.

**Example 2.74:** Let $G = G_1 \times G_2 \times G_3 \times G_4$ where $G_1 = \{g \mid g^3 = 1\}$, $G_2 = \{Z_2\}$, $G_3 = Z_4$ under addition modulo 4 and $G_4 = \{g \mid g^5 = 1\}$ is the external direct product of 4 distinct groups. G is of order 60.

The normal subgroups of G are $H_1 = \{(1001), (g\,001)\,(g^2 001)\}$ $H_2 = \{(1001)\,(1101)\}$, $H_3 = \{(1001)\,(1011)\,(1021)\,(1031)\}$ and $H_4 = \{(1001)\,(100g)\,(100g^2)\,(100g^3)\,(100g^4)\}$ We see $H_i \cap H_j = (1001)$ for all $i \neq j$; $1 \leq i, j \leq 4$.

However $G \neq \bigcup_{i=1}^{4} H_i$. Now G is two colourable normal bad group. For the elements of $H_1$ can be given one colour. $H_2$ and $H_3$ another colour and $H_4$ the colour as that of $H_1$. The corresponding special identity subgraph is given by the following diagram.

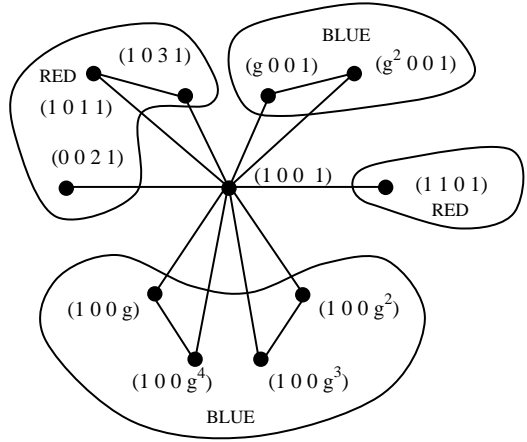

Figure 2.67

(1001) can be given red or blue colourable. Thus G is a 2 colourable normal bad group.

**THEOREM 2.14:** *Let $G = G_1 \times G_2 \times \ldots \times G_n$ be the direct product of n-groups, n an even number. Then G is a 2 colourable normal bad group.*

*Proof:* Given $G = G_1 \times G_2 \times \ldots \times G_n$ direct product of n-distinct / different groups n even. We see $\overline{G_i} = \{e_1\} \times \{e_2\} \times \ldots \times G_i \times \ldots \times \{e_n\}$ is a subgroup in G, which is normal in G and is isomorphic with $G_i$. Here $\{e_i\}$ is the identity element of $G_i$, $i = 1, 2, \ldots, n$.

Thus G has n normal subgroups which are disjoint.

Clearly $\overline{G_i} \cap \overline{G_j} = \{$identity element of G$\}$. Now since n is even we see this group G is a 2 colourable normal bad group.

It is left as a simple problem for the reader to find the number of colours needed to colour the group $G = G_1 \times \ldots \times G_n$, when n is odd, where $G_1, \ldots, G_n$ are n distinct groups in the external direct product.

Can one say n be even or odd in $G = G_1 \times \ldots \times G_n$. G is 2 colourable normal bad group.

In case of k-colourable bad group or k-colourable group of k-colourable normal bad group we see k is always less than or equal to 3. Thus more than study of these k-colourable property the interesting features would be their graphic representation for it would interest the students to look at it and concretely view atleast finite groups.

Next we see about the p-sylow subgroup of a group and their colouring.

We want to show by colouring how many p-sylow subgroups are conjugate etc. This is mainly to attract the students about groups and their properties. So we are not bothered about colouring with minimum number of colours but we are bothered about mainly how we can make easy the understanding of Sylow theorem or Cauchy theorem and so on.

First we show how to colour the p-sylow subgroups of a group G.

*Example 2.75:* Let $S_3 = \{1, p_1, p_2, p_3, p_4, p_5\}$ be the symmetric group of degree 3. The order of $S_3$ is 6. $6 = 2.3$. $S_3$ has only 2-sylow subgroups and 3-sylow subgroups. The 2-sylow subgroups of $S_3$ are $H_1 = \{1, p_1\}$, $H_2 = \{1, p_2\}$ and $H_3 = \{1, p_3\}$. The 3-sylow subgroup of $S_3$ is $H_3 = \{1, p_4, p_5\}$.

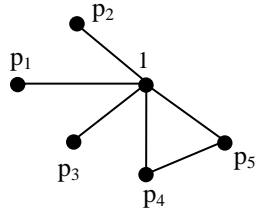

Figure 2.68

$p_1$, $p_2$, $p_3$ are given one colour for they are 2-Sylow subgroups of $S_3$. The 3-sylow subgroup $H_3$ is given another colour to $p_4$ and $p_5$. Thus the vertices are given different colours for only different types of p-sylow subgroups.

Thus $S_3$ is 2 colourable p-sylow subgroups.

*Example 2.76:* Let $S_4$ be the permutation group of degree 4. Clearly $o(S_4) = 1.2.3.4 = 24 = 2^3.3$. Thus $S_4$ has only 2-sylow subgroups and 3-sylow subgroups. Hence $S_4$ is a 2 colourable p-sylow subgroups.

*Example 2.77:* Let $S_5$ be the permutation group of (12345). $o(S_5) = \underline{|5} = 1.2.3.4.5 = 2^3. 3.5$. $S_5$ has only 2-sylow subgroups, 3-sylow subgroups and 5-sylow subgroups. Thus $S_5$ is a 3-colourable p-sylow subgroup.

In view of the above examples we have the following theorem.

**THEOREM 2.15:** *Let $S_n$ be the permutation group of degree n. The group $S_n$ is a m-colourable p- sylow subgroup, where m is the number of primes less than or equal to n.*

*Proof:* Let $S_n$ be the symmetric group of degree n. $o(S_n) = 1.2.3…n = \underline{|n}$. Now p is a prime such that $p \leq n$, then $p / o(S_n)$ as $o(S_n) = \underline{|n}$. So $S_n$ has p-sylow subgroups. This is true of all primes p such that $p \leq n$ as $o(S_n) = \underline{|n}$. So if m is the number of distinct primes which divide $\underline{|n}$ or equivalently m is the number of primes which is less than or equal to n then we have m-distinct sylow subgroups of different orders. Thus the group $S_n$ is a m-colourable p-sylow subgroups.

*Example 2.78:* Consider $S_{31}$, the symmetric group of degree 31. $o(S_{31}) = \underline{|31} = 1.2.3.4.5.6.7.8 \ldots 31$. The primes which are less than or equal to 31 are 2,3, 5, 7, 11, 13, 17, 19, 23, 29 and 31. Thus $S_{31}$ has 11-distinct sylow subgroups of different orders. Thus $S_{31}$ is a 11-colourable p-sylow subgroups.

Another interesting way of colouring the vertices of a group is by colouring the center of a non commutative group G. So by looking at the colours the student can know the centre of the group G.

*Example 2.79:* Let $G = \{a, b \;/\; a^2 = b^5 = 1; bab = a\} = \{1, a, b, ab, ab^2, ab^3, ab^4, b^2b^3, b^4\}$ C (G) = 1 only 1 is given a colour and rest of the vertices in the graph $G_i$ of G remain uncoloured.

Thus the graph of G is

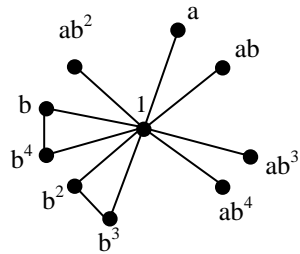

Figure 2.69

*Example 2.80:* Let $G = \{a, b \;/\; a^2 = b^4 = 1, bab = a\} = \{1, a, b, b^2, b^3, ab, ab^2, ab^3\}$. The graph of G is as follows.

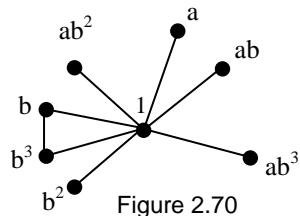

Figure 2.70

In this graph 1 and $b^2$ are given a colour, rest of the vertices remain uncoloured.

***Example 2.81:*** Let $G = S_3 \times \{g \,/\, g^4 = 1\}$ be the group of order 24. To find the center of G. Given $G = \{(e\ 1), (p_1, 1), (p_2, 1), (p_3,1), (p_4,1), (p_5,1), (e, g), (p_1, g), (p_2, g), (p_3, g), (p_4, g), (p_5, g), (e,g^2), (p_1,g^2), (p_2,g^2), (p_3,g^2), (p_4,g^2), (p_5,g^2), (p_1,g^3), (e_1,g^3), (p_2,g^3), (p_3,g^3), (p_4,g^3), (p_5,g^3)\}$. The center of G is given by $\{(e, 1), (e, g), (e, g^2), (e, g^3)\}$

The graph of G is given below

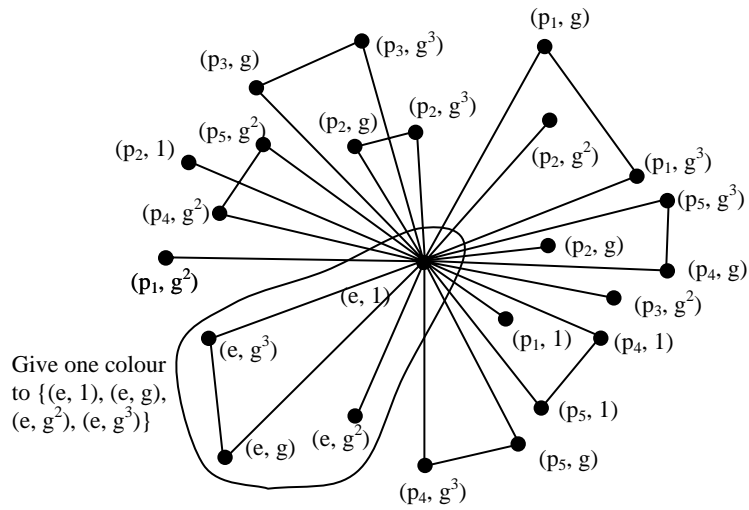

Figure 2.71

The rest of the vertices are not given any colour.

So by looking at the coloured graph $G_i$ one can find the center and inverse of each elements.

***Example 2.82:*** Let $G = S_3 \times D_{2.11} \times H$ where $H = \{g \mid g^{14} = e_1\}$ be a group of order $n = 6 \times 22 \times 14$. The center of G is given by $C(G) = \{(e, 1, g) \,/\, g \in H\} \subseteq G$.

The graph $G_i$ associated with G is coloured as follows:

Out of the 1848 vertices of the group G only the 14 vertices of $C(G) = \{e, 1, e_1) (e, 1, g), …, (e, 1, g^{13})\}$ are given one colour, rest are uncoloured.

Now all the vertices which form the group elements are adjoined with identity. Two elements are adjacent if and only if one is the inverse of the other and they are joined. If an element is a such that its square is identity that is self inversed then only the vertex is joined with the identity in the centre. All the group elements which form the vertex set are joined with the identity.

**Remark:** By looking at the identity graph of a group one can immediately say the number elements x in the group G which are such that $x^2 = e$ (e-identity element of G)

Now we proceed onto give the matrix representation of the identity graph of G which is called as the graph – matrix representation of a group G.

The identity graph in the case of the group G can be given the corresponding adjacency matrix or connection matrix which is known as the graph –matrix representation of the group.

**DEFINITION 2.8:** *Let G be a group with elements e, $g_1$, …, $g_n$. Clearly the order of G is n + 1. Let $G_i$ be the identity graph of G. The adjacency matrix of $G_i$ is a $(n + 1) \times (n + 1)$ matrix $X = (x_{ij})$ in which diagonal terms are zero i.e., $x_{ii} = 0$ for i = 1, 2, …, n + 1 the first row and first column are one except, $x_{ij} = 1$ if the element $g_i$ is the inverse of $g_j$ in which case $x_{ij} = x_{ji} = 1$ if $(i \neq j)$. We call the matrix $X = (x_{ij})$ to be $(n + 1) \times (n + 1)$ the identity graph matrix of the group G.*

We shall illustrate this situation by some examples.

*Example 2.82:* Let $Z_{10} = \{0, 1, 2, …, 9\}$ be the group under addition modulo 10. The identity graph of $Z_{10}$ is given by the following figure.

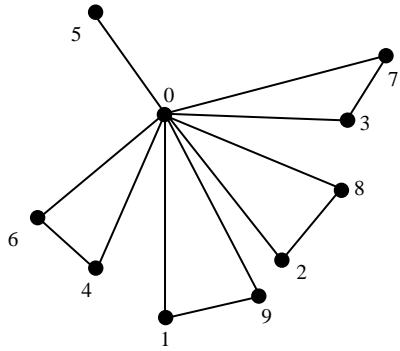

Figure 2.72

The identity graph matrix of $Z_{10}$ is a $10 \times 10$ matrix X.

$$X = \begin{array}{c|cccccccccc} & 0 & 1 & 2 & 3 & 4 & 5 & 6 & 7 & 8 & 9 \\ \hline 0 & 0 & 1 & 1 & 1 & 1 & 1 & 1 & 1 & 1 & 1 \\ 1 & 1 & 0 & 0 & 0 & 0 & 0 & 0 & 0 & 0 & 1 \\ 2 & 1 & 0 & 0 & 0 & 0 & 0 & 0 & 0 & 1 & 0 \\ 3 & 1 & 0 & 0 & 0 & 0 & 0 & 0 & 1 & 0 & 0 \\ 4 & 1 & 0 & 0 & 0 & 0 & 0 & 1 & 0 & 0 & 0 \\ 5 & 1 & 0 & 0 & 0 & 0 & 0 & 0 & 0 & 0 & 0 \\ 6 & 1 & 0 & 0 & 0 & 1 & 0 & 0 & 0 & 0 & 0 \\ 7 & 1 & 0 & 0 & 1 & 0 & 0 & 0 & 0 & 0 & 0 \\ 8 & 1 & 0 & 1 & 0 & 0 & 0 & 0 & 0 & 0 & 0 \\ 9 & 1 & 1 & 0 & 0 & 0 & 0 & 0 & 0 & 0 & 0 \end{array}$$

***Example 2.84:*** Let $G = \{a, b \mid a^2 = b^4 = 1, bab = a\} = \{1, a, b, b^2, b^3, ab, ab^2, ab^3\}$ be the dihedral group of order 8.

The identity graph of G denoted by $G_i$ is as follows.

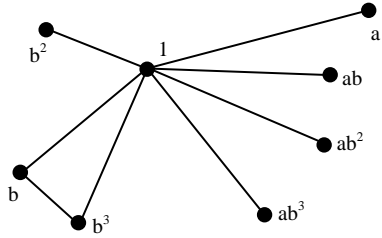

Figure 2.73

The identity graph matrix X of $G_i$ is as follows.

|       | 1 | a | b | $b^2$ | $b^3$ | ab | $ab^2$ | $ab^3$ |
|-------|---|---|---|-------|-------|----|--------|--------|
| 1     | 0 | 1 | 1 | 1     | 1     | 1  | 1      | 1      |
| a     | 1 | 0 | 0 | 0     | 0     | 0  | 0      | 0      |
| b     | 1 | 0 | 0 | 0     | 1     | 0  | 0      | 0      |
| $b^2$ | 1 | 0 | 0 | 0     | 0     | 0  | 0      | 0      |
| $b^3$ | 1 | 0 | 1 | 0     | 0     | 0  | 0      | 0      |
| ab    | 1 | 0 | 0 | 0     | 0     | 0  | 0      | 0      |
| $ab^2$| 1 | 0 | 0 | 0     | 0     | 0  | 0      | 0      |
| $ab^3$| 1 | 0 | 0 | 0     | 0     | 0  | 0      | 0      |

*Example 2.85:* Let

$$A_4 = \{1, g_1 = \begin{pmatrix} 1 & 2 & 3 & 4 \\ 2 & 1 & 4 & 3 \end{pmatrix}, g_2 = \begin{pmatrix} 1 & 2 & 3 & 4 \\ 3 & 4 & 1 & 2 \end{pmatrix},$$

$$g_3 = \begin{pmatrix} 1 & 2 & 3 & 4 \\ 4 & 3 & 2 & 1 \end{pmatrix}, g_4 = \begin{pmatrix} 1 & 2 & 3 & 4 \\ 1 & 3 & 4 & 2 \end{pmatrix},$$

$$g_5 = \begin{pmatrix} 1 & 2 & 3 & 4 \\ 1 & 4 & 2 & 3 \end{pmatrix}, g_6 = \begin{pmatrix} 1 & 2 & 3 & 4 \\ 3 & 2 & 4 & 1 \end{pmatrix},$$

$$g_7 = \begin{pmatrix} 1 & 2 & 3 & 4 \\ 4 & 2 & 1 & 3 \end{pmatrix}, g_8 = \begin{pmatrix} 1 & 2 & 3 & 4 \\ 2 & 4 & 3 & 1 \end{pmatrix},$$

$$g_9 = \begin{pmatrix} 1 & 2 & 3 & 4 \\ 4 & 1 & 3 & 2 \end{pmatrix}, g_{10} = \begin{pmatrix} 1 & 2 & 3 & 4 \\ 2 & 3 & 1 & 4 \end{pmatrix},$$

$$g_{11} = \begin{pmatrix} 1 & 2 & 3 & 4 \\ 3 & 1 & 2 & 4 \end{pmatrix} \Big\}$$

be the alternating subgroup of $S_4$.
The graph $G_i$ of $A_4$ is as follows.

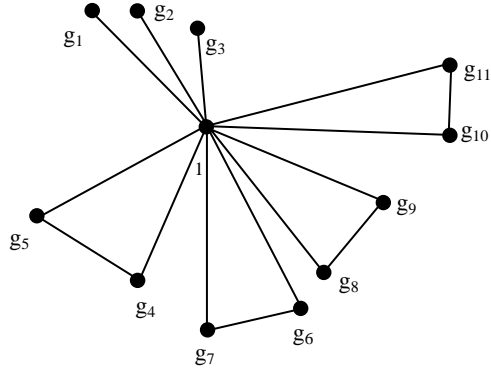

Figure 2.74

The corresponding identity graph – matrix X of $A_4$ is as follows:

|     | 1 | $g_1$ | $g_2$ | $g_3$ | $g_4$ | $g_5$ | $g_6$ | $g_7$ | $g_8$ | $g_9$ | $g_{10}$ | $g_{11}$ |
|---|---|---|---|---|---|---|---|---|---|---|---|---|
| 1 | 0 | 1 | 1 | 1 | 1 | 1 | 1 | 1 | 1 | 1 | 1 | 1 |
| $g_1$ | 1 | 0 | 0 | 0 | 0 | 0 | 0 | 0 | 0 | 0 | 0 | 0 |
| $g_2$ | 1 | 0 | 0 | 0 | 0 | 0 | 0 | 0 | 0 | 0 | 0 | 0 |
| $g_3$ | 1 | 0 | 0 | 0 | 0 | 0 | 0 | 0 | 0 | 0 | 0 | 0 |
| $g_4$ | 1 | 0 | 0 | 0 | 0 | 1 | 0 | 0 | 0 | 0 | 0 | 0 |
| $g_5$ | 1 | 0 | 0 | 0 | 1 | 0 | 0 | 0 | 0 | 0 | 0 | 0 |
| $g_6$ | 1 | 0 | 0 | 0 | 0 | 0 | 0 | 1 | 0 | 0 | 0 | 0 |
| $g_7$ | 1 | 0 | 0 | 0 | 0 | 0 | 1 | 0 | 0 | 0 | 0 | 0 |
| $g_8$ | 1 | 0 | 0 | 0 | 0 | 0 | 0 | 0 | 0 | 1 | 0 | 0 |
| $g_9$ | 1 | 0 | 0 | 0 | 0 | 0 | 0 | 0 | 1 | 0 | 0 | 0 |
| $g_{10}$ | 1 | 0 | 0 | 0 | 0 | 0 | 0 | 0 | 0 | 0 | 0 | 1 |
| $g_{11}$ | 1 | 0 | 0 | 0 | 0 | 0 | 0 | 0 | 0 | 0 | 1 | 0 |

(1) We see X is a symmetric matrix with diagonal entries to be zero.
(2) Further if the row $g_i$ has only one 1 at the 1$^{st}$ column then $g_i \in G$ is such that $g_i^2 = 1$.
(3) If a row $g_j$ has two ones and the rest zero then for the $g_j$ we have a $g_k$ row that has two ones and $g_j g_k = g_k g_j = 1$.

This observation is also true for columns.

We now proceed onto study or define a graph of a group by its conjugate elements. We assume the groups are non commutative.

**DEFINITION 2.9:** *Let G be a non abelian group. The equivalence classes of G be denoted by [e], [$g_1$], ..., [$g_n$]. Then each element $h_i$ in an equivalence class [$g_i$], is joined with $g_i$, i = 1, 2, ..., n.*

*This graph will be known as conjugate graph of the conjugacy classes of a non commutative group.*

We illustrate this situation by the following examples.

***Example 2.86:*** Let $S_3 = \{e, p_1, p_2, p_3, p_4, p_5\}$ be the symmetric group of degree 3. The conjugacy classes of $S_3$ are [e], [$p_1$] and [$p_4$]. The conjugacy graph of $S_3$ is

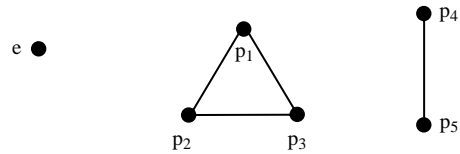

Figure 2.75

***Example 2.87:*** Let $D_{24} = \{a, b \,/\, a^2 = b^4 = 1, bab = a\} = \{1, a, b, b^2, b^3, ab, ab^2, ab^3\}$. The conjugacy class of $D_{24}$ are $\{1\}$, $\{a, ab^2, b^2\} = \{a\}$, $\{b\} = \{b, b^3\}$ and $\{ab^3\} = \{ab, ab^3\}$.

The conjugacy graph of $D_{24}$ is

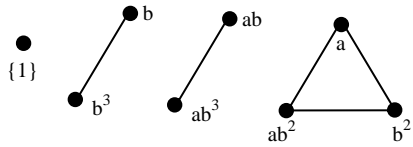

Figure 2.76

We now proceed onto find the conjugacy graph of $D_{27}$.

***Example 2.88:*** Let $G = D_{27} = \{a, b / a^2 = b^7 = 1, bab = a\} = \{1, a, b, b^2, \ldots, b^6, ab, ab^2, \ldots, ab^6\}$ be the dihedral group of order 14.

The conjugacy classes of $D_{27}$ are $\{1\}$, $\{a\} = \{a, ab^3, ab, ab^5, ab^4, ab^6, ab^2\}$, $\{b\} = \{b, b^6\}$, $\{b^2\} = \{b^2, b^5\}$ and $\{b^3\} = \{b^3, b^4\}$.

The conjugacy graph associated with $D_{27}$ is

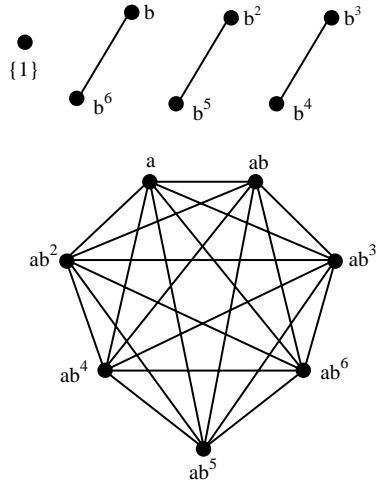

Figure 2.77

We see the graphs of different type for $D_{2n}$ when n is a prime and n a non prime.

***Example 2.89:*** Let $D_{26} = \{a, b \mid a^2 = b^6 = 1, bab = a\} = \{1, a, b, b^2, \ldots, b^5, ab, ab^2, \ldots, ab^5\}$ be a group of order 12. The conjugacy classes of $D_{2.6}$ is $\{1\}$, $\{b^2\} = \{b^2, b^4\}$, $\{a\} = \{a, ab^2,$ $ab^4\}$, $\{b^3\} = \{b^3\}$, $\{b\} = \{b, b^5\}$, $\{ab\} = \{ab, ab^3, ab^5\}$. The conjugacy graph of $D_{26}$ is

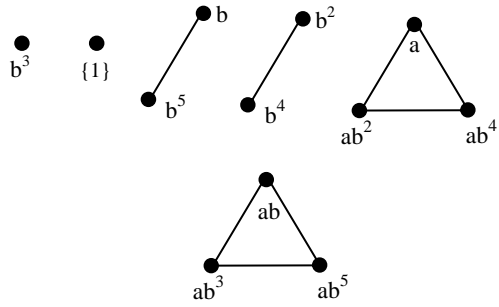

Figure 2.78

***Example 2.90:*** Let $D_{2.9} = \{a, b \mid a^2 = b^9 = 1, bab = a\}$ be the group of order 18. The conjugacy classes of $D_{29}$ is $\{1\}$, $\{a\} = \{a, ab^2, ab^4, ab^6, ab^8, ab, ab^3, ab^5, ab^7\}$, $\{b\} = \{b, b^8\}$ $b^2 = \{b^2, b^7\}$ $\{b^3\} = \{b^3, b^6\}$, $\{b^4\} = \{b^4, b^5\}$. The conjugacy graph of $D_{2.9}$ is

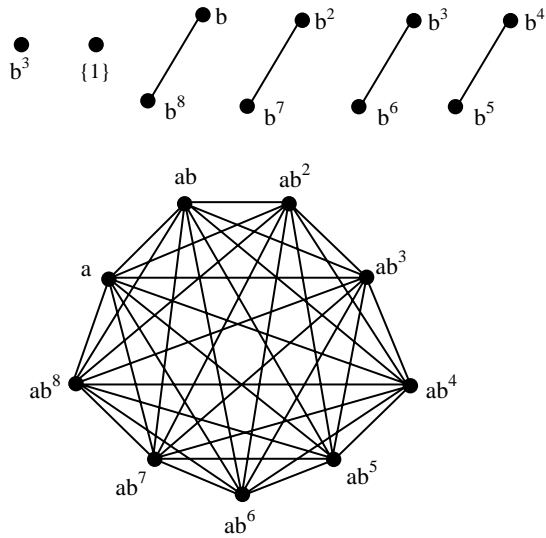

Figure 2.79

We have the following theorem.

**THEOREM 2.16:** *Let $Z_{2p} = \{a, b \,/\, a^2 = b^p = 1, bab = a\}$ be a dihedral group of order 2p, p a prime. The conjugacy classes of $Z_{2p}$ forms a collection of complete graph with $p - 1\,/\,2$ complete graphs with two vertices and one complete graph with p vertices.*

*Proof:* Let $Z_{2p} = \{a, b \,/\, a^2 = b^p = 1\, bab = a\} = \{1, a, b, b^2, \ldots, b^{p-1}, ab, ab^2, \ldots, ab^{p-1}\}$ where p is a prime. The conjugacy classes of $Z_{2p}$ are $\{1\}$, $\{a\} = \{a, ab, ab^2, \ldots, ab^{p-1}\}$, $\{b\} = \{b, b^{p-1}\}$, $\{b^2\} = \{b^2, b^{p-2}\}$ … $\left\{b^{\frac{p-1}{2}}\right\} = \left\{b^{\frac{p-1}{2}}, b^{\frac{p+1}{2}}\right\}$.

Clearly the conjugacy graph associated with this groups consists of a point graph, $p - 1/2$ number complete graphs with two vertices and one complete graph with p vertices, which is indicated below

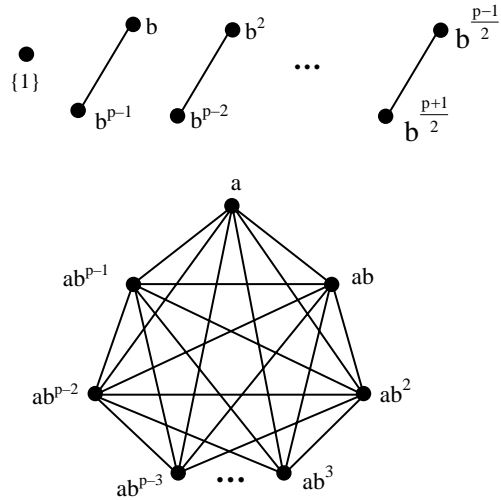

Figure 2.80

***Example 2.91:*** Let $Z_{2.10} = \{a, b \,/\, a^2 = b^{10} = 1, bab = a\}$ be the dihedral group of order 20.

The conjugacy classes of $Z_{2.10}$ are $\{1\} = \{1\}$, $\{a\} = \{a, ab^2, ab^4, ab^6, ab^8\}$, $\{ab^3\} = \{ab^7, ab^3, ab, ab^5, ab^9\}$, $\{b\} = \{b, b^9\}$, $\{b^2\} = \{b^2, b^8\}$, $\{b^3\} = \{b^3, b^7\}$, $\{b^4\} = \{b^4, b^6\}$.

The conjugacy graph of $D_{20}$ is given below.

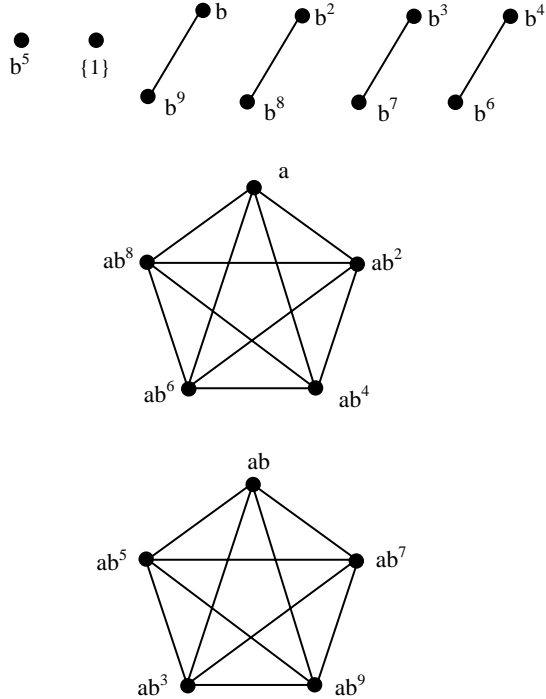

Figure 2.81

We see in this case we have 2 complete graphs with five vertices.

*Example 2.92:* Consider the dihedral group of order 24 given by $D_{2.12} = \{a, b \ / \ a^2 = b^{12} = 1, bab = a\}$. The conjugacy classes of $D_{2.12}$ is as follows.

$\{1\}$, $\{a\} = \{a, ab^2, ab^4, ab^6, ab^8, ab^{10}\}$, $\{ab\} = \{ab, ab^9, ab^7, ab^5, ab^3, ab^{11}\}$ $\{b\} = \{b, b^{11}\}$, $\{b^2\} = \{b^2, b^{10}\}$, $\{b^3\} = \{b^3, b^9\}$, $\{b^4\} = \{b^4, b^8\}$, $\{b^5\} = \{b^5, b^7\}$ and $\{b^6\}$. The conjugacy graph of $D_{2.12}$ is as follows.

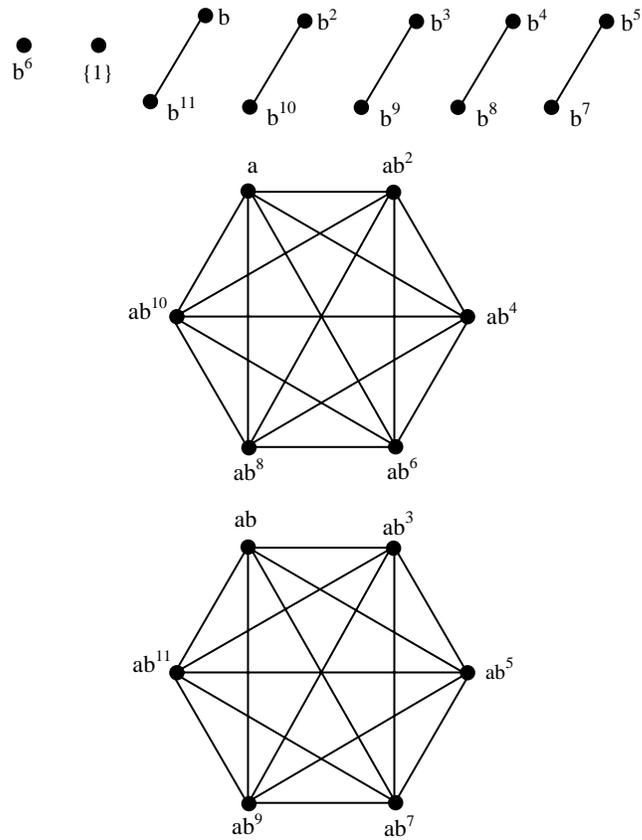

Figure 2.82

We see the conjugacy graph of $D_{2.12}$ contains 2 points graphs 5 complete graph with two vertices and 2 complete graphs with 6 vertices.

***Example 2.93:*** Let $Z_{2.18} = \{a, b \mid a^2 = b^{18} = 1, bab = a\}$ be the dihedral group of order 36. The conjugacy classes of $Z_{2.18}$ are $\{1\}$, $\{a\} = \{a, ab^2, ab^4, ab^6, ab^8, ab^{10}, ab^{12}, ab^{14}, ab^{18}\}$, $\{ab\} = \{ab, ab^3, ab^5, ab^7, ab^9, ab^{11}, ab^{13}, ab^{15}, ab^{17}\}$ $\{b\} = \{b, b^{17}\}$, $\{b^2\} = \{b^2, b^{16}\}$, $b^3 = \{b^3, b^{15}\}$, $b^4 = \{b^4, b^{14}\}$, $b^5 = \{b^5, b^{13}\}$, $\{b^6\} = \{b^6, b^{12}\}$, $\{b^7\} = \{b^7, b^{11}\}$, $\{b^8\} = \{b^8, b^{10}\}$ and $\{b^9\}$. The conjugacy graph of $Z_{2.18}$ is as follows.

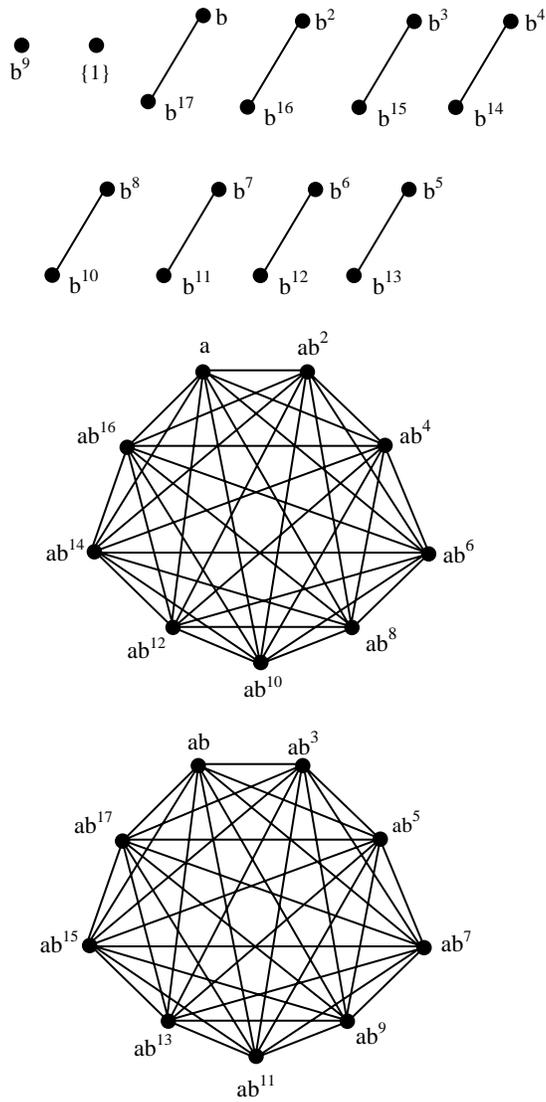

Figure 2.83

In view of this we have the following theorem

**THEOREM 2.17:** *Let $D_{2n} = \{a, b \,/\, a^2 = b^n = 1, bab = a\}$ where n is of the form 2r. Then the conjugacy graph of $D_{2n}$ has two point graphs. (n – 2)/2 number of complete graph with two vertices and two complete graphs with r vertices.*

*Proof:* Given $D_{2n} = \{a, b \,/\, a^2 = b^n = 1\; bab = a\}$ where n = 2r is a dihedral group of order 2n. The conjugacy classes of $D_{2n}$ are $\{1\}, \{r\}, \{b, b^{2r-1}\}, \{b^2, b^{2r-2}\}, \{b^3, b^{2r-3}\}, \ldots, \{b^{r-1}, b^{2r-(r-1)}\}$,
    $\{a, ab^2, ab^4, \ldots, ab^{2r-2}\}$ and
    $\{a, ab^3, ab^5, \ldots, ab^{2r-1}\}$.
    The conjugacy graph of $D_{2n}$ consists of 2r-2/2 number of complete 2 vertices graphs and two complete graph with r vertices. Hence the claim.

*Example 2.94:* Let $D_{2.8} = \{a, b \,/\, a^2 = b^8 = 1\; bab = a\}$ be the dihedral group of order 16. The conjugacy classes of $D_{28}$ are $\{1\}, \{a, ab^2, ab^4, ab^6\}, \{b, b^7\}, \{b^2, b^6\}, \{b^3, b^5\}, \{b^4\} \{ab, ab^5, ab^3, ab^7\}$. The conjugacy graph associated with $D_{28}$ is as follows.

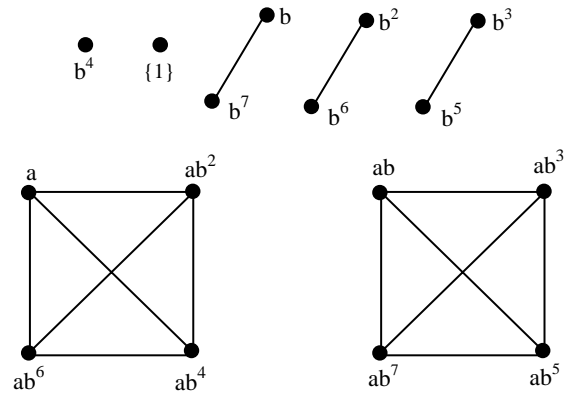

Figure 2.84

*Example 2.95:*

$$A_4 = \left\{\begin{pmatrix}1 & 2 & 3 & 4\\1 & 2 & 3 & 4\end{pmatrix}, \begin{pmatrix}1 & 2 & 3 & 4\\2 & 1 & 4 & 3\end{pmatrix}, \begin{pmatrix}1 & 2 & 3 & 4\\4 & 3 & 2 & 1\end{pmatrix}, \right.$$

$$\begin{pmatrix} 1 & 2 & 3 & 4 \\ 3 & 4 & 1 & 2 \end{pmatrix}, \begin{pmatrix} 1 & 2 & 3 & 4 \\ 1 & 3 & 4 & 2 \end{pmatrix}, \begin{pmatrix} 1 & 2 & 3 & 4 \\ 1 & 4 & 2 & 3 \end{pmatrix},$$

$$\begin{pmatrix} 1 & 2 & 3 & 4 \\ 3 & 2 & 4 & 1 \end{pmatrix}, \begin{pmatrix} 1 & 2 & 3 & 4 \\ 4 & 2 & 1 & 3 \end{pmatrix}, \begin{pmatrix} 1 & 2 & 3 & 4 \\ 2 & 4 & 3 & 1 \end{pmatrix},$$

$$\begin{pmatrix} 1 & 2 & 3 & 4 \\ 4 & 1 & 3 & 2 \end{pmatrix}, \begin{pmatrix} 1 & 2 & 3 & 4 \\ 2 & 3 & 1 & 4 \end{pmatrix}, \begin{pmatrix} 1 & 2 & 3 & 4 \\ 3 & 1 & 2 & 4 \end{pmatrix} \}$$

be the alternating group of $S_4$.

The conjugacy classes of $A_4$ are

$$\left\{ \begin{pmatrix} 1 & 2 & 3 & 4 \\ 1 & 2 & 3 & 4 \end{pmatrix}, \begin{pmatrix} 1 & 2 & 3 & 4 \\ 2 & 1 & 4 & 3 \end{pmatrix}, \begin{pmatrix} 1 & 2 & 3 & 4 \\ 4 & 3 & 2 & 1 \end{pmatrix}, \begin{pmatrix} 1 & 2 & 3 & 4 \\ 3 & 4 & 1 & 2 \end{pmatrix} \right\},$$

$$\left\{ \begin{pmatrix} 1 & 2 & 3 & 4 \\ 1 & 3 & 4 & 2 \end{pmatrix}, \begin{pmatrix} 1 & 2 & 3 & 4 \\ 3 & 1 & 2 & 4 \end{pmatrix}, \begin{pmatrix} 1 & 2 & 3 & 4 \\ 4 & 2 & 1 & 3 \end{pmatrix}, \begin{pmatrix} 1 & 2 & 3 & 4 \\ 2 & 4 & 3 & 1 \end{pmatrix} \right\},$$

$$\left\{ \begin{pmatrix} 1 & 2 & 3 & 4 \\ 1 & 4 & 2 & 3 \end{pmatrix}, \begin{pmatrix} 1 & 2 & 3 & 4 \\ 3 & 2 & 4 & 1 \end{pmatrix}, \begin{pmatrix} 1 & 2 & 3 & 4 \\ 2 & 3 & 1 & 4 \end{pmatrix}, \begin{pmatrix} 1 & 2 & 3 & 4 \\ 4 & 1 & 3 & 2 \end{pmatrix} \right\}$$

The conjugacy graph associated with $A_4$ is as follows.

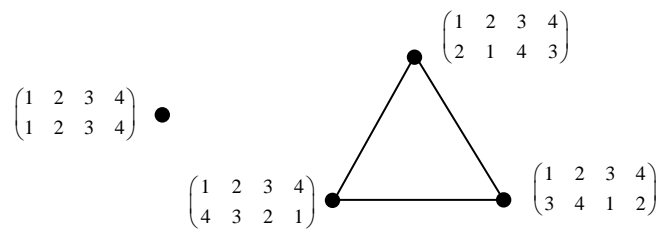

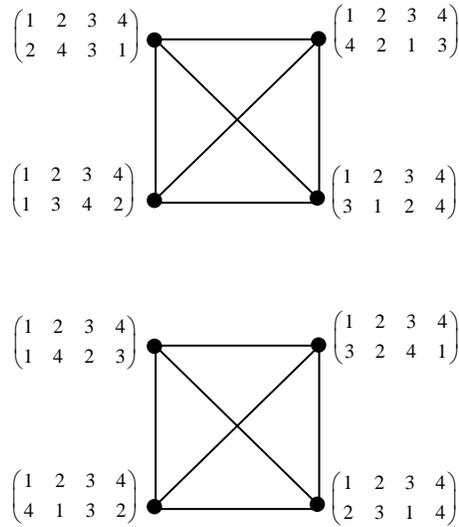

Figure 2.85

**THEOREM 2.18:** *Let G be any non commutative group of finite order. The conjugacy graph of G is always a collection of complete graphs.*

*Proof:* If x is conjugate with p elements say $g_1, \ldots, g_p$ then each $g_i$ is conjugate with p elements resulting in a complete graph with p + 1 vertices. Hence the claim.

**Chapter Three**

# IDENTITY GRAPHS OF SOME ALGEBRAIC STRUCTURES

For the first time we introduce the identity graphs of some algebraic structures like semigroups, loops and commutative rings. This chapter has three sections. In section one we study the identity graph of semigroups and S-semigroups. In section two we study the graphs of loops and commutative groupoids. In the final section the identity graph of a commutative ring is studied.

## 3.1 Identity graphs of semigroups

Now we consider the identity graph of the semigroup S which is taken under multiplication. Let (S, *) be a commutative semigroup with identity 1, i.e., a monoid, we say an element x ∈ S has an inverse y in S if x * y = y * x = 1. If y = x then x * x = $x^2$ = 1 we say x ∈ S is a self inversed element of S.

*Example 3.1.1:* $Z_{12}$ the set of modulo integers 12 is a semigroup under multiplication modulo 12. We see $Z_{12}$ is a commutative monoid.

We can draw identity graph for semigroups with units. We say the element and its inverse are adjoined by an edge. Like zero divisor graphs for semigroups we draw identity graphs for semigroups which are commutative with unit.

*Example 3.1.2:* Let $Z_6 = \{0, 1, 2, 3, 4, 5\}$ be the semigroup under multiplication modulo 6. The identity graph of $Z_6$ is just a line graph joining 1 and 5.

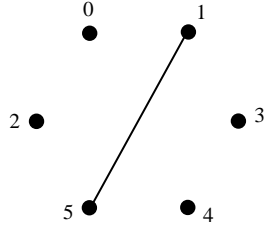

Figure 3.1.1

*Example 3.1.3:* Let $Z_5 = \{0, 1, 2, 3, 4\}$ be the semigroup under multiplication modulo 5. The identity graph is

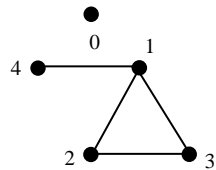

Figure 3.1.2

*Example 3.1.4:* Let $Z_{12} = \{0, 1, 2, \ldots, 11\}$ be the semigroup under multiplication modulo 12. The identity graph of $Z_{12}$ is

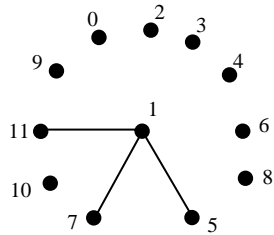

Figure 3.1.3

***Example 3.1.5:*** Let
$$S(2) = \left\{ \begin{pmatrix} 1 & 2 \\ 1 & 1 \end{pmatrix}, \begin{pmatrix} 1 & 2 \\ 1 & 2 \end{pmatrix}, \begin{pmatrix} 1 & 2 \\ 2 & 1 \end{pmatrix}, \begin{pmatrix} 1 & 2 \\ 2 & 2 \end{pmatrix} \right\}$$
be the symmetric semigroup under composition of maps. The identity graph of $S(2)$ is

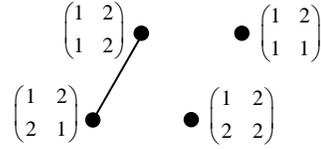

Figure 3.1.4

***Example 3.1.6:*** Let $Z_{15} = \{0, 1, 2, …, 14\}$ be the semigroup under multiplication modulo 15. The identity graph of $Z_{15}$ is

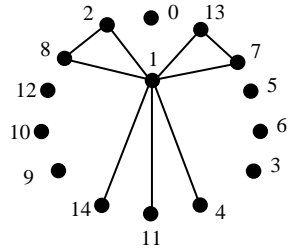

Figure 3.1.5

***Example 3.1.7:*** Let $Z_{14} = \{0, 1, 2, …, 13\}$ be the semigroup under multiplication modulo 14. The identity graph of $Z_{14}$ is

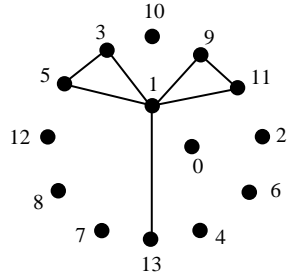

Figure 3.1.6

***Example 3.1.8:*** Let $Z_{10} = \{0, 1, 2, \ldots, 9\}$ be the semigroup under multiplication modulo 10. The identity graph of $Z_{10}$ is

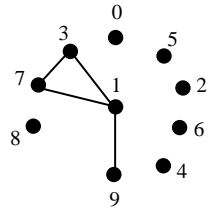

Figure 3.1.7

***Example 3.1.9:*** Let $Z_{21} = \{0, 1, 2, \ldots, 20\}$ be the semigroup under multiplication modulo 21. The identity graph of $Z_{21}$ is

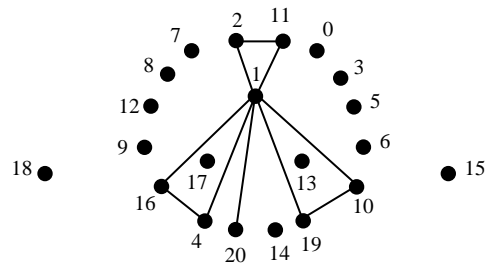

Figure 3.1.8

***Example 3.1.10:*** Let $Z_{18} = \{0, 1, 2, \ldots, 17\}$ be the semigroup under multiplication modulo 18. The identity graph of $Z_{18}$ is as follows:

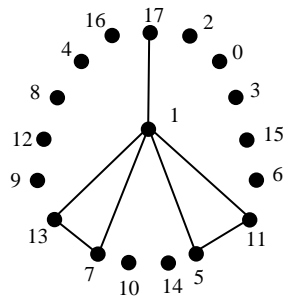

Figure 3.1.9

***Example 3.1.11:*** Let $Z_{30}$ be the semigroup under multiplication modulo 30. The identity graph of $Z_{30}$ is

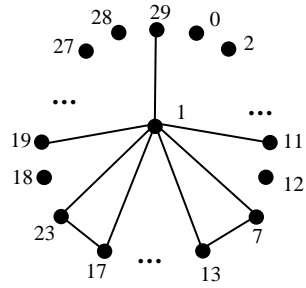

Figure 3.1.10

***Example 3.1.12:*** Let

$$V = \left\{ \begin{pmatrix} 1 & 0 \\ 0 & 0 \end{pmatrix}, \begin{pmatrix} 0 & 0 \\ 0 & 0 \end{pmatrix}, \begin{pmatrix} 0 & 1 \\ 0 & 0 \end{pmatrix}, \begin{pmatrix} 0 & 0 \\ 1 & 0 \end{pmatrix}, \right.$$

$$\begin{pmatrix} 0 & 0 \\ 0 & 1 \end{pmatrix}, \begin{pmatrix} 1 & 1 \\ 0 & 0 \end{pmatrix}, \begin{pmatrix} 1 & 0 \\ 1 & 0 \end{pmatrix}, \begin{pmatrix} 0 & 0 \\ 1 & 1 \end{pmatrix}, \begin{pmatrix} 0 & 1 \\ 0 & 1 \end{pmatrix},$$

$$\begin{pmatrix} 1 & 0 \\ 0 & 1 \end{pmatrix}, \begin{pmatrix} 0 & 1 \\ 1 & 0 \end{pmatrix}, \begin{pmatrix} 1 & 1 \\ 1 & 0 \end{pmatrix},$$

$$\left. \begin{pmatrix} 0 & 1 \\ 1 & 1 \end{pmatrix}, \begin{pmatrix} 1 & 0 \\ 1 & 1 \end{pmatrix}, \begin{pmatrix} 1 & 1 \\ 0 & 1 \end{pmatrix}, \begin{pmatrix} 1 & 1 \\ 1 & 1 \end{pmatrix} \right\}$$

be the semigroup under multiplication; elements of V are from $Z_2 = \{0, 1\}$. Clearly V is a semigroup of order 16.

The identity graph associated with V is as follows.

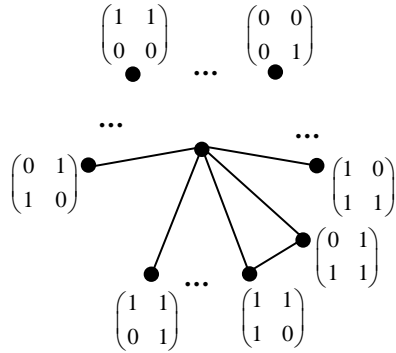
Figure 3.1.11

However the semigroup V is non commutative since
$$ab = ba = \begin{pmatrix} 1 & 0 \\ 0 & 1 \end{pmatrix}$$
is true for the inverse elements we need not be bothered about commutativity, as only thing we should guarantee is that the presence of unique inverse for $a \in V$ which is such that
$$ab = ba = \begin{pmatrix} 1 & 0 \\ 0 & 1 \end{pmatrix}.$$

However this work is left for the reader to prove.

We know the zero divisor graph of the semigroup has been studied extensively in [7]. We now give some examples of the zero divisor graph and compare it with the identity graph.

*Example 3.1.13:* Let $Z_6 = \{0, 1, 2, 3, 4, 5\}$ be the semigroup under multiplication modulo 6. The zero divisors in $Z_6$ are $2.3 \equiv 0 \pmod 6$, $4.3 \equiv 0 \pmod 6$. The zero divisor graph of $Z_6$ is as follows:

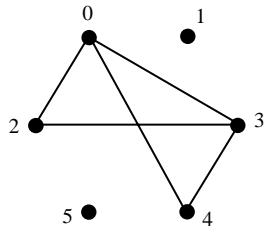
Figure 3.1.12

***Example 3.1.14:*** Let $Z_8 = \{0, 1, 2, 3, \ldots, 7\}$ be the semigroup under multiplication modulo 8. The zero divisor graph of $Z_8$ is as follows:

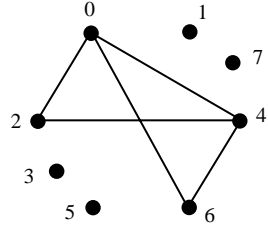

Figure 3.1.13

***Example 3.1.15:*** Let $Z_{10} = \{0, 1, 2, 3, 4, 5, \ldots, 8\}$ be the semigroup under multiplication modulo 10. The zero divisor graph of $Z_{10}$ is as follows:

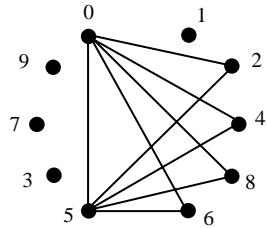

Figure 3.1.14

***Example 3.1.16:*** Let $Z_{18} = \{0, 1, 2, 3, \ldots, 17\}$ be the semigroup under multiplication modulo 18. The zero divisor graph of $Z_{18}$ is as follows:

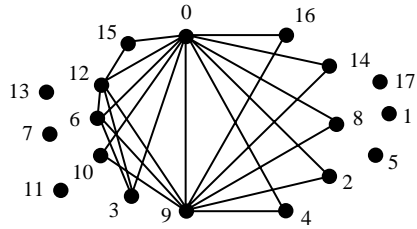

Figure 3.1.15

Now having seen the zero divisor graph and the identity graph of a semigroup S, we now proceed on to define the notion of identity-zero combined graph of a semigroup G.

**DEFINITION 3.1.1:** *Let $S = \{0, 1, s_1, \ldots, s_n\}$ be the commutative semigroup where 0 and $1 \in S$ i.e., semigroup is a monoid which has both zero divisors and units. The combined graph of S with zero divisor graph and identity graph will be known as the combined identity-zero graph of the semigroup.*

We illustrate this by some examples.

***Example 3.1.17:*** Let $Z_6 = \{0, 1, 2, \ldots, 5\}$ be the semigroup under multiplication modulo 6. The identity-zero combined graph (combined identity-zero graph) of $Z_6$ is as follows:

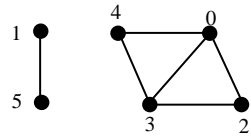

Figure 3.1.16

***Example 3.1.18:*** Let $S = \{0, 1, 2, \ldots, 7\} = Z_8$ be the semigroup under multiplication modulo 8. The combined identity-zero graph of $Z_8$ is as follows:

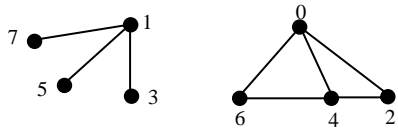

Figure 3.1.17

***Example 3.1.19:*** Let $Z_9 = \{0, 1, 2, \ldots, 8\}$ be the semigroup under multiplication modulo 9. The identity-zero combined graph of $Z_9$ is as follows:

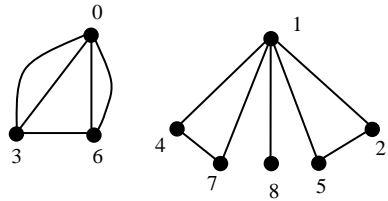

Figure 3.1.18

***Example 3.1.20:*** The combined identity-zero graph of $Z_{10} = \{0, 1, 2, 3, \ldots, 9\}$ under multiplication modulo 10 is as follows:

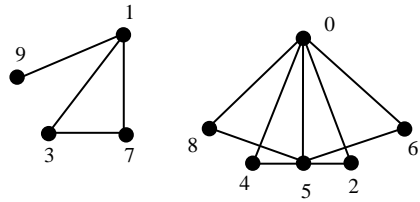

Figure 3.1.19

***Example 3.1.21:*** The identity-zero combined graph of the semigroup $Z_{12} = \{0, 1, 2, 3, \ldots, 11\}$ under multiplication modulo 12 is as follows:

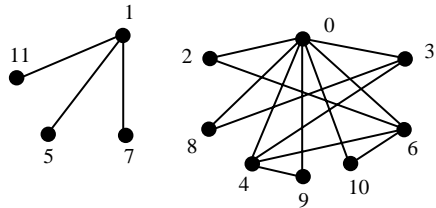

Figure 3.1.20

***Example 3.1.22:*** Let $Z_{15} = \{0, 1, 2, 3, \ldots, 14\}$ be the semigroup under multiplication modulo 15. The combined identity-zero graph of $Z_{15}$ is as follows:

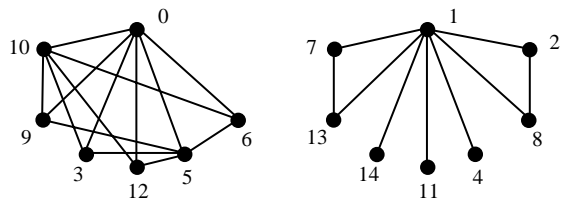

Figure 3.1.21

***Example 3.1.23:*** Let $Z_{16} = \{0, 1, 2, 3, \ldots, 15\}$ be the semigroup under multiplication modulo 16. The combined identity-zero graph of $Z_{16}$ is as follows:

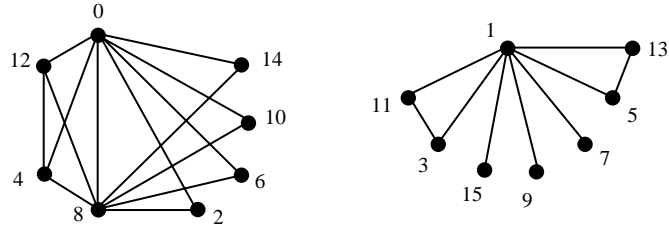

Figure 3.1.22

***Example 3.1.24:*** Let $Z_{30} = \{0, 1, 2, 3, \ldots, 29\}$ be the semigroup under multiplication modulo 30. The combined identity-zero graph of $Z_{30}$ is as follows.

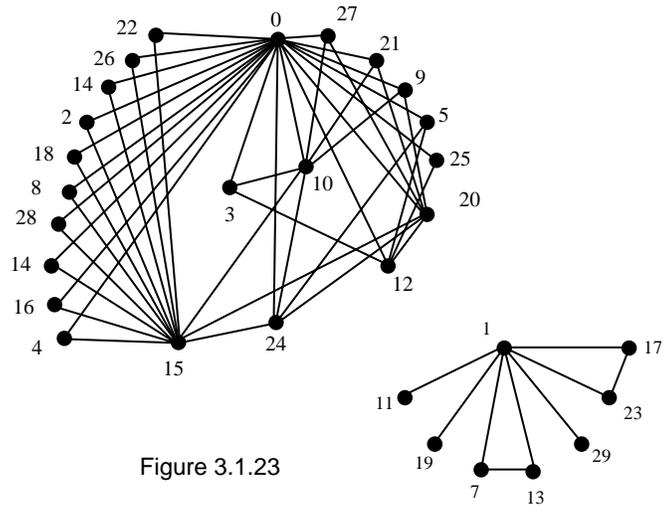

Figure 3.1.23

Now having seen some examples of monoids with zero divisors. We illustrate the three adjacency matrices associated with the identity graph, zero divisor graph and the combined identity-zero graph.

**DEFINITION 3.1.2:** *Let S be a semigroup. The adjacency matrix associated with identity-zero combined graph $S_i$ is defined to be the identity-zero combined adjacency matrix of $S_i$.*

It is assumed that the reader is familiar with adjacency matrix of the identity graph and the zero graph. If the semigroup has no zero divisors they we will not have the notion of zero graph or the combined identity-zero graph.

If the semigroup has no identity then the semigroup has no identity graph associated with it.

We will illustrate this situation before we proceed to define some more new notions.

***Example 3.1.25:*** Let $Z_{12} = \{0, 1, 2, 3, \ldots, 11\}$ be the semigroup under multiplication modulo 12.
The combined identity-zero graph of $Z_{12}$ is as follows.

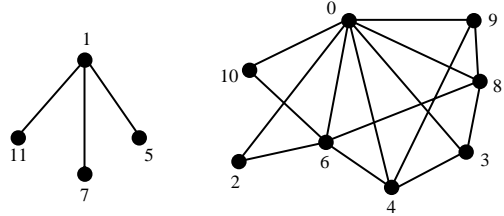

Figure 3.1.24

This graph is the associated with the following adjacency matrix

$$A = \begin{array}{c} \\ 0 \\ 1 \\ 2 \\ 3 \\ 4 \\ 5 \\ 6 \\ 7 \\ 8 \\ 9 \\ 10 \\ 11 \end{array} \begin{array}{c} \begin{array}{cccccccccccc} 0 & 1 & 2 & 3 & 4 & 5 & 6 & 7 & 8 & 9 & 10 & 11 \end{array} \\ \left[ \begin{array}{cccccccccccc} 0 & 0 & 1 & 1 & 1 & 0 & 1 & 0 & 1 & 1 & 1 & 0 \\ 0 & 0 & 0 & 0 & 0 & 1 & 0 & 1 & 0 & 0 & 0 & 1 \\ 1 & 0 & 0 & 0 & 0 & 0 & 1 & 0 & 0 & 0 & 0 & 0 \\ 1 & 0 & 0 & 0 & 1 & 0 & 0 & 0 & 1 & 0 & 0 & 0 \\ 1 & 0 & 0 & 1 & 0 & 0 & 1 & 0 & 0 & 1 & 0 & 0 \\ 0 & 1 & 0 & 0 & 0 & 0 & 0 & 0 & 0 & 0 & 0 & 0 \\ 1 & 0 & 1 & 0 & 1 & 0 & 0 & 0 & 1 & 0 & 1 & 0 \\ 0 & 1 & 0 & 0 & 0 & 0 & 0 & 0 & 0 & 0 & 0 & 0 \\ 1 & 0 & 0 & 1 & 0 & 0 & 1 & 0 & 0 & 1 & 0 & 0 \\ 1 & 0 & 0 & 0 & 1 & 0 & 0 & 0 & 1 & 0 & 0 & 0 \\ 1 & 0 & 0 & 0 & 0 & 0 & 1 & 0 & 0 & 0 & 0 & 0 \\ 0 & 1 & 0 & 0 & 0 & 0 & 0 & 0 & 0 & 0 & 0 & 0 \end{array} \right] \end{array}$$

The zero graph associated with $Z_{12}$ is given in the following.

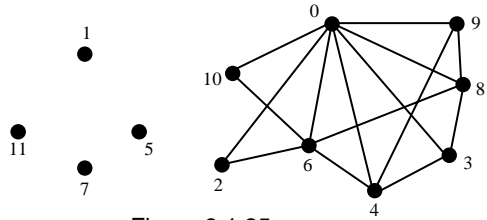

Figure 3.1.25

The adjacency matrix of the zero graph is as follows.

$$B = \begin{array}{c|cccccccccccc} & 0 & 1 & 2 & 3 & 4 & 5 & 6 & 7 & 8 & 9 & 10 & 11 \\ \hline 0 & 0 & 0 & 1 & 1 & 1 & 0 & 1 & 0 & 1 & 1 & 1 & 0 \\ 1 & 0 & 0 & 0 & 0 & 0 & 0 & 0 & 0 & 0 & 0 & 0 & 0 \\ 2 & 1 & 0 & 0 & 0 & 0 & 0 & 1 & 0 & 0 & 0 & 0 & 0 \\ 3 & 1 & 0 & 0 & 0 & 1 & 0 & 0 & 0 & 1 & 0 & 0 & 0 \\ 4 & 1 & 0 & 0 & 1 & 0 & 0 & 1 & 0 & 0 & 1 & 0 & 0 \\ 5 & 0 & 0 & 0 & 0 & 0 & 0 & 0 & 0 & 0 & 0 & 0 & 0 \\ 6 & 1 & 0 & 1 & 0 & 1 & 0 & 0 & 0 & 1 & 0 & 1 & 0 \\ 7 & 0 & 0 & 0 & 0 & 0 & 0 & 0 & 0 & 0 & 0 & 0 & 0 \\ 8 & 1 & 0 & 0 & 1 & 0 & 0 & 1 & 0 & 0 & 1 & 0 & 0 \\ 9 & 1 & 0 & 0 & 0 & 1 & 0 & 0 & 0 & 1 & 0 & 0 & 0 \\ 10 & 1 & 0 & 0 & 0 & 0 & 0 & 1 & 0 & 0 & 0 & 0 & 0 \\ 11 & 0 & 0 & 0 & 0 & 0 & 0 & 0 & 0 & 0 & 0 & 0 & 0 \end{array}$$

The identity graph of the semigroup $Z_{12}$ is as follows.

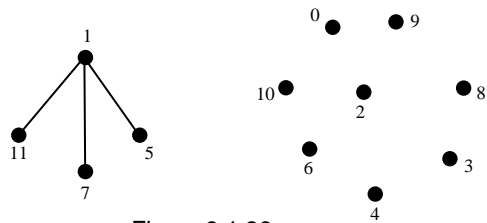

Figure 3.1.26

The adjacency matrix of the identity graph is as follows.

$$C = \begin{matrix} & \begin{matrix} 0 & 1 & 2 & 3 & 4 & 5 & 6 & 7 & 8 & 9 & 10 & 11 \end{matrix} \\ \begin{matrix} 0 \\ 1 \\ 2 \\ 3 \\ 4 \\ 5 \\ 6 \\ 7 \\ 8 \\ 9 \\ 10 \\ 11 \end{matrix} & \begin{bmatrix} 0 & 0 & 0 & 0 & 0 & 0 & 0 & 0 & 0 & 0 & 0 & 0 \\ 0 & 0 & 0 & 0 & 0 & 1 & 0 & 1 & 0 & 0 & 0 & 1 \\ 0 & 0 & 0 & 0 & 0 & 0 & 0 & 0 & 0 & 0 & 0 & 0 \\ 0 & 0 & 0 & 0 & 0 & 0 & 0 & 0 & 0 & 0 & 0 & 0 \\ 0 & 0 & 0 & 0 & 0 & 0 & 0 & 0 & 0 & 0 & 0 & 0 \\ 0 & 1 & 0 & 0 & 0 & 0 & 0 & 0 & 0 & 0 & 0 & 0 \\ 0 & 0 & 0 & 0 & 0 & 0 & 0 & 0 & 0 & 0 & 0 & 0 \\ 0 & 1 & 0 & 0 & 0 & 0 & 0 & 0 & 0 & 0 & 0 & 0 \\ 0 & 0 & 0 & 0 & 0 & 0 & 0 & 0 & 0 & 0 & 0 & 0 \\ 0 & 0 & 0 & 0 & 0 & 0 & 0 & 0 & 0 & 0 & 0 & 0 \\ 0 & 0 & 0 & 0 & 0 & 0 & 0 & 0 & 0 & 0 & 0 & 0 \\ 0 & 1 & 0 & 0 & 0 & 0 & 0 & 0 & 0 & 0 & 0 & 0 \end{bmatrix} \end{matrix}$$

It is verified that the combined identity-zero adjacency matrix of the combined identity-zero graph can be, got as the sum of the adjacency matrix of the special identity and adjacency matrix of the zero divisor graph, i.e., $A = B + C$.

***Example 3.1.26:*** Let $Z_6 = \{0, 1, 2, 3, 4, 5\}$ be the semigroup under multiplication modulo 6.
The zero divisor graph of $Z_6$ is

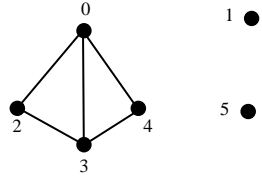

Figure 3.1.27

The associated adjacency matrix D of the zero divisor graph

$$D = \begin{array}{c} \\ 0 \\ 1 \\ 2 \\ 3 \\ 4 \\ 5 \end{array} \begin{array}{c} 0 \ 1 \ 2 \ 3 \ 4 \ 5 \\ \begin{bmatrix} 0 & 0 & 1 & 1 & 1 & 0 \\ 0 & 0 & 0 & 0 & 0 & 0 \\ 1 & 0 & 0 & 1 & 0 & 0 \\ 1 & 0 & 1 & 0 & 1 & 0 \\ 1 & 0 & 0 & 1 & 0 & 0 \\ 0 & 0 & 0 & 0 & 0 & 0 \end{bmatrix} \end{array}$$

The special identity graph of $Z_6$ is

Figure 3.1.28

The adjacency matrix B of the special identity graph is as follows:

$$B = \begin{array}{c} \\ 0 \\ 1 \\ 2 \\ 3 \\ 4 \\ 5 \end{array} \begin{array}{c} 0 \ 1 \ 2 \ 3 \ 4 \ 5 \\ \begin{bmatrix} 0 & 0 & 0 & 0 & 0 & 0 \\ 0 & 0 & 0 & 0 & 0 & 1 \\ 0 & 0 & 0 & 0 & 0 & 0 \\ 0 & 0 & 0 & 0 & 0 & 0 \\ 0 & 0 & 0 & 0 & 0 & 0 \\ 0 & 1 & 0 & 0 & 0 & 0 \end{bmatrix} \end{array}.$$

The combined identity-zero divisor graph of $Z_6$ is given below.

Figure 3.1.29

The combined adjacency matrix C of the identity-zero divisor combined graph C is as follows:

$$C = \begin{array}{c|cccccc} & 0 & 1 & 2 & 3 & 4 & 5 \\ \hline 0 & 0 & 0 & 1 & 1 & 1 & 0 \\ 1 & 0 & 0 & 0 & 0 & 0 & 1 \\ 2 & 1 & 0 & 0 & 1 & 0 & 0 \\ 3 & 1 & 0 & 1 & 0 & 1 & 0 \\ 4 & 1 & 0 & 0 & 1 & 0 & 0 \\ 5 & 0 & 1 & 0 & 0 & 0 & 0 \end{array}$$

We see

$$C = \begin{bmatrix} 0 & 0 & 1 & 1 & 1 & 0 \\ 0 & 0 & 0 & 0 & 0 & 1 \\ 1 & 0 & 0 & 1 & 0 & 0 \\ 1 & 0 & 1 & 0 & 1 & 0 \\ 1 & 0 & 0 & 1 & 0 & 0 \\ 0 & 1 & 0 & 0 & 0 & 0 \end{bmatrix}$$

$$= \begin{bmatrix} 0 & 0 & 1 & 1 & 1 & 0 \\ 0 & 0 & 0 & 0 & 0 & 0 \\ 1 & 0 & 0 & 1 & 0 & 0 \\ 1 & 0 & 1 & 0 & 1 & 0 \\ 1 & 0 & 0 & 1 & 0 & 0 \\ 0 & 0 & 0 & 0 & 0 & 0 \end{bmatrix} + \begin{bmatrix} 0 & 0 & 0 & 0 & 0 & 0 \\ 0 & 0 & 0 & 0 & 0 & 1 \\ 0 & 0 & 0 & 0 & 0 & 0 \\ 0 & 0 & 0 & 0 & 0 & 0 \\ 0 & 0 & 0 & 0 & 0 & 0 \\ 0 & 1 & 0 & 0 & 0 & 0 \end{bmatrix}$$

$$= D + B.$$

***Example 3.1.27:*** Let $Z_8 = \{0, 1, 2, \ldots, 7\}$ be the semigroup under multiplication modulo 8.

The zero divisor graph of $Z_8$ is as follows:

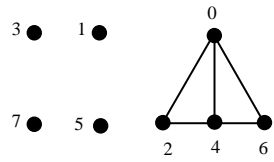

Figure 3.1.30

The corresponding adjacency matrix of the above graph is

$$B = \begin{array}{c} \\ 0 \\ 1 \\ 2 \\ 3 \\ 4 \\ 5 \\ 6 \\ 7 \end{array} \begin{array}{c} \begin{matrix} 0 & 1 & 2 & 3 & 4 & 5 & 6 & 7 \end{matrix} \\ \begin{bmatrix} 0 & 0 & 1 & 0 & 1 & 0 & 1 & 0 \\ 0 & 0 & 0 & 0 & 0 & 0 & 0 & 0 \\ 1 & 0 & 0 & 0 & 1 & 0 & 0 & 0 \\ 0 & 0 & 0 & 0 & 0 & 0 & 0 & 0 \\ 1 & 0 & 1 & 0 & 0 & 0 & 1 & 0 \\ 0 & 0 & 0 & 0 & 0 & 0 & 0 & 0 \\ 1 & 0 & 0 & 0 & 1 & 0 & 0 & 0 \\ 0 & 0 & 0 & 0 & 0 & 0 & 0 & 0 \end{bmatrix} \end{array}$$

The special identity graph of $Z_8$ is as follows:

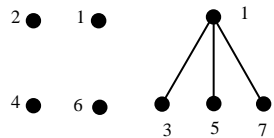

Figure 3.1.31

The adjacency matrix of the above graph is

$$A = \begin{array}{c} \\ 0 \\ 1 \\ 2 \\ 3 \\ 4 \\ 5 \\ 6 \\ 7 \end{array} \begin{array}{c} 0\ 1\ 2\ 3\ 4\ 5\ 6\ 7 \\ \begin{bmatrix} 0 & 0 & 0 & 0 & 0 & 0 & 0 & 0 \\ 0 & 0 & 0 & 1 & 0 & 1 & 0 & 1 \\ 0 & 0 & 0 & 0 & 0 & 0 & 0 & 0 \\ 0 & 1 & 0 & 0 & 0 & 1 & 0 & 0 \\ 0 & 0 & 0 & 0 & 0 & 0 & 0 & 0 \\ 0 & 1 & 0 & 0 & 0 & 0 & 0 & 0 \\ 0 & 0 & 0 & 0 & 0 & 0 & 0 & 0 \\ 0 & 1 & 0 & 0 & 0 & 0 & 0 & 0 \end{bmatrix} \end{array}$$

Now we find

$$B + A = \begin{array}{c} \\ 0 \\ 1 \\ 2 \\ 3 \\ 4 \\ 5 \\ 6 \\ 7 \end{array} \begin{array}{c} 0\ 1\ 2\ 3\ 4\ 5\ 6\ 7 \\ \begin{bmatrix} 0 & 0 & 1 & 0 & 1 & 0 & 1 & 0 \\ 0 & 0 & 0 & 1 & 0 & 1 & 0 & 1 \\ 1 & 0 & 0 & 0 & 1 & 0 & 0 & 0 \\ 0 & 1 & 0 & 0 & 0 & 0 & 0 & 0 \\ 1 & 0 & 1 & 0 & 0 & 0 & 1 & 0 \\ 0 & 1 & 0 & 0 & 0 & 0 & 0 & 0 \\ 1 & 0 & 0 & 0 & 1 & 0 & 0 & 0 \\ 0 & 1 & 0 & 0 & 0 & 0 & 0 & 0 \end{bmatrix} \end{array}.$$

We find the special identity-zero graph of $Z_8$.

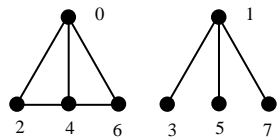

Figure 3.1.32

The adjacency matrix associated with the graph is as follows:

$$\begin{array}{c} \phantom{0} \\ 0 \\ 1 \\ 2 \\ 3 \\ 4 \\ 5 \\ 6 \\ 7 \end{array} \begin{array}{c} 0\ 1\ 2\ 3\ 4\ 5\ 6\ 7 \\ \begin{bmatrix} 0 & 0 & 1 & 0 & 1 & 0 & 1 & 0 \\ 0 & 0 & 0 & 1 & 0 & 1 & 0 & 1 \\ 1 & 0 & 0 & 0 & 1 & 0 & 0 & 0 \\ 0 & 1 & 0 & 0 & 0 & 0 & 0 & 0 \\ 1 & 0 & 1 & 0 & 0 & 0 & 1 & 0 \\ 0 & 1 & 0 & 0 & 0 & 0 & 0 & 0 \\ 1 & 0 & 0 & 0 & 1 & 0 & 0 & 0 \\ 0 & 1 & 0 & 0 & 0 & 0 & 0 & 0 \end{bmatrix} \end{array} = B + A.$$

***Example 3.1.28:*** Let $Z_5 = \{0, 1, 2, \ldots, 4\}$ be the semigroup under multiplication modulo 5.

The zero divisor graph of $Z_5$ is

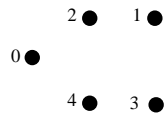

Figure 3.1.33

Thus we see as $Z_5$ has no non trivial zero divisors; no edges in the zero divisor graph but only vertices. Thus the related adjacency matrix would only be a zero matrix as we do not consider $ab = 0$ with $a = 0$ or $b = 0$ as a non trivial zero divisor.

The special identity graph of $Z_5$ is

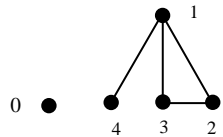

Figure 3.1.34

The related adjacency matrix is as follows.

$$\begin{array}{c c}
 & \begin{array}{c c c c c} 0 & 1 & 2 & 3 & 4 \end{array} \\
\begin{array}{c} 0 \\ 1 \\ 2 \\ 3 \\ 4 \end{array} & \begin{bmatrix} 0 & 0 & 0 & 0 & 0 \\ 0 & 0 & 1 & 1 & 1 \\ 0 & 1 & 0 & 1 & 0 \\ 0 & 1 & 1 & 0 & 0 \\ 0 & 1 & 0 & 0 & 0 \end{bmatrix}
\end{array}.$$

We see the combined identity-zero graph of $Z_5$ is the same as the special identity graph of $Z_5$.

Further the zero divisor graph has only 5 vertices and no edges.

In view of this we have the following theorem.

**THEOREM 3.1.1:** *Let $Z_p = \{0, 1, 2, \ldots, p-1\}$ be the semigroup under multiplication modulo p, p a prime. The zero divisor graph has no edges so the related adjacency matrix is a zero matrix and the special identity graph is the same as the combined identity-zero graph.*

*Proof:* Since in the semigroup $Z_p = \{0, 1, 2, \ldots, p - 1\}$, p a prime, we see $Z_p \setminus \{0\}$ is a group so $Z_p$ has no non trivial zero divisors. Hence the zero divisor graph has no edges hence the associated adjacency matrix is a p × p zero matrix.

Now $Z_p \setminus \{0\}$ is a group so every element $x \in Z_p \setminus \{0\}$ has inverse so in the special identity graph of $Z_p$ all elements in $Z_p \setminus \{0\}$ are adjacent with one and 0 alone is left with no element adjacent with it. Thus we see the matrix associated with the zero divisor graph is just a zero matrix.

Like wise if the semigroup has no identity then this semigroup will not have the special identity graph associated with it. Thus in this case also we will not have the combined identity-zero matrix associated with it.

We give a few examples.

***Example 3.1.29:*** Let $S = 2Z_{12} = \{0, 2, 4, 6, 8, 10\}$ be the semigroup under multiplication modulo 12. Clearly $1 \notin 2Z_{12}$. Thus the zero divisor graph associated with $2Z_{12}$ is

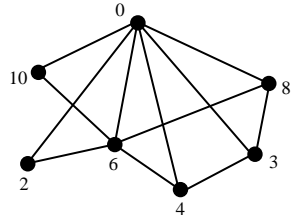

Figure 3.1.35

This has no special identity graph associated with it as $1 \notin 2Z_{12}$.

***Example 3.1.30:*** Let $3Z_{15} = \{0, 3, 6, 12\}$ be the semigroup. The semigroup too has no identity. The zero graph associated with $3Z_{15}$ is as follows.

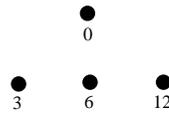

Figure 3.1.36

No zero divisors so no graph can be associated with it.
Also this semigroup cannot have the special identity graph associated with it as $1 \notin 3Z$.

Thus we can have semigroups which has no zero divisors and no special identity graph as it does not contain 1 so such semigroups cannot be given any graph representation. This property is a major difference between the groups and semigroups.

***Example 3.1.31:*** Let $2Z_{30} = \{0, 2, 4, 6, 8, 10, 12, 14, 16, 18, 20, 22, 24, 26, 28\}$ be the semigroup under multiplication modulo 30. Since $1 \notin 2Z_{30}$ the question of its special identity graph does not arise.

Now the zero divisor graph of $2Z_{30}$ is as follows.

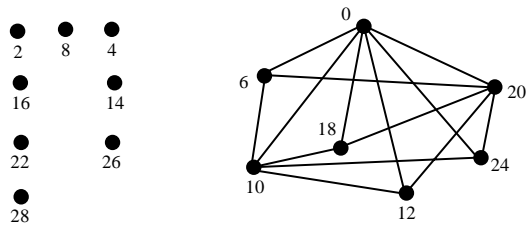
Figure 3.1.37

The associated adjacency matrix M with this graph is as follows.

$$M = \begin{array}{c} \\ 0 \\ 2 \\ 4 \\ 6 \\ 8 \\ 10 \\ 12 \\ 14 \\ 16 \\ 18 \\ 20 \\ 22 \\ 24 \\ 26 \\ 28 \end{array} \begin{array}{c} \begin{array}{cccccccccccccccc} 0 & 2 & 4 & 6 & 8 & 10 & 12 & 14 & 16 & 18 & 20 & 22 & 24 & 26 & 28 \end{array} \\ \left[ \begin{array}{ccccccccccccccc} 0 & 0 & 0 & 1 & 0 & 1 & 1 & 0 & 0 & 1 & 1 & 0 & 1 & 0 & 0 \\ 0 & 0 & 0 & 0 & 0 & 0 & 0 & 0 & 0 & 0 & 0 & 0 & 0 & 0 & 0 \\ 0 & 0 & 0 & 0 & 0 & 0 & 0 & 0 & 0 & 0 & 0 & 0 & 0 & 0 & 0 \\ 1 & 0 & 0 & 0 & 0 & 1 & 0 & 0 & 0 & 0 & 1 & 0 & 0 & 0 & 0 \\ 0 & 0 & 0 & 0 & 0 & 0 & 0 & 0 & 0 & 0 & 0 & 0 & 0 & 0 & 0 \\ 1 & 0 & 0 & 1 & 0 & 0 & 1 & 0 & 0 & 1 & 0 & 0 & 1 & 0 & 0 \\ 1 & 0 & 0 & 0 & 0 & 1 & 0 & 0 & 0 & 0 & 1 & 0 & 0 & 0 & 0 \\ 0 & 0 & 0 & 0 & 0 & 0 & 0 & 0 & 0 & 0 & 0 & 0 & 0 & 0 & 0 \\ 0 & 0 & 0 & 0 & 0 & 0 & 0 & 0 & 0 & 0 & 0 & 0 & 0 & 0 & 0 \\ 1 & 0 & 0 & 0 & 0 & 1 & 0 & 0 & 0 & 0 & 1 & 0 & 0 & 0 & 0 \\ 1 & 0 & 0 & 1 & 0 & 0 & 1 & 0 & 0 & 1 & 0 & 0 & 1 & 0 & 0 \\ 0 & 0 & 0 & 0 & 0 & 0 & 0 & 0 & 0 & 0 & 0 & 0 & 0 & 0 & 0 \\ 1 & 0 & 0 & 0 & 0 & 1 & 0 & 0 & 0 & 0 & 1 & 0 & 0 & 0 & 0 \\ 0 & 0 & 0 & 0 & 0 & 0 & 0 & 0 & 0 & 0 & 0 & 0 & 0 & 0 & 0 \\ 0 & 0 & 0 & 0 & 0 & 0 & 0 & 0 & 0 & 0 & 0 & 0 & 0 & 0 & 0 \end{array} \right] \end{array}$$

This semigroup has no associated special identity graph as $1 \notin 2 Z_{30}$.

***Example 3.1.32:*** Let S(3) be the semigroup of the maps of (1 2 3) to (1 2 3). This semigroup has no zero divisors hence no zero divisor graph associated with it. However the special identity graph associated with it is as follows:

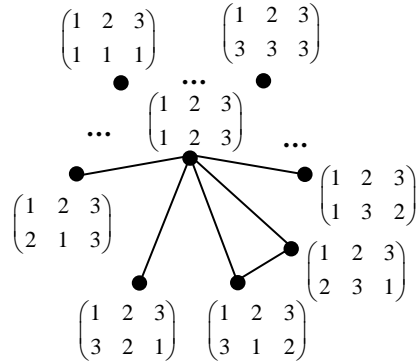

Figure 3.1.38

Thus this semigroup does not have a combined identity-zero divisor graph associated with it.

**THEOREM 3.1.2:** *The class of symmetric semigroups S(n) has only special identity graphs associated with it and with no zero divisor graph hence cannot have the combined identity-zero graph associated with it.*

*Proof:* S(n) is a semigroup of order $n^n$. Clearly S(n) has no zero divisors. So S(n) cannot have any zero divisor graph associated with it.

Further S(n) contains the subset $S_n$ which is the symmetric group of degree n. Now associated with $S_n$ is the special identity graph. Thus with S(n) we have an associated special identity graph and no zero divisor graph. Hence S(n) cannot have the combined identity-zero graph associated with it.

We illustrate this by an example.

***Example 3.1.33:*** Let S(4) be the symmetric semigroup got from the maps of (1 2 3 4) to itself. S(4) is a semigroup under the operation of composition of mappings.

$S_4$ is the permutation group of degree 4 is a proper subset of S(4). Clearly S(4) has no zero divisors. Thus S(4) has no zero divisor graph associated with it. Now $S_4 \subseteq S(4)$ and since $S_4$ is group; every element in $S_4$ has an inverse. Hence $S_4 \subseteq S(4)$ has a special identity graph associated with it.

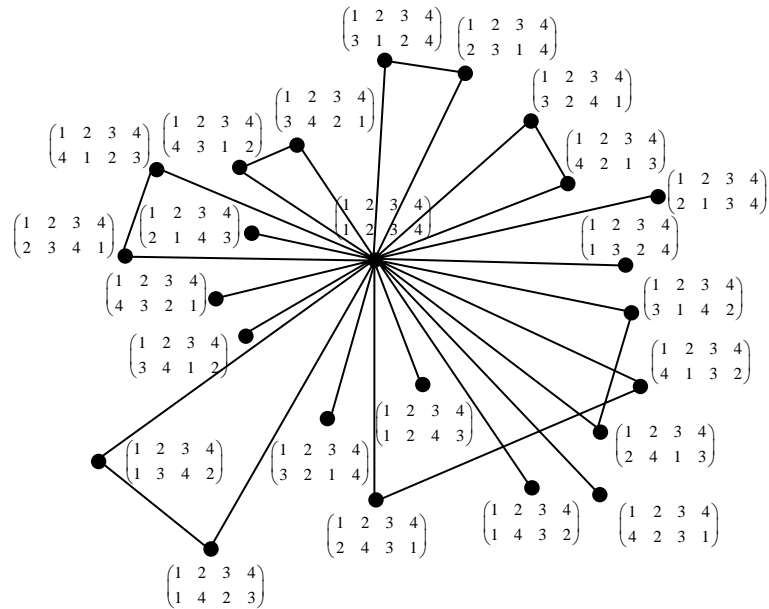

Figure 3.1.39

Thus the semigroup too cannot have a combined identity-zero divisor graph associated with it.

We give a theorem which shows we have a class of semigroups which cannot have the special identity graph associated with it.

**THEOREM 3.1.3:** *Let $p_i Z_n = \{0, p_i, 2p_i, …, p_i(n – 1)\}$ be a semigroup under multiplication modulo n where $p_i / n$ and $n = p_1 p_2 … p_t$, where each $p_i$ is a prime $1 \leq i \leq t$, $t \geq 2$. These* *semigroups do not contain identity so these classes of semigroups have no special identity graphs associated with them.*

*Proof:* Given $p_i Z_n = \{0, p_i, \ldots, p_i(n-1)\}$ is a semigroup under multiplication modulo n, where $n = p_1 p_2 \ldots p_t$, $t \geq 2$ and $p_1 \ldots p_t$ are distinct primes $1 \leq i \leq t$. Clearly $1 \notin p_i Z_n$ so $p_i Z_n$ cannot contain units hence $p_i Z_n$ is a semigroup for which one cannot associate special identity graph with it.

Thus we have a class of semigroups which has no special identity graph associated with it. Thus this class of semigroups cannot have combined identity-zero divisor graph associated with it.

Next we give a class of semigroups which has combined identity-zero divisor graphs associated with it.

**THEOREM 3.1.4:** *Let $Z_n = \{0, 1, 2, \ldots, n-1\}$ be the semigroup under multiplication modulo n where n is a composite number. This semigroup has combined identity-zero divisor graph.*

*Proof:* Given $Z_n = \{0, 1, 2, \ldots, n-1\}$ is a semigroup of order n, n a composite number under multiplication modulo n. Clearly $Z_n$ has zero divisors as well as units. Thus $Z_n$ has a zero divisor graph and a special identity graph associated with it. Hence $Z_n$ has a combined identity-zero divisor graph associated with it. Thus we have a class of semigroups for which we have an associated combined identity-zero graph.

We illustrate this by some examples before we proceed onto define some more new notions.

*Example 3.1.34:* Let $3Z_{24} = \{0, 3, 6, 9, 12, 15, 18, 21\}$ be the semigroup under multiplication modulo 24. We see $3Z_{24}$ has no units but only zero divisors. Also $3Z_{24}$ is not a monoid as $1 \notin 3Z_{24}$. The zero divisor graph associated with $3Z_{24}$ is as follows:

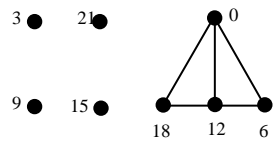

Figure 3.1.40

The adjacency matrix of the zero divisor graph is

$$\begin{array}{c|cccccccc} & 0 & 3 & 6 & 9 & 12 & 15 & 18 & 21 \\ \hline 0 & 0 & 0 & 1 & 0 & 1 & 0 & 1 & 0 \\ 3 & 0 & 0 & 0 & 0 & 0 & 0 & 0 & 0 \\ 6 & 1 & 0 & 0 & 0 & 1 & 0 & 0 & 0 \\ 9 & 0 & 0 & 0 & 0 & 0 & 0 & 0 & 0 \\ 12 & 1 & 0 & 1 & 0 & 0 & 0 & 1 & 0 \\ 15 & 0 & 0 & 0 & 0 & 0 & 0 & 0 & 0 \\ 18 & 1 & 0 & 0 & 0 & 1 & 0 & 0 & 0 \\ 21 & 0 & 0 & 0 & 0 & 0 & 0 & 0 & 0 \end{array}$$

***Example 3.1.35:*** Let $4Z_{24} = \{0, 4, 8, 12, 16, 20\}$ be the semigroup under multiplication modulo 24. This semigroup too has no units but only zero divisors. Infact $1 \notin 4Z_{24}$ so no units, that is cannot have the special identity graph associated with it. The zero divisor graph associated with $4Z_{24}$ is as follows:

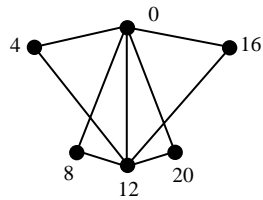

Figure 3.1.41

The zero divisor matrix associated with this graph is as follows;

$$\begin{array}{c|cccccc} & 0 & 4 & 8 & 12 & 16 & 20 \\ \hline 0 & 0 & 1 & 1 & 1 & 1 & 1 \\ 4 & 1 & 0 & 0 & 1 & 0 & 0 \\ 8 & 1 & 0 & 0 & 1 & 0 & 0 \\ 12 & 1 & 1 & 1 & 0 & 1 & 1 \\ 16 & 1 & 0 & 0 & 1 & 0 & 0 \\ 20 & 1 & 0 & 0 & 1 & 0 & 0 \end{array}$$

**Example 3.1.36:** Consider the semigroup $2Z_{24} = \{0, 2, 4, 6, 8, 10, 12, 14, 16, 18, 20, 22\}$ under multiplication modulo 24. This has no unit so cannot have a special identity graph associated with it. The zero divisor graph associated with $2Z_{24}$ is as follows.

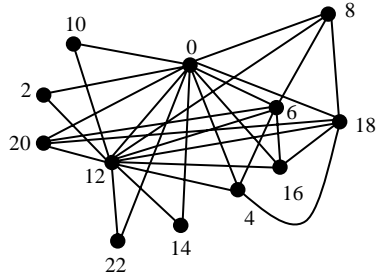

Figure 3.1.42

The adjacency matrix associated with the zero divisor graph is as follows.

$$\begin{array}{c|cccccccccccc} & 0 & 2 & 4 & 6 & 8 & 10 & 12 & 14 & 16 & 18 & 20 & 22 \\ \hline 0 & 0 & 1 & 1 & 1 & 1 & 1 & 1 & 1 & 1 & 1 & 1 & 1 \\ 2 & 1 & 0 & 0 & 0 & 0 & 0 & 1 & 0 & 0 & 0 & 0 & 0 \\ 4 & 1 & 0 & 0 & 1 & 0 & 0 & 1 & 0 & 0 & 1 & 0 & 0 \\ 6 & 1 & 0 & 1 & 0 & 1 & 0 & 1 & 0 & 1 & 0 & 1 & 0 \\ 8 & 1 & 0 & 0 & 1 & 0 & 0 & 1 & 0 & 0 & 1 & 0 & 0 \\ 10 & 1 & 0 & 0 & 0 & 0 & 0 & 1 & 0 & 0 & 0 & 0 & 0 \\ 12 & 1 & 1 & 1 & 1 & 1 & 1 & 0 & 1 & 1 & 1 & 1 & 1 \\ 14 & 1 & 0 & 0 & 0 & 0 & 0 & 1 & 0 & 0 & 0 & 0 & 0 \\ 16 & 1 & 0 & 0 & 1 & 0 & 0 & 1 & 0 & 0 & 1 & 0 & 0 \\ 18 & 1 & 0 & 1 & 0 & 1 & 0 & 1 & 0 & 1 & 0 & 1 & 0 \\ 20 & 1 & 0 & 0 & 1 & 0 & 0 & 1 & 0 & 0 & 1 & 0 & 0 \\ 22 & 1 & 0 & 0 & 0 & 0 & 0 & 1 & 0 & 0 & 0 & 0 & 0 \end{array}$$

**Example 3.1.37:** Let $8Z_{24} = \{0, 8, 16\}$ be the semigroup under multiplication modulo 24. For this semigroup we can take 24 as the identity. This is evident from the following table

| × | 8 | 16 |
|---|---|---|
| 8 | 16 | 8 |
| 16 | 8 | 16 |

Thus the special unit matrix of $8Z_{24}$ is

$$\begin{array}{c}0\\8\\16\end{array}\begin{bmatrix}0 & 0 & 0\\0 & 0 & 1\\0 & 1 & 0\end{bmatrix}$$

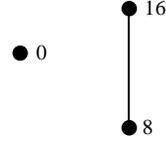

Figure 3.1.43

*Example 3.1.38:* Let $6Z_{24} = \{0, 6, 12, 18\}$ be the semigroup under multiplication modulo 24. The zero divisor graph of $6Z_{24}$

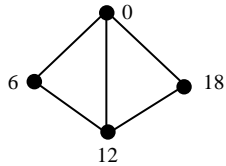

Figure 3.1.44

The matrix of the zero divisor graph is as follows:

$$\begin{array}{c}\\0\\6\\12\\18\end{array}\begin{array}{cccc}0 & 6 & 12 & 18\end{array}\\\begin{bmatrix}0 & 1 & 1 & 1\\1 & 0 & 1 & 0\\1 & 1 & 0 & 1\\1 & 0 & 1 & 0\end{bmatrix}$$

This sort of study with a different element as identity is interesting.

*Example 3.1.39:* Let $Z_{20} = \{0, 1, 2, \ldots, 19\}$ be the semigroup under multiplication modulo 20. The zero divisor graph of $Z_{20}$ is as follows.

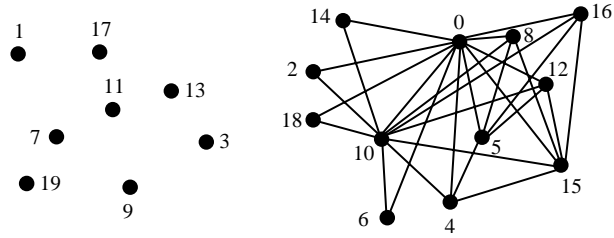
Figure 3.1.45

$$\begin{bmatrix} 0 & 0 & 1 & 0 & 1 & 1 & 1 & 0 & 1 & 0 & 1 & 0 & 1 & 0 & 1 & 1 & 1 & 0 & 1 & 0 \\ 0 & 0 & 0 & 0 & 0 & 0 & 0 & 0 & 0 & 0 & 0 & 0 & 0 & 0 & 0 & 0 & 0 & 0 & 0 & 0 \\ 1 & 0 & 0 & 0 & 0 & 0 & 0 & 0 & 0 & 0 & 1 & 0 & 0 & 0 & 0 & 0 & 0 & 0 & 0 & 0 \\ 0 & 0 & 0 & 0 & 0 & 0 & 0 & 0 & 0 & 0 & 0 & 0 & 0 & 0 & 0 & 0 & 0 & 0 & 0 & 0 \\ 1 & 0 & 0 & 0 & 0 & 1 & 0 & 0 & 0 & 0 & 1 & 0 & 0 & 0 & 0 & 1 & 0 & 0 & 0 & 0 \\ 1 & 0 & 0 & 0 & 1 & 0 & 0 & 0 & 1 & 0 & 0 & 0 & 1 & 0 & 0 & 0 & 1 & 0 & 0 & 0 \\ 1 & 0 & 0 & 0 & 0 & 0 & 0 & 0 & 0 & 0 & 1 & 0 & 0 & 0 & 0 & 0 & 0 & 0 & 0 & 0 \\ 0 & 0 & 0 & 0 & 0 & 0 & 0 & 0 & 0 & 0 & 0 & 0 & 0 & 0 & 0 & 0 & 0 & 0 & 0 & 0 \\ 1 & 0 & 0 & 0 & 0 & 1 & 0 & 0 & 0 & 0 & 1 & 0 & 0 & 0 & 0 & 1 & 0 & 0 & 0 & 0 \\ 0 & 0 & 0 & 0 & 0 & 0 & 0 & 0 & 0 & 0 & 0 & 0 & 0 & 0 & 0 & 0 & 0 & 0 & 0 & 0 \\ 1 & 0 & 1 & 0 & 1 & 0 & 1 & 0 & 1 & 0 & 1 & 0 & 1 & 0 & 1 & 0 & 1 & 0 & 1 & 0 \\ 0 & 0 & 0 & 0 & 0 & 0 & 0 & 0 & 0 & 0 & 0 & 0 & 0 & 0 & 0 & 0 & 0 & 0 & 0 & 0 \\ 1 & 0 & 0 & 0 & 0 & 1 & 0 & 0 & 0 & 0 & 1 & 0 & 0 & 0 & 0 & 1 & 0 & 0 & 0 & 0 \\ 0 & 0 & 0 & 0 & 0 & 0 & 0 & 0 & 0 & 0 & 0 & 0 & 0 & 0 & 0 & 0 & 0 & 0 & 0 & 0 \\ 1 & 0 & 0 & 0 & 0 & 0 & 0 & 0 & 0 & 0 & 1 & 0 & 0 & 0 & 0 & 0 & 0 & 0 & 0 & 0 \\ 1 & 0 & 0 & 0 & 1 & 0 & 0 & 0 & 1 & 0 & 0 & 0 & 1 & 0 & 1 & 0 & 0 & 0 & 0 & 0 \\ 1 & 0 & 0 & 0 & 0 & 1 & 0 & 0 & 0 & 0 & 1 & 0 & 0 & 0 & 0 & 1 & 0 & 0 & 0 & 0 \\ 0 & 0 & 0 & 0 & 0 & 0 & 0 & 0 & 0 & 0 & 0 & 0 & 0 & 0 & 0 & 0 & 0 & 0 & 0 & 0 \\ 1 & 0 & 0 & 0 & 0 & 0 & 0 & 0 & 0 & 0 & 1 & 0 & 0 & 0 & 0 & 0 & 0 & 0 & 0 & 0 \\ 0 & 0 & 0 & 0 & 0 & 0 & 0 & 0 & 0 & 0 & 0 & 0 & 0 & 0 & 0 & 0 & 0 & 0 & 0 & 0 \end{bmatrix}$$

The special identity graph associated with the semigroup is

Figure 3.1.46

The matrix associated with the special identity graph is as follows:

$$\begin{bmatrix}
0 & 0 & 0 & 0 & 0 & 0 & 0 & 0 & 0 & 0 & 0 & 0 & 0 & 0 & 0 & 0 & 0 & 0 & 0 & 0 \\
0 & 0 & 0 & 1 & 0 & 0 & 0 & 1 & 0 & 1 & 0 & 1 & 0 & 1 & 0 & 0 & 0 & 1 & 0 & 1 \\
0 & 0 & 0 & 0 & 0 & 0 & 0 & 0 & 0 & 0 & 0 & 0 & 0 & 0 & 0 & 0 & 0 & 0 & 0 & 0 \\
0 & 1 & 0 & 0 & 0 & 0 & 0 & 1 & 0 & 0 & 0 & 0 & 0 & 0 & 0 & 0 & 0 & 0 & 0 & 0 \\
0 & 0 & 0 & 0 & 0 & 0 & 0 & 0 & 0 & 0 & 0 & 0 & 0 & 0 & 0 & 0 & 0 & 0 & 0 & 0 \\
0 & 0 & 0 & 0 & 0 & 0 & 0 & 0 & 0 & 0 & 0 & 0 & 0 & 0 & 0 & 0 & 0 & 0 & 0 & 0 \\
0 & 0 & 0 & 0 & 0 & 0 & 0 & 0 & 0 & 0 & 0 & 0 & 0 & 0 & 0 & 0 & 0 & 0 & 0 & 0 \\
0 & 0 & 0 & 1 & 0 & 0 & 0 & 0 & 0 & 0 & 0 & 0 & 0 & 0 & 0 & 0 & 0 & 0 & 0 & 0 \\
0 & 0 & 0 & 0 & 0 & 0 & 0 & 0 & 0 & 0 & 0 & 0 & 0 & 0 & 0 & 0 & 0 & 0 & 0 & 0 \\
0 & 1 & 0 & 0 & 0 & 0 & 0 & 0 & 0 & 0 & 0 & 0 & 0 & 0 & 0 & 0 & 0 & 0 & 0 & 0 \\
0 & 0 & 0 & 0 & 0 & 0 & 0 & 0 & 0 & 0 & 0 & 0 & 0 & 0 & 0 & 0 & 0 & 0 & 0 & 0 \\
0 & 1 & 0 & 0 & 0 & 0 & 0 & 0 & 0 & 0 & 0 & 0 & 0 & 0 & 0 & 0 & 0 & 0 & 0 & 0 \\
0 & 0 & 0 & 0 & 0 & 0 & 0 & 0 & 0 & 0 & 0 & 0 & 0 & 0 & 0 & 0 & 0 & 0 & 0 & 0 \\
0 & 1 & 0 & 0 & 0 & 0 & 0 & 0 & 0 & 0 & 0 & 0 & 0 & 0 & 0 & 0 & 1 & 0 & 0 \\
0 & 0 & 0 & 0 & 0 & 0 & 0 & 0 & 0 & 0 & 0 & 0 & 0 & 0 & 0 & 0 & 0 & 0 & 0 & 0 \\
0 & 0 & 0 & 0 & 0 & 0 & 0 & 0 & 0 & 0 & 0 & 0 & 0 & 0 & 0 & 0 & 0 & 0 & 0 & 0 \\
0 & 0 & 0 & 0 & 0 & 0 & 0 & 0 & 0 & 0 & 0 & 0 & 0 & 0 & 0 & 0 & 0 & 0 & 0 & 0 \\
0 & 1 & 0 & 0 & 0 & 0 & 0 & 0 & 0 & 0 & 0 & 0 & 1 & 0 & 0 & 0 & 0 & 0 & 0 \\
0 & 0 & 0 & 0 & 0 & 0 & 0 & 0 & 0 & 0 & 0 & 0 & 0 & 0 & 0 & 0 & 0 & 0 & 0 & 0 \\
0 & 1 & 0 & 0 & 0 & 0 & 0 & 0 & 0 & 0 & 0 & 0 & 0 & 0 & 0 & 0 & 0 & 0 & 0 & 0
\end{bmatrix}$$

One can using both the matrices get the combined special identity-zero divisor graph which is as follows.

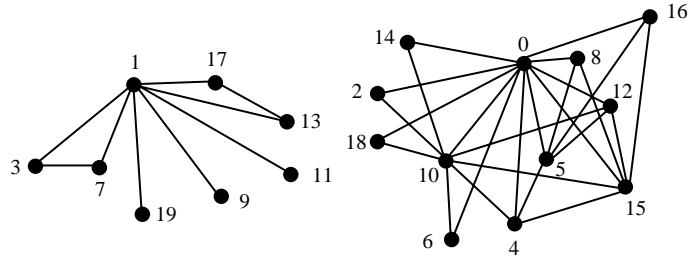

Figure 3.1.47

The related combined adjacency matrix is

$$\begin{bmatrix}
0 & 0 & 1 & 0 & 1 & 1 & 1 & 0 & 1 & 0 & 1 & 0 & 1 & 0 & 1 & 1 & 1 & 0 & 1 & 0 \\
0 & 0 & 0 & 1 & 0 & 0 & 0 & 1 & 0 & 1 & 0 & 1 & 0 & 1 & 0 & 0 & 0 & 1 & 0 & 1 \\
1 & 0 & 0 & 0 & 0 & 0 & 0 & 0 & 0 & 0 & 1 & 0 & 0 & 0 & 0 & 0 & 0 & 0 & 0 & 0 \\
0 & 1 & 0 & 0 & 0 & 0 & 0 & 1 & 0 & 0 & 0 & 0 & 0 & 0 & 0 & 0 & 0 & 0 & 0 & 0 \\
1 & 0 & 0 & 0 & 0 & 1 & 0 & 0 & 0 & 0 & 1 & 0 & 0 & 0 & 0 & 1 & 0 & 0 & 0 & 0 \\
1 & 0 & 0 & 0 & 1 & 0 & 0 & 0 & 1 & 0 & 0 & 0 & 1 & 0 & 0 & 0 & 1 & 0 & 0 & 0 \\
1 & 0 & 0 & 0 & 0 & 0 & 0 & 0 & 0 & 0 & 1 & 0 & 0 & 0 & 0 & 0 & 0 & 0 & 0 & 0 \\
0 & 1 & 0 & 1 & 0 & 0 & 0 & 0 & 0 & 0 & 0 & 0 & 0 & 0 & 0 & 0 & 0 & 0 & 0 & 0 \\
1 & 0 & 0 & 0 & 0 & 1 & 0 & 0 & 0 & 0 & 1 & 0 & 0 & 0 & 0 & 1 & 0 & 0 & 0 & 0 \\
0 & 1 & 0 & 0 & 0 & 0 & 0 & 0 & 0 & 0 & 0 & 0 & 0 & 0 & 0 & 0 & 0 & 0 & 0 & 0 \\
1 & 0 & 1 & 0 & 1 & 0 & 1 & 0 & 1 & 0 & 0 & 0 & 1 & 0 & 1 & 0 & 1 & 0 & 1 & 0 \\
0 & 1 & 0 & 0 & 0 & 0 & 0 & 0 & 0 & 0 & 0 & 0 & 0 & 0 & 0 & 0 & 0 & 0 & 0 & 0 \\
1 & 0 & 0 & 0 & 0 & 1 & 0 & 0 & 0 & 0 & 1 & 0 & 0 & 0 & 0 & 1 & 0 & 0 & 0 & 0 \\
0 & 1 & 0 & 0 & 0 & 0 & 0 & 0 & 0 & 0 & 0 & 0 & 0 & 0 & 0 & 0 & 1 & 0 & 0 & 0 \\
1 & 0 & 0 & 0 & 0 & 0 & 0 & 0 & 0 & 0 & 1 & 0 & 0 & 0 & 0 & 0 & 0 & 0 & 0 & 0 \\
1 & 0 & 0 & 0 & 1 & 0 & 0 & 0 & 1 & 0 & 0 & 0 & 1 & 0 & 0 & 1 & 0 & 0 & 0 & 0 \\
1 & 0 & 0 & 0 & 0 & 1 & 0 & 0 & 0 & 0 & 1 & 0 & 0 & 0 & 1 & 0 & 0 & 0 & 0 & 0 \\
0 & 1 & 0 & 0 & 0 & 0 & 0 & 0 & 0 & 0 & 0 & 0 & 1 & 0 & 0 & 0 & 0 & 0 & 0 & 0 \\
1 & 0 & 0 & 0 & 0 & 0 & 0 & 0 & 0 & 1 & 0 & 0 & 0 & 0 & 0 & 0 & 0 & 0 & 0 & 0 \\
0 & 1 & 0 & 0 & 0 & 0 & 0 & 0 & 0 & 0 & 0 & 0 & 0 & 0 & 0 & 0 & 0 & 0 & 0 & 0
\end{bmatrix}$$

Now we proceed onto define the zero divisor graph of a Smarandache semigroup (S-semigroup). Since every S-semigroup is also a semigroup we have the same definition of zero divisor graph and special identity graph to hold good. However we have the following to be true in case of S-semigroups.

We define for S-semigroups the special group semigroup identity graphs.

**DEFINITION 3.1.2:** *Let S be a S-semigroup. Let P be a proper subset of G such that P is a group under the operations of G. Then we have a special identity graph associated with P. This graph will be known as the special group semigroup identity graph of S.*

*Note:* A S-semigroup has atleast one special group semigroup identity graph. It is pertinent to note in general a semigroup S need not have a special group – semigroup identity graph. We first give some examples of these structures, before we proceed on to define more properties about them.

*Example 3.1.40:* Let $S = \{0, 1, 2, \ldots, 14\}$ be a semigroup under multiplication modulo 15. $P = \{1, 14\}$ is a proper subgroup in S. The special group - semigroup identity graph of S is given by

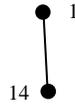

Figure 3.1.48

Take $P_1 = \{5, 10\}$ is again a subgroup in $P_1$ with 10 as the identity. The special group semigroup identity graph associated with it is;

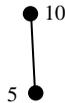

Take the subset $P_2 = \{3, 6, 9, 12\}$ in S. Clearly $P_2$ is a subgroup of S. The special group semigroup identity graph associated with $P_2$ as a group is as follows. For $P_2$, 6 acts as the identity element of the group.

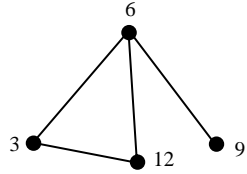

Thus we have seen 3 subgroups in S and their related special group semigroup identity graphs. Here also two groups are isomorphic. Now we see the position of these groups in the combined special identity-zero divisor graph of S.

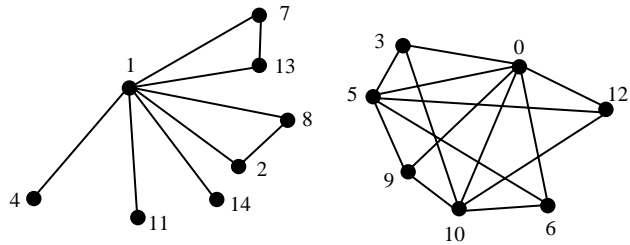

Figure 3.1.49

From this example it is surprising to see the special group semigroup identity graphs are from the zero divisor group also. However one group is from the special identity graph of S.

***Example 3.1.41:*** Let $S = \{0, 1, 2, \ldots, 7\}$ be a semigroup under multiplication modulo 8. Clearly S is a S-semigroup for $7^2 = 1$ (mod 8) and $P = \{1, 7\}$ forms a group. Also $P_1 = \{1, 3\}$ forms a group and $P_2 = \{1, 5\}$ forms a group.
Thus the special group semigroup identity graphs of S related to P, $P_1$ and $P_2$ is as follows:

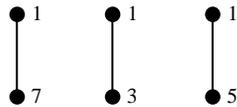

Figure 3.1.50

The combined special identity-zero divisor graph of S is as follows:

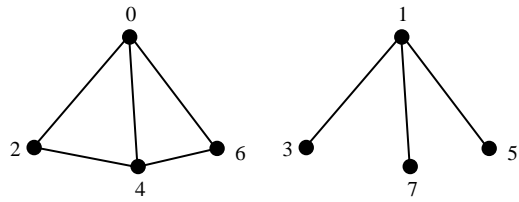

Figure 3.1.51

We see also $P_3 = \{1, 3, 5, 7\}$ is a subgroup of S and its special group semigroup identity graph is the special identity graph of S. However no group has been found from the zero divisor graph of S.

It is pertinent to mention here that when $S = Z_{15}$ and $S = Z_8$ the behavior of these two semigroups under modulo multiplication behaves differently.

***Example 3.1.42:*** Let $S = Z_9 = \{0, 1, 2, \ldots, 8\}$ be a semigroup under multiplication modulo 9. Clearly $Z_9$ is a S-semigroup.

The combined special identity-zero divisor of $Z_9$ is as follows:

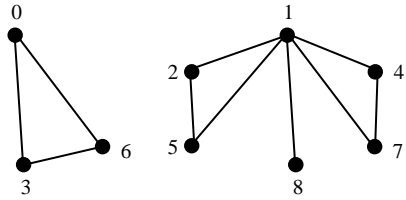

Figure 3.1.52

$P_1 = \{1, 8\}$ is a subgroup of S.
$P_2 = \{1, 7, 4\}$ is again a subgroup of S.
$P_3 = \{1, 2, 5, 8, 7, 4\}$ is again a subgroup of S.

Thus vertices of the special identity graph of S is again a group. However the elements of S which forms the zero divisor graph does not yield to any subgroups of the semigroup.

***Example 3.1.43:*** Let $S = Z_{25} = \{0, 1, 2, \ldots, 24\}$ be a semigroup under multiplication modulo 25.

The special combined identity-zero divisor graph of $Z_{25}$ is as follows:

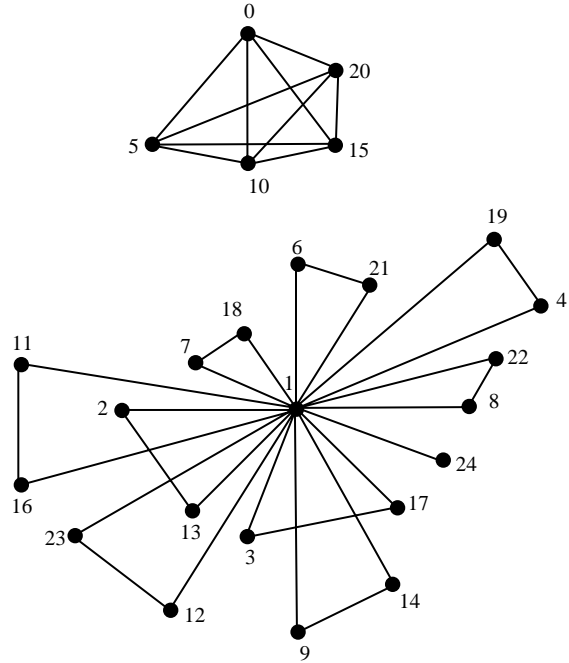

Figure 3.1.53

Now we find the subsets of S which are subgroups under multiplication modulo 25.

$P_1 = \{1, 24\}$ is a subgroup of S.
$P_2 = \{1, 2, 3, 4, 6, 7, 8, 9, 10, 11, 12, 13, 14, 15, 16, 17, 18, 19, 20, 21, 22, 23, 24\}$ is a subgroup.

The special group semigroup identity graph of $S = Z_{25}$ is the whole of the special identity graph given by the following diagram.

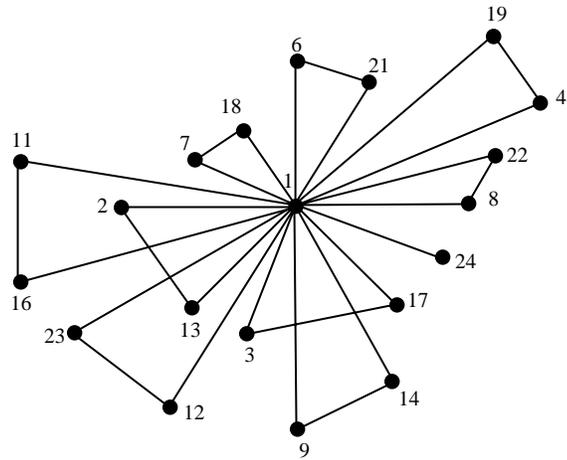

Now we see yet another example.

***Example 3.1.44:*** Let G = {0, 1, 2, 3} be the S-semigroup under multiplication modulo 4. Clearly the combined special identity-zero divisor graph of G is given by

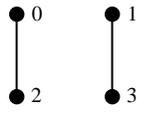

Figure 3.1.54

The special group semigroup identity graph is given by

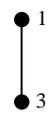

***Example 3.1.45:*** Let S = {0, 1, 2, 3, 4, 5} be the semigroup under multiplication modulo 6. The combined special identity-zero divisor graph of S is

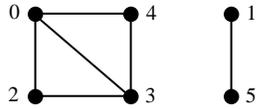
Figure 3.1.55

The special group semigroup identity graph of S is

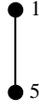

***Example 3.1.46:*** Let $Z_{20} = \{0, 1, 2, \ldots, 19\}$ be the S-semigroup under multiplication modulo 20. The combined special identity-zero divisor graph of $Z_{20}$ is as follows.

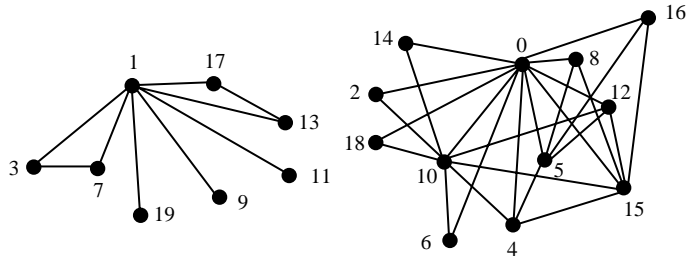
Figure 3.1.56

The subgroup of $Z_{20} = S$ are $P_1 = \{1, 9\}$, $P_2 = \{1, 11\}$, $P_3 = \{1, 19\}$, $P_4 = \{1, 3, 7, 9\}$, $P_5 = \{1, 13, 17, 9\}$, $P_6 = \{1, 3, 7, 9, 11, 13, 17, 19\}$ and $P_7 = \{1, 9, 11, 19\}$.

Thus the related graphs of these subgroups are given in the following diagrams.

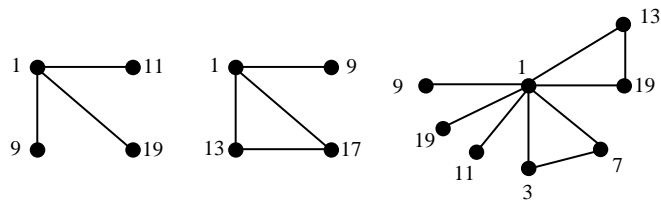

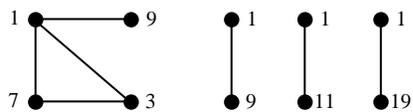

Consider another example.

***Example 3.1.47:*** Let $Z_{30} = \{0, 1, 2, \ldots, 29\}$ be the semigroup under multiplication modulo 30.

The zero divisor graph of $Z_{30}$ is as follows

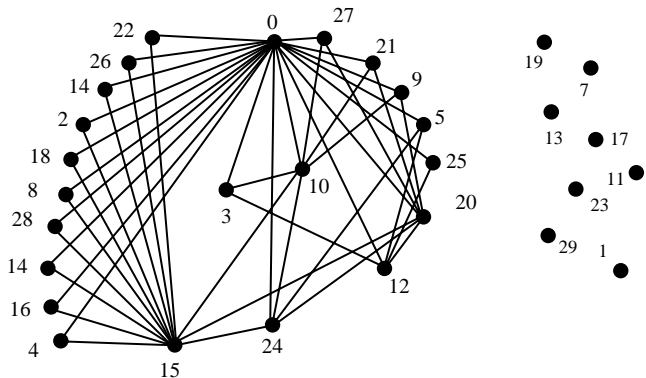

Figure 3.1.57

The special identity graph of $Z_{30}$ is

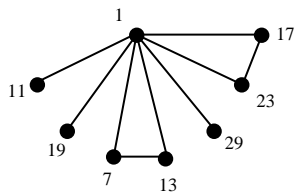

Cleary $Z_{30}$ is also a S-semigroup.

**THEOREM 3.1.5:** *Let S be a S-semigroup. Then S has atleast one nontrivial special identity graph.*

*Proof:* Since every S-semigroup S has a proper subset P which is a group under the operations of S we see the special identity graph of P gives the non trivial special identity graph of S. Hence the claim.

We now show we have classes of S-semigroups which satisfy the above theorem.

The class of symmetric semigroups S(n) where n denotes the set $(a_1, \ldots, a_n)$ or $(1, 2, 3, \ldots, n)$ and S (n) is the set of all maps of the set $(1, 2, 3, \ldots, n)$ to itself. Clearly $S_n \subseteq S(n)$ and $S_n$ is the symmetric group got from the one to one maps of $(1, 2, 3, \ldots, n)$.

Thus S(n) is a S-semigroup and $S_n$ gives the special identity graph. Infact S(n) will have several identity graphs depending on the proper subgroups of $S_n$ including $S_n$.

We illustrate this situation by the following example.

*Example 3.1.48:* Let S(4) be the set of all maps of (1 2 3 4) to itself. Clearly S(4) is a S-semigroup as S(4) contains the symmetric group of degree 4 viz. $S_4$. Some of the subgroups of S(4) are as follows:

$$A_4, H_1 = \left\{ \begin{pmatrix} 1 & 2 & 3 & 4 \\ 1 & 2 & 4 & 3 \end{pmatrix}, e = \begin{pmatrix} 1 & 2 & 3 & 4 \\ 1 & 2 & 3 & 4 \end{pmatrix} \right\}$$

$$H_2 = \left\{ \begin{pmatrix} 1 & 2 & 3 & 4 \\ 1 & 2 & 3 & 4 \end{pmatrix} = e, \begin{pmatrix} 1 & 2 & 3 & 4 \\ 2 & 3 & 4 & 1 \end{pmatrix}, \begin{pmatrix} 1 & 2 & 3 & 4 \\ 3 & 4 & 1 & 2 \end{pmatrix}, \begin{pmatrix} 1 & 2 & 3 & 4 \\ 4 & 1 & 2 & 3 \end{pmatrix} \right\}$$

$$H_3 = \left\{ \begin{pmatrix} 1 & 2 & 3 & 4 \\ 1 & 2 & 3 & 4 \end{pmatrix}, \begin{pmatrix} 1 & 2 & 3 & 4 \\ 2 & 1 & 4 & 3 \end{pmatrix}, \begin{pmatrix} 1 & 2 & 3 & 4 \\ 4 & 3 & 2 & 1 \end{pmatrix}, \begin{pmatrix} 1 & 2 & 3 & 4 \\ 3 & 4 & 1 & 2 \end{pmatrix} \right\}$$ and so on.

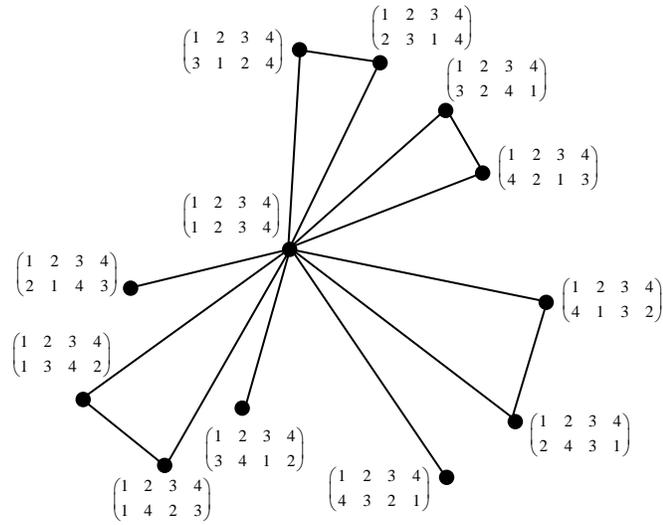

Figure 3.1.58

The special identity graphs associated with the subgroups $A_4$, $H_1$, $H_2$ and $H_3$ are as follows.

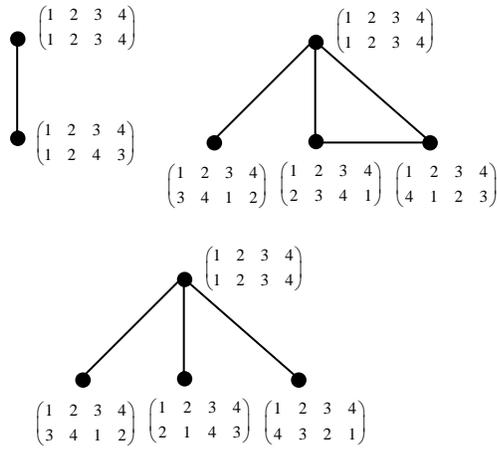

**THEOREM 3.1.6:** *Let $Z_n$ be the semigroup under multiplication modulo n, $n \in N$. $Z_n$ has atleast one special identity graph associated with it.*

*Proof:* The result follows from the fact that every $Z_n$ is a S-semigroup hence every $Z_n$ has atleast one proper subset P which is a group; P being a group one can find the special identity graph associated with it.

***Example 3.1.49:*** Let $3Z_{15} = \{0, 3, 6, 9, 12\}$ be a semigroup under multiplication. Clearly $P = \{3, 6, 9, 12\}$ is a group under multiplication modulo 15. P is given by the following table.

|    | 6  | 3  | 9  | 2  |
|----|----|----|----|----|
| 6  | 6  | 3  | 9  | 12 |
| 3  | 3  | 9  | 12 | 6  |
| 9  | 9  | 12 | 6  | 3  |
| 12 | 12 | 6  | 3  | 9  |

The special identity graph associated with P is as follows.

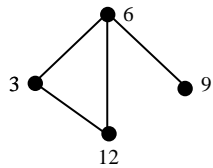

Figure 3.1.59

***Example 3.1.50:*** Let $3Z_{24} = \{0, 3, 6, 9, 12, 15, 18, 21\}$ be the semigroup under multiplication modulo 24. $P = \{9, 3, 15, 21\}$ is a proper subset of $3Z_{24}$ and is a group under multiplication modulo 24 with 9 acting as the identity. The special identity graph for this group is as follows.

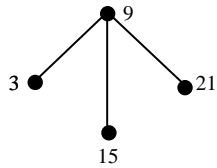

Figure 3.1.60

Now having seen the graphs associated with semigroups we now proceed onto define or extend these notions to loops and commutative groupoids.

## 3.2 Special Identity Graphs of Loops

We define the special identity graph of a loop identical to that of a group as the associativity operation has no role to play. Thus taking verbatim the definition for groups to be true for loops we proceed onto give only examples.

*Example 3.2.1:* Let L be the loop given by the following table.

| * | e | $a_1$ | $a_2$ | $a_3$ | $a_4$ | $a_5$ |
|---|---|---|---|---|---|---|
| e | e | $a_1$ | $a_2$ | $a_3$ | $a_4$ | $a_5$ |
| $a_1$ | $a_1$ | e | $a_4$ | $a_2$ | $a_5$ | $a_3$ |
| $a_2$ | $a_2$ | $a_4$ | e | $a_5$ | $a_3$ | $a_1$ |
| $a_3$ | $a_3$ | $a_2$ | $a_5$ | e | $a_1$ | $a_4$ |
| $a_4$ | $a_4$ | $a_5$ | $a_3$ | $a_1$ | e | $a_2$ |
| $a_5$ | $a_5$ | $a_3$ | $a_1$ | $a_4$ | $a_2$ | e |

The special identity graph related with L is as follows.

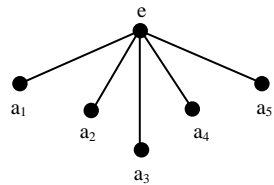

Figure 3.2.1

*Example 3.2.2:* Consider the loop $L_5(2) = \{e, 1, 2, 3, 4, 5\}$. The composition table for $L_5(2)$ is given below

| * | e | 1 | 2 | 3 | 4 | 5 |
|---|---|---|---|---|---|---|
| e | e | 1 | 2 | 3 | 4 | 5 |
| 1 | 1 | e | 3 | 5 | 2 | 4 |
| 2 | 2 | 5 | e | 4 | 1 | 3 |
| 3 | 3 | 4 | 1 | e | 5 | 2 |
| 4 | 4 | 3 | 5 | 2 | e | 1 |
| 5 | 5 | 2 | 4 | 1 | 3 | e |

The special identity graph is as follows.

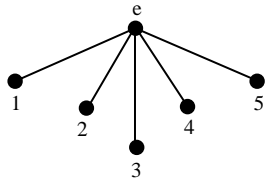

Figure 3.2.2

*Example 3.2.3:* $L_9(8)$ be the loop given below by the following table

| * | e | 1 | 2 | 3 | 4 | 5 | 6 | 7 | 8 | 9 |
|---|---|---|---|---|---|---|---|---|---|---|
| e | e | 1 | 2 | 3 | 4 | 5 | 6 | 7 | 8 | 9 |
| 1 | 1 | e | 9 | 8 | 7 | 6 | 5 | 4 | 3 | 2 |
| 2 | 2 | 3 | e | 1 | 9 | 8 | 7 | 6 | 5 | 4 |
| 3 | 3 | 5 | 4 | e | 2 | 1 | 9 | 8 | 7 | 6 |
| 4 | 4 | 7 | 6 | 5 | e | 3 | 2 | 1 | 9 | 8 |
| 5 | 5 | 9 | 8 | 7 | 6 | e | 4 | 3 | 2 | 1 |
| 6 | 6 | 2 | 1 | 9 | 8 | 7 | e | 5 | 4 | 3 |
| 7 | 7 | 4 | 3 | 2 | 1 | 9 | 8 | e | 6 | 5 |
| 8 | 8 | 6 | 5 | 4 | 3 | 2 | 1 | 9 | e | 7 |
| 9 | 9 | 8 | 7 | 6 | 5 | 4 | 3 | 2 | 1 | e |

The special graph identity of $L_9(8)$ is given below:

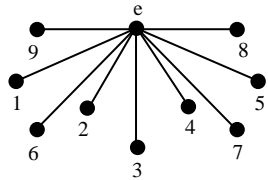

Figure 3.2.3

**THEOREM 3.2.1:** *Let $L_n(m)$ be the loop where $n > 3$, $n$ is odd and $m$ is a positive integer such that $(m, n) = 1$ and $(m – 1, n) =*

*1 with m < n. Then the special identity graph of these loops are rooted trees of special type with (n + 1) vertices.*

*Proof:* We know $L_n(m) = \{e, 1, 2, \ldots, n\}$ is a loop of order n+1 with every element $x \in L_n(m)$ a self inversed element of $L_n(m)$. Thus we see the special identity graph of these loops are special rooted trees with (n + 1) vertices given below.

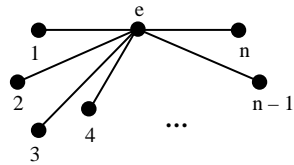

Figure 3.2.4

We can have for the new classes of loops only rooted trees of special form. However for general loop this may not be true.

Now we proceed onto define the special identity graph the zero divisor graph and the combined special identity-zero divisor graph in case of commutative monoids with identity. If the commutative monoids do not contain 1 then we do not have with it the associated special identity graph consequently the notion of combined special identity-zero divisor graph does not exist.

Thus throughout this book we only assume all the groupoids are commutative groupoids.

The notion of special identity graph, zero divisor graph and the combined special identity-zero divisor graph are defined for commutative groupoids as in the case of commutative semigroups.

We illustrate these situations by some examples.

*Example 3.2.4:* Let G be a groupoid given by the following table.

| * | 0 | 1 | 2 | 3 | 4 |
|---|---|---|---|---|---|
| 0 | 0 | 2 | 4 | 1 | 3 |
| 1 | 2 | 4 | 1 | 3 | 0 |
| 2 | 4 | 1 | 3 | 0 | 2 |
| 3 | 1 | 3 | 0 | 2 | 4 |
| 4 | 3 | 0 | 2 | 4 | 1 |

This groupoid contains the elements 0 and 1 but however under the operation * described in the table this groupoid is commutative but

$0 * x = x * 0 = 0$ does not hold good for all $x \in G$.

Also $1 * x = x * 1 = x$ does not hold good for all $x \in G$.

Thus we have class of groupoids which are commutative but for which no graph can be associated. This is explained by the following theorem.

**THEOREM 3.2.2:** *Let $Z_n = \{0, 1, 2, \ldots, n-1\}$ $n \geq 3$; $n < \infty$. Define * on $Z_n$ as $a * b = ta + tb$, $t < n$, $t \in Z_n$. Then $(Z_n, *)$ is a commutative groupoid which has no zero divisors graph or special identity graphs associated with it.*

*Proof:* These groupoids by the very binary operation defined on it are commutative and 0 is such that $0 * x = x * 0 = 0$ does not hold good. For any $a \in Z_n$; $a * 0 = ta + 0t = ta = at \neq 0$ as $a \neq 0$ and $t \neq 0$

Also if $a \in Z_n$ then
$$\begin{aligned} 1 * a &= a * 1 \\ &= ta + t \\ &= t(a+1) \\ &= t + at \neq 1 \end{aligned}$$
as $t \neq 1$ and $a \neq 1$.

So these groupoid do not contain zero or identity. Hence we cannot associate with these groupoids the notion of zero divisor graphs or special identity graphs.

However we can define commutative groupoids with zero divisors and units.

***Example 3.2.5:*** Let G be a groupoid given by the following table:

|   | 0 | 1 | 2 | 3 | 4 | 5 | 6 | 7 |
|---|---|---|---|---|---|---|---|---|
| 0 | 0 | 0 | 0 | 0 | 0 | 0 | 0 | 0 |
| 1 | 0 | 1 | 2 | 3 | 4 | 5 | 6 | 7 |
| 2 | 0 | 2 | 5 | 7 | 1 | 4 | 3 | 3 |
| 3 | 0 | 3 | 7 | 6 | 5 | 4 | 2 | 0 |
| 4 | 0 | 4 | 1 | 5 | 5 | 2 | 2 | 4 |
| 5 | 0 | 5 | 4 | 4 | 2 | 1 | 5 | 6 |
| 6 | 0 | 6 | 3 | 2 | 2 | 5 | 2 | 0 |
| 7 | 0 | 7 | 3 | 0 | 4 | 6 | 0 | 0 |

Clearly G is commutative groupoid. The zero divisor graph associated with G is as follows:

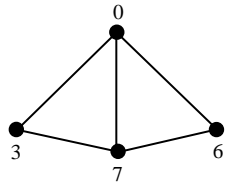

Figure 3.2.5

The special identity graph of G is given below.

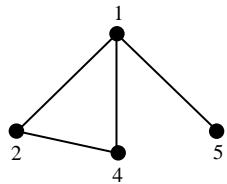

Figure 3.2.6

Here it is pertinent to mention in case of groupoids there are elements which are either zero divisors or units. Interested reader can further study about graphs related with commutative groupoids and loops.

## 3.3 The identity graph of a finite commutative ring with unit

Let R be a finite commutative ring with unit we define here the notion of identity graph of R. The graph of R is formed with elements of R which are units in R. This notion will bring in a nice relation between the graphs and commutative ring with unit. We know already a study relating to zero divisors of a ring with graphs were introduced in 1988 by Beck.

This study was fancied as the vertex coloring of a commutative ring. Beck defined this notion as follows: "A commutative ring with unit is considered as a simple graph R whose vertices are all elements of R such that two different elements x and y in R are adjacent if and only if $x.y = 0$; ($x \neq 0$, $y \neq 0$). The '0' is adjacent with every element in R.

In a similar way we define the new notion of identity graph of a commutative ring with 1 of finite order.

**DEFINITION 3.3.1:** *Let R be a finite commutative ring with 1. We take U(R) the set of units in R (clearly $U(R) \neq \phi$ as $1 \in U(R)$). Now the elements of U(R) form the vertices of the simple graph. Two elements x and y in R are adjacent if and only if $x.y = 1$. We assume that 1 is adjacent with every unit in R. The graph associated with U(R) is defined to be the unit graph of R.*

**Remark:** In case of zero divisor graph we take for the simple graph the vertices as all the elements of R. Here for the identity or unit graph of R we take the vertices as its unit elements of R.

**Remark:** If R has no element other than 1, i.e., U(R) = {1} then the identity or unit graph is just a point.

● 1

***Example 3.3.1:*** The identity or unit graph of $Z_8$ is given below.

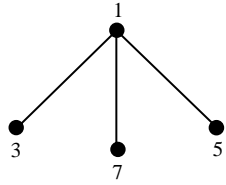

Figure 3.3.1

Now we can compare the unit or identity graph of $Z_8$ with the zero graph of $Z_8$.

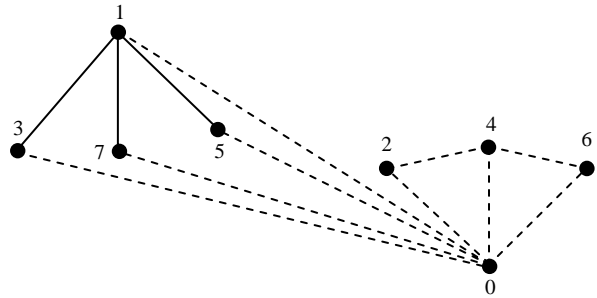

Figure 3.3.2

Thus we see the zero divisor graph alone is

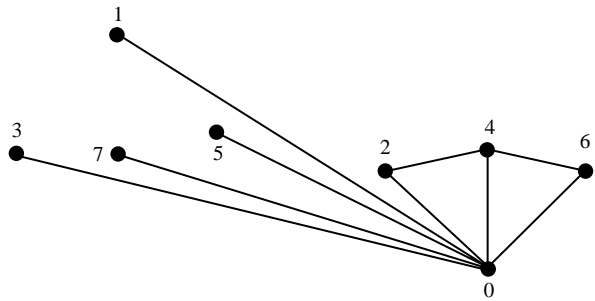

Figure 3.3.3

***Example 3.3.2:*** Let $Z_{12} = \{0, 1, 2, \ldots, 11\}$ be the ring of integers modulo 12. The identity or unit graph associated with $Z_{12}$ is as follows:

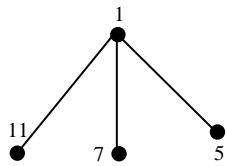

Figure 3.3.4

However the zero divisor graph of $Z_{12}$ is

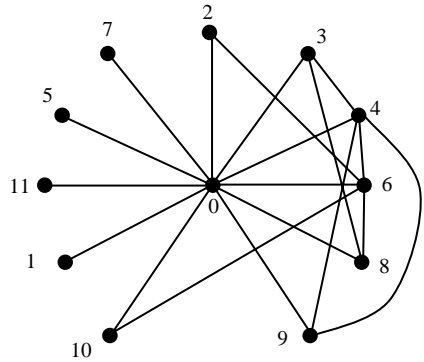

Figure 3.3.5

*Example 3.3.3:* Let $Z_{10} = \{0, 1, 2, \ldots, 9\}$ be the ring of modulo integers 10. The identity graph of $Z_{10}$ is

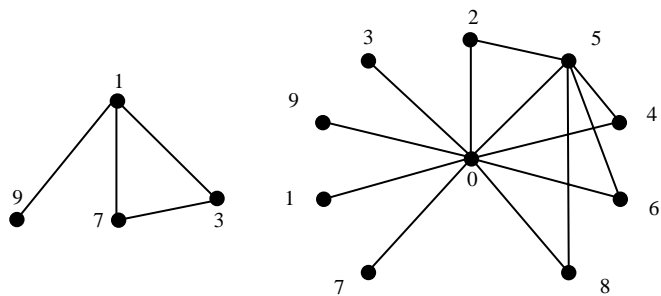

Figure 3.3.6

***Example 3.3.4:*** Let $Z_9 = \{0, 1, 2, \ldots, 8\}$ be the ring of integers modulo 9.
The identity graph of $Z_9$ is

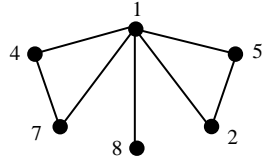

Figure 3.3.7

The zero divisor graph of $Z_9$ is

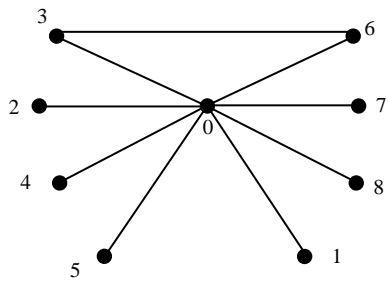

Figure 3.3.8

***Example 3.3.5:*** Let $Z_{15} = \{0, 1, 2, \ldots, 14\}$ be the ring of integers modulo 15.

The identity graph of $Z_{15}$ is

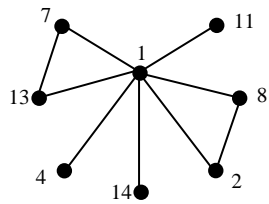

Figure 3.3.9

The zero divisor graph of $Z_{15}$ is given in the following:

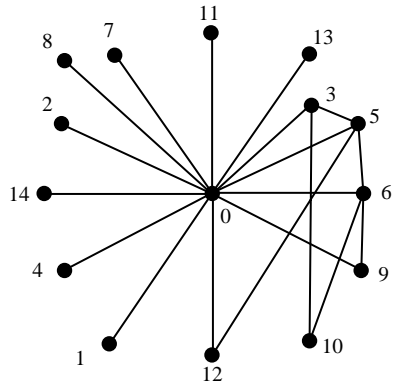

Figure 3.3.10

**Example 3.3.6:** Let $Z_{16} = \{0, 1, 2, \ldots, 15\}$ be the ring of integers modulo 16. The unit graph of $Z_{16}$ is as follows:

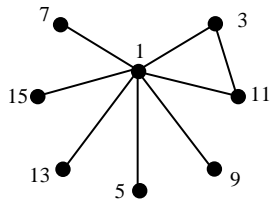

Figure 3.3.11

The zero divisor graph of $Z_{16}$ is as follows:

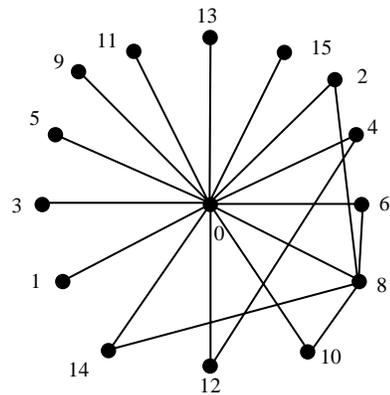

Figure 3.3.12

***Example 3.3.7:*** Let $Z_{25} = \{0, 1, 2, \ldots, 24\}$ be the ring of integers modulo 25. The identity graph of $Z_{25}$ is as follows:

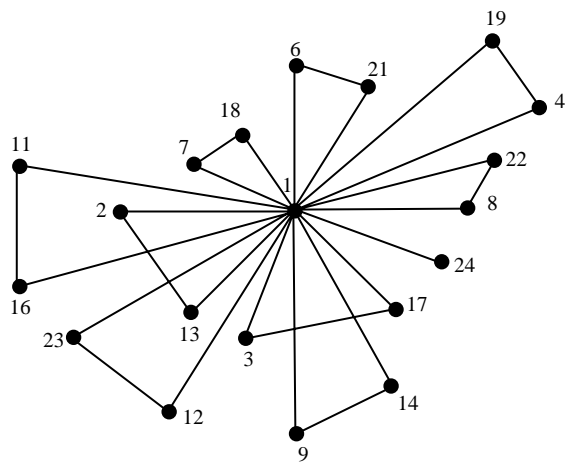

Figure 3.1.13

The unit center is just 1. The zero divisor graph of $Z_{25}$ is

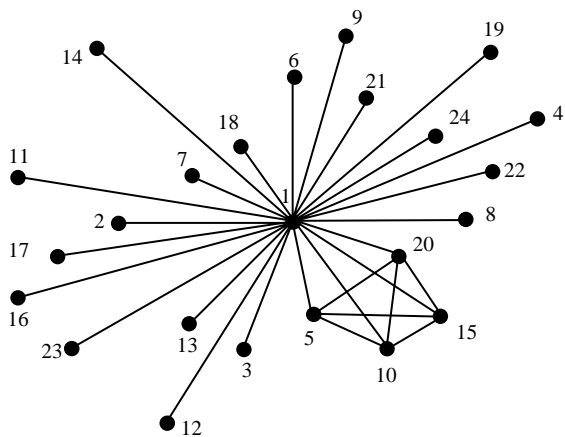

Figure 3.1.14

Now we proceed onto define the notion of combined identity zero divisor graph for a commutative ring with unit.

**DEFINITION 3.3.2:** *Let R be a commutative ring with unit. The combined identity zero divisor graph when it exists is the two graphs namely the zero divisor graph of R and the special identity or unit graph of R.*

We illustrate this situation by some simple examples and also justify the definition.

***Example 3.3.8:*** Let $Z_7 = \{0, 1, 2, \ldots, 6\}$ be the ring under multiplication and addition modulo 7. The zero divisor graph of $Z_7$ does not exist as $x.y = 0 \pmod 7$ is impossible for any $x, y \in Z_7 \setminus \{0\}$.

The special identity or unit graph of $Z_7$ is as follows:

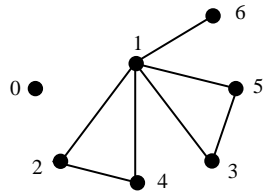

Figure 3.3.15

The unit center is 1.

***Example 3.3.9:*** Let $Z_{11} = \{0, 1, 2, \ldots, 10\}$ be the ring of integers modulo 11. This ring too has no zero divisor graph only this ring has special identity graph associated with it which is as follows:

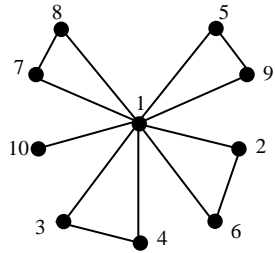

Figure 3.3.16

In view of these examples we can prove the following theorem.

**THEOREM 3.3.1:** *Let $Z_p = \{0, 1, 2, ..., p - 1\}$ be the ring of integers modulo p. $Z_p$ has no zero divisor graph only special identity graph. Hence $Z_p$ has no combined special identity zero divisor graph associated with it.*

*Proof:* We know $Z_p$ is a field hence, $Z_p$ has no nontrivial zero divisors. So $Z_p$ cannot be associated with a zero divisor graph. Since every element in $Z_p \setminus \{0\}$ has inverse, $Z_p$ has a special identity graph with $p - 1$ vertices.

*Example 3.3.10:* Let $Z_{12} = \{0, 1, 2, ..., 11\}$ be the ring of integers modulo 12.

The identity graph of $Z_{12}$ is as follows:

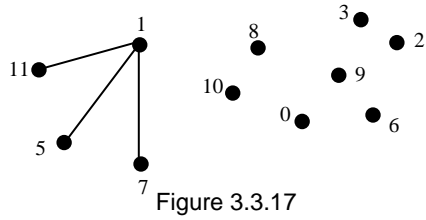

Figure 3.3.17

The unit center is 1.
The zero divisor graph of $Z_{12}$ is

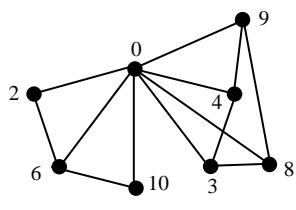

Figure 3.3.18

The zero center of $Z_{12}$ is 0 only.

Thus we see $Z_{12}$ has both zero divisor graph as well as the special identity graph.

Further it is interesting to note that these properly divides $Z_{12}$ into two disjoint classes.

***Example 3.3.11:*** Let $Z_{10} = \{0, 1, 2, \ldots, 9\}$ be the ring of integers modulo 10.

The zero divisor graph of $Z_{10}$ is as follows:

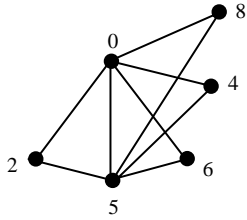

Figure 3.3.19

The zero centers are 0 and 5.
The special identity graph of $Z_{10}$ is as follows:

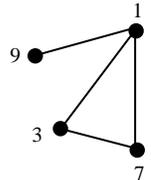

Figure 3.3.20

We see $Z_{10}$ also has both the zero divisor graph and the two graphs are disjoint.
 The unit center of $Z_{10}$ is just one.

***Example 3.3.12:*** Let $Z_{15} = \{0, 1, 2, \ldots, 14\}$ be the ring of integers modulo 15.

The zero divisor graph of $Z_{15}$ is

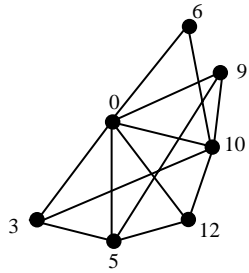
Figure 3.3.21

This graph too has only 0 as its zero center. The special identity graph of $Z_{15}$ is as follows:

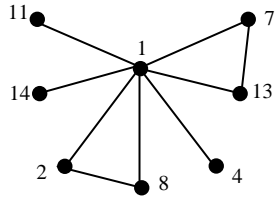
Figure 3.3.22

The unit center of the graph is just 1.

***Example 3.3.13:*** Let $Z_{14} = \{0, 1, 2, \ldots, 13\}$ be the ring of integers modulo 14.

The zero divisor graph of $Z_{14}$ is

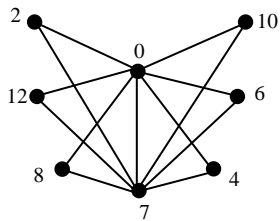
Figure 3.3.23

The zero center of $Z_{14}$ is 0 and 7.
The special identity graph of $Z_{14}$ is as follows:

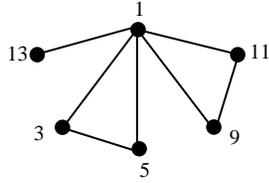
Figure 3.3.24

The unit center of $Z_4$ is just one.

***Example 3.3.14:*** Let $Z_{22} = \{0, 1, 2, \ldots, 21\}$ be the ring of integers modulo 21. The zero divisor graph associated with $Z_{22}$ is

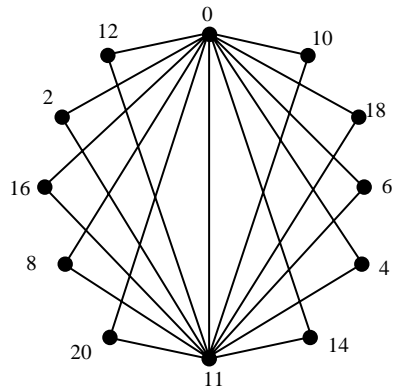
Figure 3.3.25

The zero center of this zero divisor graph is 0 and 11.
  The special identity graph of $Z_{22}$ is as follows:

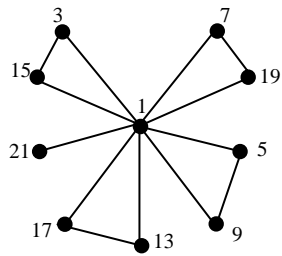
Figure 3.3.26

Clearly the unit center of $Z_{22}$ is 1.

***Example 3.3.15:*** Let $Z_{30} = \{0, 1, 2, \ldots, 29\}$ be the ring of integers modulo 30. The zero divisor graph of $Z_{30}$ is as follows:

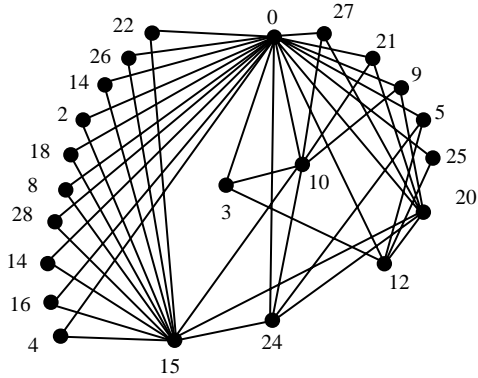

Figure 3.3.27

This graph too has only 0 to be its zero center. Now we draw the unit graph of $Z_{30}$.

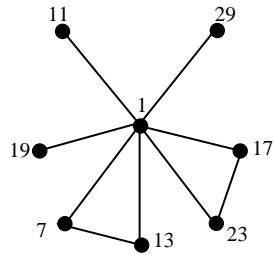

Figure 3.3.28

The unit center of $Z_{30}$ is just 1.

In view of all these examples we have the following theorem.

**THEOREM 3.3.2:** *Let $Z_{2n\,=\,m} = \{0, 1, 2, \ldots, 2n-1 = m-1\}$ be the ring of integers modulo $2n = m$ where n is a prime number. Then $Z_{2n}$ has both the zero divisor graph as well as the special* *identity graph such that the zero center of the zero divisor graph is n and 0 and the unit center of $Z_{2n}$ is just 1.*

*Proof:* We see the zero divisors are contributed by the even numbers and the number of even number is n and they also contribute to zero divisors as n as zero center as well as 0 as zero center. Hence the vertices 0 and p have same number of edges going out of them. The case of unit center is obvious.

**Remark:** We see this is not true when n is a non prime. Also this theorem does not hold good if m = 3p where p is again a prime such (3, p) = 1.

All these claims in the remark are substantiated by examples.

**DEFINITION 3.3.3:** *Let R be a commutative ring or a non commutative ring. The additive inverse graph of R is the special identity graph of R using 0 as the additive identity.*

*Example 3.3.16:* Let $Z_8$ be the ring of integers modulo 8. The additive inverse graph of $Z_8$ is as follows:

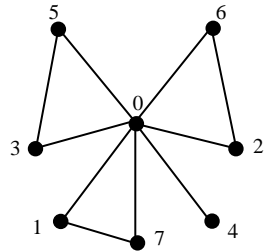

Figure 3.3.29

Clearly the identity (unit) center of $Z_8$ is 0.

*Example 3.3.17:* Let $Z_7$ = {0, 1, 2, …, 6} be the ring of integers modulo 7. The additive inverse graph of $Z_7$ is 0 which is the unit center of $Z_7$.

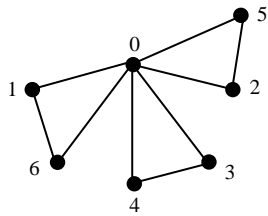

Figure 3.3.30

**Example 3.3.18:** Let $Z_{15} = \{0, 1, 2, \ldots, 14\}$ be the ring of integers modulo 15. The additive inverse graph of $Z_{15}$ is

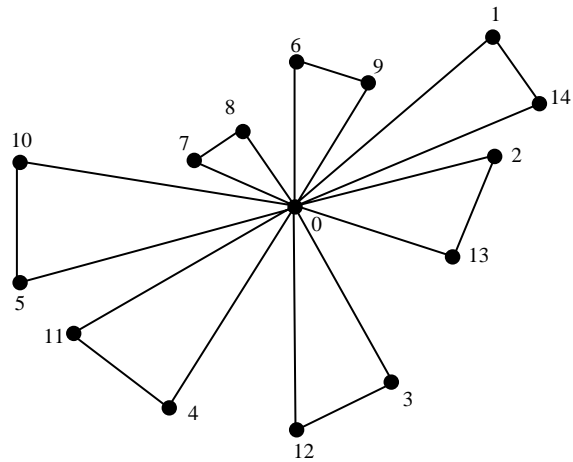

Figure 3.3.31

Clearly 0 is the unit (inverse) center of $Z_{15}$.

Thus we see with a ring we can in general have three graphs associated with them.

(1) additive inverse graphs which always exists and all vertices are included.
(2) zero divisor graph, it may or may not exist. For $Z_p$ has no zero divisors for p a prime.

(3) The special identity graph may or may not exist for if 1 ∉ R then the question of special identity graph becomes superfluous.
(4) Additive inverse graphs always exists for a ring be it commutative or other wise.

We now find the 3 graphs for the following rings.

***Example 3.3.19:*** Let $Z_6 = \{0, 1, 2, 3, 4, 5\}$ be the ring of integers modulo 6.
The zero divisor graph of $Z_6$ is

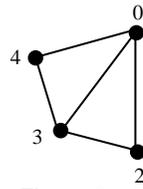

Figure 3.3.32

The special identity graph of $Z_6$ is

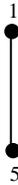

The additive zero graph of $Z_6$ is

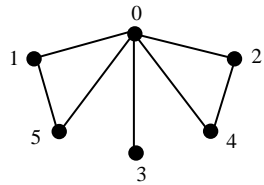

Figure 3.3.33

***Example 3.3.20:*** Let $Z_9 = \{0, 1, 2, \ldots, 8\}$ be the ring of integers modulo 9.
The additive zero graph of $Z_9$ is

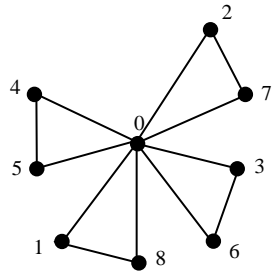

Figure 3.3.34

The zero divisor graph of $Z_9$ is

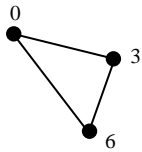

Figure 3.3.35

The special unit graph of $Z_9$ is

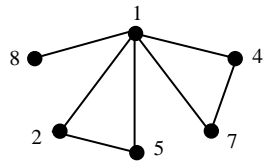

Figure 3.3.36

It is important to note that the zero divisor graph of $Z_9$ happens to be the subgraph of the additive zero graph of $Z_9$.

***Example 3.3.21:*** Let $Z_{18} = \{0, 1, 2, \ldots, 17\}$ be the ring of integers modulo 18.
The zero divisor graph of $Z_{18}$ is as follows:

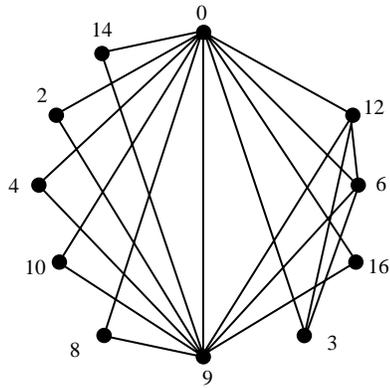

Figure 3.3.37

The additive inverse graph of $Z_{18}$ is as follows:

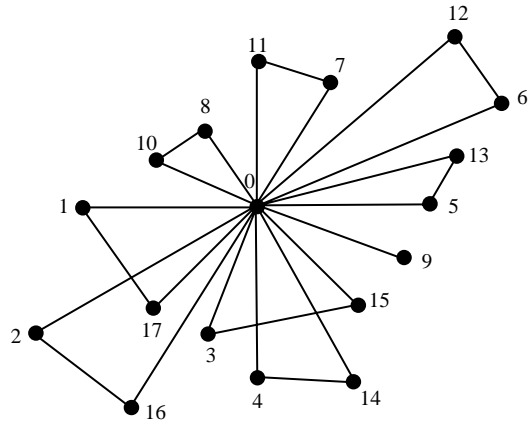

Figure 3.3.38

The special identity graph of $Z_{18}$ is as follows:

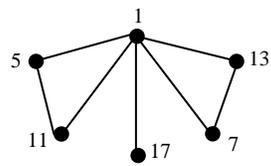

Figure 3.3.39

Thus we see rings $Z_n$ when n is a composite number have all the three graphs associated with it.

Next we proceed onto define the notion of special identity graph and the additive inverse graph of finite fields. It is pertinent to mention at this juncture that fields have no zero divisor graphs associated with them.

The zero divisor graph and the special identity graph of a field are defined as in case of commutative rings. So we illustrate them with examples.

***Example 3.3.22:*** Let $Z_2 = \{0, 1\}$ be the field of characteristic two.
The special identity graph is just a point

The additive graph of $Z_2$ is just a point

***Example 3.3.23:*** Let $Z_3 = \{0, 1, 2\}$ be the prime field of characteristic three.
The additive inverse graph is

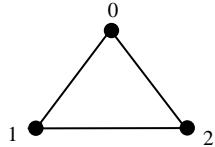

Figure 3.3.40

The special identity graph is

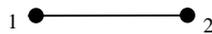

Figure 3.3.41

We see the additive inverse graph has p vertices. The special identity graph has (p – 1) vertices.

**Example 3.3.24:** Let $Z_5 = \{0, 1, 2, \ldots, 4\}$ be the prime field of characteristic five. The additive inverse graph of $Z_5$ is a follows:

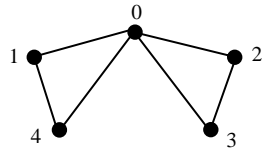

Figure 3.3.42

The special identity graph of $Z_5$ is

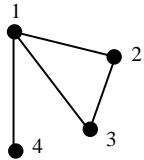

Figure 3.3.43

**Example 3.3.25:** Let $Z_7 = \{0, 1, 2, \ldots, 6\}$ be the prime field of characteristic 7. The identity graph of $Z_7$ is

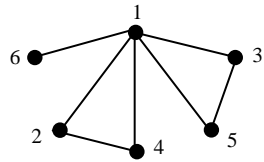

Figure 3.3.44

The additive inverse graph of $Z_7$ is

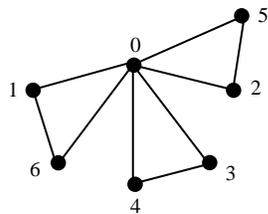

Figure 3.3.45

We see both the graphs in case of fields have only one unit center as 1 in case of special identity graph and 0 in case of additive inverse graph.

Now we proceed onto define the special identity graph and additive inverse graph of a S-ring which we call as the Smarandache special identity graph and Smarandache special additive graph.

**DEFINITION 3.3.4:** *Let R be a S-ring F be a proper subset of R which is the field. The special identity graph of F will be called as the Smarandache special identity graph (S-special identity graph) of R.*

*The additive inverse graph of F will be known as the Smarandache special additive inverse graph (S-special additive inverse graph) of R.*

We first illustrate this situation by some examples.

*Example 3.3.26:* Let $Z_{12} = \{0, 1, 2, \ldots, 11\}$ be the ring of integers modulo 12. The proper subset $F = \{0, 4, 8\} \subseteq Z_{12}$ is a field isomorphic to $Z_3$. 4 acts as the unit element. Thus the Smarandache additive inverse graph of $Z_{12}$ is

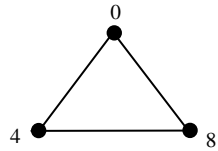

Figure 3.3.46

and the S-special identity graph

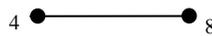

Figure 3.3.47

*Example 3.3.27:* Let $Z_{10} = \{0, 1, 2, \ldots, 9\}$ be the ring of integers modulo 10. $Z_{10}$ is a S-ring. For take $F = \{0, 2, 4, 6, 8\}$ is a field isomorphic to $Z_5$. Here 6 acts as the unit or multiplicative identity.

The S-additive inverse graph of $Z_{10}$ is

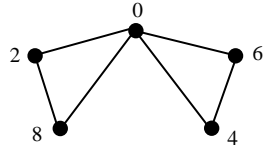

Figure 3.3.48

The S-special identity graph of $Z_{10}$ is

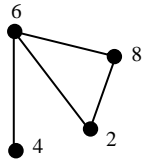

Figure 3.3.49

Now take P = {0, 6} a field isomorphic to $Z_2$. The S- additive inverse graph of $Z_{10}$ is

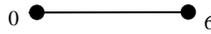

Figure 3.3.50

The S-special identity graph of $Z_{10}$ is

6

Thus we see the ring $Z_{10}$ has two S-special identity graph and two S-additive inverse graph. Thus a S-ring can have in general more than one S-special identity graph and S-additive inverse graph.

*Example 3.3.28:* Let $Z_{30}$ = {0, 1, 2, …, 29} be the ring of integers modulo 30. This is a S-ring. For take $F_1$ = {0, 10, 20} is a field isomorphic to $Z_3$, 10 acts as identity.

This S-additive inverse graph of $Z_{30}$ is

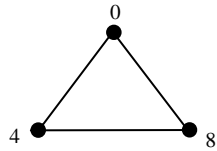

Figure 3.3.51

The S-special identity graph of $Z_{30}$ is $F_2 = \{0, 6, 12, 18, 24\} \subseteq Z_{30}$ is the field isomorphic to $Z_5$ with 6 acting as the multiplicative identity.

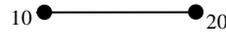

Figure 3.3.52

The S-additive inverse graph of $Z_{30}$ is

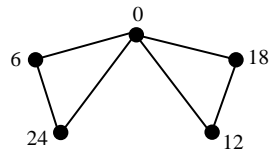

Figure 3.3.53

The S- special identity graph of $Z_{30}$ is

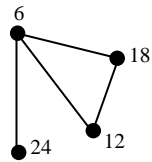

Figure 3.3.54

We have the following theorem.

**THEOREM 3.3.3:** *Let $Z_n = \{0, 1, 2, ..., n – 1\}$ be the ring of integers modulo n. If $n = p_1, p_2, ..., p_t$, t-distinct primes then $Z_n$* *has t-number of S-special identity graphs and S-additive inverse graphs associated with it.*

*Proof:* Given $Z_n$ is the ring of integers modulo n where $n = p_1 p_2 \ldots p_t$, t distinct primes.

Let $m_1 = p_2 \ldots p_t$
$m_2 = p_1 p_3 \ldots p_t$ and so on
$m_t = p_1 p_2 \ldots p_{t-1}$.

Take $m_1 Z_n = \{0, p_1, 2p_1, \ldots, (p_1 - 1)m_1\}$ clearly $m_1 Z_n$ is isomorphic to the field $Z_{p_1}$.

Likewise $m_i Z_n = \{0, p_i, 2p_i, \ldots, (p_i - 1)m_i\}$ is a field isomorphic to $Z_{p_i}$, $1 \leq i \leq t$. Thus relative to each field $Z_{p_1}, \ldots, Z_{p_t}$ we have t number of S-additive inverse graphs and S-special unit identity graphs associated with $Z_n$. Hence the claim.

Chapter Four

# SUGGESTED PROBLEMS

In this chapter we suggest over 50 problems for the reader to solve them.

1. Characterize all groups which are k-colourable normal good groups.

2. Does there exist a group G which is a k-colourable normal good group? (The authors think that there does not exist a group G which is a k-colourable normal good group, i.e., G is a group such that $G = \bigcup_i N_i$; $N_i \subseteq G$ with $N_i \cap N_j = \{1\}$, $N_i$ a normal subgroup of G).

3. Find the special identity graph of $S_4$.

4. Find the special identity graph of $A_5$.

5. Find the special identity graph of $G = \langle g \mid g^{25} = 1 \rangle$.

6. Find the special identity graph of $G = \langle g \mid g^{36} = 1 \rangle$.

7. Prove or disprove isomorphic groups have identical special identity graphs.

8. Are the special identity graphs of $D_{23}$ and $S_3$ same? Justify your claim.

9. Find the special identity graph of $G = A_4 \times S_3$.

10. Find the unit center of the special identity graph $D_{2\,15}$.

11. Find the special identity graph of $D_{2\,16}$.

12. Is it true that the unit centre of the special identity graph of a group is always the vertex which is the identity element of the group G?

13. Can a generalized special identity graph of $S_n$ be given, n any natural number?

14. Find the special identity graph of $S_{25}$.

15. Find the special identity graph of $S_{24}$. (Compare the graphs of $S_{25}$ and $S_{24}$.)

16. Find the special identity graph of $G = A_6 \times A_3$.

17. If $G = G_1 \times G_2$ is the direct product of two groups. What can we say about the graph $G = G_1 \times G_2$? Is it the union of the graphs associated with the two groups? Or is it the sum of the graphs associated with the two groups? Or there exists no relation between the graphs of G and $G_1$ and $G_2$.

18. Let $G = S_3 \times A_4 \times D_{27}$ be the direct product group of $S_3$, $A_4$ and $D_{27}$. Obtain the special identity graph of G. Find the special identity graphs of $A_4$, $S_3$ and $D_{27}$ and find any possible relations between them.

19. Find the special identity graph and the zero graph of the semigroup, $Z_{32}$ under multiplication modulo 32.

20. Find the zero divisor graph of the semigroup $Z_{30} \times Z_{12} = S$.

21. Obtain some interesting results about special identity graphs of the semigroup.

22. Let $G = S_5$ be the symmetric group of degree 5. Find the conjugate graph of $S_5$.

23. Find the conjugate graph of $D_{29}$.

24. How many complete graphs does the conjugate graph of the group $D_{2\,30}$ contain?

25. Can a generalization of the conjugate graph of $S_n$ be made for any n?

26. Find the conjugate graph of the alternating group $A_5$.

27. Obtain all the conjugate graph of the group $A_n$; $n \in N$.

28. Obtain some interesting results about the conjugate graph of the group $S_n$.

29. Find the special identity graph of $S_6$.

30. Find the special identity graph of $S_8$ and compare it with the special identity graph of $S_{27}$.

31. Find the special identity graph of $A_8$ and compare it with the special identity graph of $A_{27}$.

32. Compare the conjugacy graphs of the groups $A_8$ and $A_{27}$.

33. Find the special identity graph of $G = S_3 \times A_4$. Find also the conjugacy graph of G and compare it with the

conjugacy graph of $S_3$ and $A_4$. Does there exist any relation between the 3 conjugacy graphs?

34. Find the zero divisor graph of the group ring $Z_2G$ where $G = \{g \mid g^9 = 1\}$. Find the unit center of the special identity graph of $Z_2G$.

35. Let $Z_3G$ be the group ring of the group $G = \{g \mid g^{25} = 1\}$ over the field $Z_3$.

 (1) Find the special identity graph of $Z_3G$.

 (2) Find the additive inverse graph of $Z_3G$.

 (3) Find the zero divisor graph of $Z_3G$.

36. Let $Z_8G$ where $G = \langle g \mid g^5 = 1 \rangle$ be the group ring of the group G over the ring $Z_8$. Find the zero divisor graph of $Z_8G$. What is the zero center?

37. Let $F = \dfrac{Z_2[x]}{\langle x^5 + x^2 + 1 \rangle = I}$, be the field. Find the additive inverse graph and the unit special identity graph of F; where I is the ideal generated by the polynomial $x^5 + x^2 + 1$, in $Z_2[x]$.

38. Let $Z_{23}$ be the field of characteristic 23. Find the additive inverse graph and the special identity (unit) graph of $Z_{23}$.

39. Characterize those rings which do not contain the zero divisor graph.

40. Characterize those rings which do not contain the special unit or identity graph.

41. Define special identity (unit) graphs for non commutative rings.

42. Find the additive inverse graph of $Z_2S_4$. Find the zero divisor graph and unit graph of $Z_2S_4$.

43. Characterize those rings which has two zero centers for the zero graph.

44. Does these exists rings with more than two zero centers for the zero graph?

45. Can a ring with more than one unit center for the unit graph exist? Justify your claim.

46. Find some interesting properties about the zero centers of the zero graphs of a ring.

47. Apply the zero divisor graphs and unit graphs of a ring to the theory net working in computers.

48. Obtain some interesting applications of these special graphs of groups and rings.

49. Find the S-additive inverse graphs and S-special identity graphs of $Z_{60}$.

50. Find the S-special identity graphs and S-additive inverse graphs of $Z_{30}G$ where $G = \langle g \mid g^{12} = 1 \rangle$.

51. Find the S-special identity graphs and S-additive inverse graphs of $Z_{18}G$ where $G = \langle g \mid g^6 = 1 \rangle$.

52. Find the S-special identity graphs and S-additive inverse graphs of $Z_{20}S_4$.

# FURTHER READING

# INDEX

## T

t-colourable normal bad group, 66
Tree, 7

## U

Unit graph of a ring, 134, 140

## Z

Zero divisors in a groupoid, 12

# ABOUT THE AUTHORS

**Dr.W.B.Vasantha Kandasamy** is an Associate Professor in the Department of Mathematics, Indian Institute of Technology Madras, Chennai. In the past decade she has guided 12 Ph.D. scholars in the different fields of non-associative algebras, algebraic coding theory, transportation theory, fuzzy groups, and applications of fuzzy theory of the problems faced in chemical industries and cement industries.

She has to her credit 646 research papers. She has guided over 68 M.Sc. and M.Tech. projects. She has worked in collaboration projects with the Indian Space Research Organization and with the Tamil Nadu State AIDS Control Society. This is her $42^{nd}$ book.

On India's 60th Independence Day, Dr.Vasantha was conferred the Kalpana Chawla Award for Courage and Daring Enterprise by the State Government of Tamil Nadu in recognition of her sustained fight for social justice in the Indian Institute of Technology (IIT) Madras and for her contribution to mathematics. (The award, instituted in the memory of Indian-American astronaut Kalpana Chawla who died aboard Space Shuttle Columbia). The award carried a cash prize of five lakh rupees (the highest prize-money for any Indian award) and a gold medal.
She can be contacted at vasanthakandasamy@gmail.com
You can visit her on the web at: http://mat.iitm.ac.in/~wbv

---

**Dr. Florentin Smarandache** is a Professor of Mathematics and Chair of Math & Sciences Department at the University of New Mexico in USA. He published over 75 books and 150 articles and notes in mathematics, physics, philosophy, psychology, rebus, literature.

In mathematics his research is in number theory, non-Euclidean geometry, synthetic geometry, algebraic structures, statistics, neutrosophic logic and set (generalizations of fuzzy logic and set respectively), neutrosophic probability (generalization of classical and imprecise probability). Also, small contributions to nuclear and particle physics, information fusion, neutrosophy (a generalization of dialectics), law of sensations and stimuli, etc. He can be contacted at smarand@unm.edu